\documentclass[3p,times]{elsarticle}

\usepackage{array}
\usepackage{color}
\usepackage{tabularx}
\usepackage{graphicx}
\usepackage{amsmath}
\usepackage{amssymb}
\usepackage{amsfonts}	
\usepackage{moreverb}
\usepackage{dsfont}
\usepackage{grffile}
\usepackage{bm}
\usepackage{multirow}
\usepackage{soul}

\newcommand\BibTeX{{\rmfamily B\kern-.05em \textsc{i\kern-.025em b}\kern-.08em
T\kern-.1667em\lower.7ex\hbox{E}\kern-.125emX}}

\newcommand{\x}{\mbf{x}}
\newcommand{\y}{\mbf{y}}

\newcommand{\mbf}[1]{\mathbf{#1}}			%

\newcommand{\Q}{\mathbf{Q}}
\renewcommand{\u}{\mathbf{u}}

\newcommand{\q}{\mathbf{q}}
\newcommand{\F}{\mathbf{F}}
\newcommand{\f}{\mathbf{f}}
\newcommand{\g}{\mathbf{g}}
\newcommand{\h}{\mathbf{h}}
\renewcommand{\v}{\mathbf{v}}

\newcommand{\PNM}{\mathbb{P}_N\mathbb{P}_M}
\newcommand{\PMM}{\mathbb{P}_N\mathbb{P}_N}

\newcommand{\be}{\begin{equation}}
\newcommand{\ee}{\end{equation}}
\newcommand{\bdm}{\begin{displaymath}}
\newcommand{\edm}{\end{displaymath}}

\newcommand{\apriori}{\textit{a priori} }
\newcommand{\aposteriori}{\textit{a posteriori} }

\newfont{\numerikEleven}{ecrm1000}
\newfont{\numerikTen}{cmss10}
\newfont{\numerikNine}{cmss9}
\newfont{\numerikEight}{cmss8}

\journal{Journal of Computational Physics}

\begin{document} 
\begin{frontmatter}
\title{A Posteriori Subcell Limiting of the Discontinuous Galerkin Finite Element Method 
for Hyperbolic Conservation Laws} 
\author[trento]{Michael Dumbser$^{*}$}
\ead{michael.dumbser@unitn.it}
\cortext[cor1]{Corresponding author}

\author[trento]{Olindo Zanotti}
\ead{olindo.zanotti@unitn.it}

\author[imt]{Rapha{\"e}l Loub{\`e}re}
\ead{raphael.loubere@math.univ-toulouse.fr}

\author[lanl]{Steven Diot}
\ead{diot@lanl.gov}

\address[trento]{Department of Civil, Environmental and Mechanical Engineering, University of Trento, Via Mesiano, 77 - 38123 Trento, Italy.}
\address[imt]{CNRS and Institut de Math\'{e}matiques de Toulouse (IMT)
 Universit{\'e} Paul-Sabatier, Toulouse, France}
\address[lanl]{Fluid Dynamics and Solid Mechanics (T-3), Los Alamos National Laboratory, NM 87545 U.S.A.}

\begin{abstract}
The purpose of this work is to propose a novel \textit{a posteriori} finite volume subcell limiter 
technique for the Discontinuous Galerkin finite element method for nonlinear systems of hyperbolic 
conservation laws in multiple space dimensions that works well for \textit{arbitrary} high order 
of accuracy in space \textit{and time} and that does \textit{not} destroy the natural \textit{subcell resolution} 
properties of the DG method. High order time discretization is achieved via a one-step ADER approach that uses 
a local space-time discontinuous Galerkin predictor method to evolve the data locally in time within each cell. 

Our new limiting strategy is based on the so-called MOOD paradigm, which \aposteriori verifies 
the validity of a discrete candidate solution against physical and numerical detection criteria 
after each time step. Here, we employ a relaxed discrete maximum principle in the sense of  
piecewise polynomials and the positivity of the numerical solution as detection criteria. Within the 
DG scheme on the main grid, the discrete solution is represented by piecewise polynomials of 
degree $N$. 
For those troubled cells that need limiting, our new limiter approach recomputes the discrete 
solution by scattering the DG polynomials at the previous time step onto a set of $N_s=2N+1$ finite 
volume subcells per space dimension. 
A robust but accurate ADER-WENO finite volume scheme then updates the subcell averages 
of the conservative variables within the detected troubled cells. The recomputed subcell averages 
are subsequently gathered back into high order cell-centered DG polynomials on the main grid via a subgrid 
reconstruction operator. The choice of $N_s=2N+1$ subcells is optimal since it allows to match the 
maximum admissible time step of the finite volume scheme on the subgrid with the maximum admissible 
time step of the DG scheme on the main grid, minimizing at the same time the local truncation error
of the subcell finite volume scheme. It furthermore provides an excellent subcell resolution of 
discontinuities. 

Our new approach is therefore radically different from classical DG limiters, where the limiter is 
using TVB or (H)WENO \textit{reconstruction} based on the discrete solution of the DG 
scheme on the \textit{main grid} at the new time level. In our case, the discrete solution is \textit{recomputed} 
within the troubled cells using a \textit{different and more robust} numerical scheme on a \textit{subgrid level}. 

We illustrate the performance of the new \aposteriori subcell ADER-WENO finite volume limiter approach for 
very high order DG methods via the simulation of numerous test cases run on Cartesian grids in two and 
three space dimensions, using DG schemes of up to \textit{tenth} order of accuracy in space and time 
($N=9$). The method is also able to run on massively parallel large scale supercomputing infrastructure, 
which is shown via one 3D test problem that uses 10 billion space-time degrees of freedom per time step. 
\end{abstract}

\end{frontmatter}



\section{Introduction} \label{sec:introduction}

The discontinuous Galerkin (DG) finite element method has been originally proposed by Reed and Hill \cite{reed}. Later, a solid theoretical framework has been 
established by Cockburn and Shu in a well-known series of papers \cite{cbs0,cbs1,cbs2,cbs3,cbs4} for the application of discontinuous Galerkin schemes to nonlinear 
hyperbolic systems of conservation laws. A very important property of DG schemes is that they satisfy a local cell entropy inequality for any polynomial degree 
$N$ used for the approximation of the discrete solution. As a consequence, this guarantees \textit{nonlinear} stability in $L_2$ norm for arbitrary high order
of accuracy, see the proof by Jiang and Shu \cite{jiangshu} for the scalar case and its subsequent extensions to systems \cite{BarthCharrier,HouLiu}. This means 
that the DG scheme is by nature very robust and clearly appropriate for the solution of nonlinear hyperbolic conservation laws. However, being a linear scheme
in the sense of Godunov \cite{godunov}, even the DG method needs some sort of nonlinear limiting to avoid the Gibbs phenomenon at shock waves or other
discontinuities. There is a vast literature on the topic of limiters for DG schemes and a non-exhaustive review on this topic will be presented later in section 
\ref{sec:apriori_limiter} of this paper. The key idea of many DG limiters is the following: first an unlimited solution is computed with the DG scheme, then an 
indicator detects so-called \textit{troubled cells}, i.e. those zones of the domain which may need limiting, see \cite{Qiu_2005} for a detailed comparison. 
For troubled cells, the degrees of freedom of the discrete solution are then modified by some sort of nonlinear \textit{reconstruction} technique, based on the 
discrete solution in the troubled cell and its neighbors.  

Concerning the time discretization, mostly explicit TVD Runge-Kutta schemes are used, which lead to the so-called Runge-Kutta DG schemes. A review on DG schemes 
can be found in \cite{CBS-book,CBS-convection-dominated}. However, explicit DG schemes suffer from a very severe time step restriction where the maximum admissible
Courant number typically scales as approximately $1/(2N+1)$, if $N$ denotes the polynomial degree of the approximation of the DG scheme. 
Alternative high order accurate explicit time discretizations for DG schemes have been explored in \cite{QiuDumbserShu} and \cite{taube_jsc,dumbser_jsc}, which, however, 
have led to an even more restrictive CFL condition. 
While the DG method is mostly used only for spatial discretization, it has been introduced as a uniform discretization of space and time in the global space-time 
DG scheme of Van der Vegt et al.  \cite{spacetimedg1,spacetimedg2,KlaijVanDerVegt}, which leads to an implicit method of theoretically arbitrary high order of 
accuracy in space and time and which is unconditionally stable. The strategy presented in \cite{spacetimedg1,spacetimedg2}, however, requires the solution of a 
global nonlinear algebraic system at each time step. In order to reduce the complexity of globally implicit space-time DG schemes, in 
\cite{DumbserEnauxToro,HidalgoDumbser,Dumbser2008} a \textit{local} space-time DG approach has been suggested, which leads only to an element-local implicit method. 
However, also in this case the final DG scheme is explicit and thus has to satisfy the typical stability condition of explicit DG schemes. 


In the finite volume context, recently a new concept has been proposed, namely the Multi-dimensional Optimal Order Detection (MOOD) approach, 
which is an \textit{a posteriori} approach to the problem of limiting. The key idea of this paradigm is to run a spatially \textit{unlimited} high-order finite 
volume scheme in order to produce a so-called \textit{candidate solution}. Then the validity of this candidate solution is tested against a set of 
predefined \textit{admissibility criteria}. Some cells are marked as 'acceptable' and are therefore valid. Some others may be locally marked as 'problematic' 
or 'troubled', if they do not pass the detection process. These cells and their neighbors are consequently \textit{locally recomputed} using polynomial 
reconstructions of a lower degree. Thus, after decrementing the polynomial degree and locally recomputing the solution, a new candidate solution is obtained. 
That solution is again tested for validity and the decrementing procedure re-applies, if necessary. 

Such order decrementing can occur several times within one time step for the same cell, but it will always halt after a finite number of steps: 
either the cell is valid for a polynomial degree greater than $0$, or the degree zero is reached. In the worst case, 
a cell is updated with a robust and stable first order accurate Godunov-type finite volume scheme, which is supposed to produce always 
valid (monotone and positivity-preserving) solutions under CFL condition. This \textit{a posteriori} check and order decrementing loop is called the 
'MOOD loop'. We refer the reader to \cite{CDL1,CDL2,CDL3,ADER_MOOD_14} for more details. 


The link between the MOOD concept developed in the finite volume framework and the typical strategy adopted in a classical DG limiter seems obvious. In both
approaches first a candidate solution is computed using an \textit{unlimited} scheme. Then, troubled cells are detected based on some criteria and the discrete 
solution is corrected. However, there are important differences: the typical DG limiters \textit{postprocess} the candidate solution using a nonlinear 
\textit{reconstruction} technique. They furthermore use only \textit{one} time level (the current one) for the detection and postprocessing step. In contrast, 
the MOOD approach first of all uses \textit{two} time levels (the old one and the current one) for the detection of troubled cells. Second, the MOOD approach  
\textit{recomputes} the solution using a \textit{different numerical scheme}, which is supposed to be more robust at shock waves. The fact of recomputing the 
solution and of looking at two different time levels  for the detection of troubled cells makes it in the notation of the authors of \cite{CDL1,CDL2,CDL3,ADER_MOOD_14} 
an \aposteriori approach.



The key innovation of the present paper is now the use of the aforementioned \aposteriori MOOD paradigm as a limiter for high order DG schemes. Since simple order 
decrementing as in the MOOD approach for finite volume schemes would obviously destroy the natural subcell resolution capability of DG, a more sophisticated strategy is 
needed here. We therefore suggest to recompute the solution of troubled cells on a finer \textit{subgrid} inside each cell, using a more robust but still very accurate 
one-step ADER-WENO finite volume scheme \cite{titarevtoro,Balsara2013934,AMR3DCL}. The data can be scattered from the main grid to the subgrid and gathered back 
via appropriate subcell projection and subcell reconstruction operators. Such operators are in principle well-known from high order finite volume schemes and spectral
finite volume schemes \cite{spectralfv2d,spectralfv3d,spectralfv.dg}. The choice of the subgrid size is very important. In this paper we suggest to choose the subgrid
size so that the local CFL number on the subgrid is as large as possible, hence the maximum admissible time step size of the finite volume scheme on the subgrid matches 
the maximum admissible time step of the DG scheme on the main grid. This leads to $N_s=2N+1$ subcells per space dimension for a DG scheme using a piecewise polynomial
approximation of degree $N$. For alternative subcell methods used as limiter for the DG method see \cite{CasoniHuerta2,Sonntag}. However, none of these 
uses our \aposteriori detection concept, nor do they recompute the solution in troubled cells via a better than second order accurate subcell finite volume method. 
Furthermore, the method proposed in \cite{Sonntag} does not use the optimal subgrid size that allows to get the maximum admissible CFL number on the subgrid. 


The rest of this paper is organized as follows. 
The one-step discontinuous Galerkin scheme used in this paper (ADER-DG) is presented in section \ref{sec:framework}. 
In section \ref{sec:apriori_limiter} we give a thorough but still non-exhaustive overview of existing \apriori limiters for the 
DG method, such as slope/moment reduction, artificial viscosity or WENO-like techniques. Then in section \ref{sec:aposteriori_limiter} we 
present all details of our new subcell-based \aposteriori approach. Next, section \ref{sec:numerics} gathers all numerical results for a 
large set of different test cases in order to assess the validity and robustness of our new \aposteriori subgrid limiter. Smooth and non-smooth 
test cases are simulated. Finally, conclusions and perspectives are drawn in section \ref{sec:conclusion}.



\section{The ADER Discontinuous Galerkin Method} \label{sec:framework}

In this paper we consider nonlinear systems of hyperbolic conservation laws in multiple space dimensions of the form 
\begin{equation}
\label{eqn.pde.nc}
\frac{\partial \Q}{\partial t} + \nabla \cdot \F\left(\Q\right) = 0,  
\qquad \x \in \Omega \subset \mathds{R}^d, \quad t \in \mathds{R}_0^+,
\end{equation}
with appropriate initial and boundary conditions
\begin{equation}
\label{eqn.pde.bc}
 \Q(\x,0) = \Q_0(\x), \quad \forall \x \in \Omega, \qquad \quad 
 \Q(\x,t) = \Q_B(\x,t) \quad \forall \x \in \partial \Omega, \quad t \in \mathds{R}_0^+,
\end{equation}
where $\Q \in \Omega_Q \subset \mathds{R}^\nu$ is the state vector of $\nu$ conserved quantities, 
and $\F(\Q) = (\f,\g, \h)$ is a non-linear flux tensor that depends on state $\Q$. 
$\Omega$ denotes the computational domain in $d$ space dimensions 
whereas $\Omega_Q$ is the space of physically admissible states, also called state space or phase-space. \\
We solve this system of equations by applying the general family of $\PNM$ methods introduced in 
\cite{Dumbser2008}, which provides high order of accuracy in both space and time.
However, in this paper we only focus on the family of pure Discontinuous Galerkin (DG) schemes, 
i.e. the $\PMM$ schemes in the context of \cite{Dumbser2008}. The numerical method is formulated 
as a one-step predictor corrector method \cite{GassnerDumbserMunz}: in the predictor step 
\eqref{eqn.pde.nc} is solved within each element \textit{in the small} (see \cite{eno}) 
by means of a locally implicit space-time discontinuous Galerkin scheme. 
The final time update of the discrete solution is explicit and is obtained by the one-step corrector. 
In the following we only summarize the main steps, while for more details the reader is referred to 
\cite{Dumbser2008,DumbserZanotti,HidalgoDumbser,GassnerDumbserMunz,Balsara2013934}.

\subsection{Space discretization and data representation}

The computational domain $\Omega$ is discretized by a Cartesian grid of conforming 
elements $T_i$, where the index $i$ ranges from 1 to the total number of elements $N_E$. 
The elements are chosen to be quadrilaterals in 2D and hexahedrons in 3D. The union of 
all elements represents the computational grid or the \textit{main grid} of the domain,  
\begin{equation}
\label{eqn.tetdef}
 \mathcal{T}_{\Omega} = \bigcup \limits_{i=1}^{N_E} T_i.
\end{equation}
We denote the cell volume by $|T_i| = \int_{T_i} d\x$.
At the beginning of each time-step, the numerical solution of equation \eqref{eqn.pde.nc} 
for the state vector $\Q$ is represented within each cell $T_i$ of the main grid  
by piecewise polynomials of maximum degree $N \geq 0$ and is denoted by 
$\u_h(\x,t^n) \in \mathcal{U}_h$,   
\begin{equation}
\label{eqn.ansatz.uh}
  \u_h(\x,t^n) = \sum_l\Phi_l(\x) \hat{\u}^n_l,  \quad \x \in T_i ,
\end{equation}
where $\u_h$ is referred to as the discrete ``representation'' of the solution. 
The space $\mathcal{U}_h$ of piecewise polynomials up to degree $N$ is spanned by 
the basis functions $\Phi_l=\Phi_l(\x)$. Throughout this paper we use the Lagrange 
interpolation polynomials passing through the tensor-product Gauss-Legendre 
quadrature points \cite{stroud} as spatial basis functions, see also 
\cite{Kopriva,GassnerKopriva,KoprivaGassner}. \\

\subsection{Local space-time predictor}
\label{sec.galerkin}

The representation polynomials $\u_h(\x,t^n)$ are now evolved in time according to a
local weak formulation of the governing PDE in space-time, see
\cite{DumbserEnauxToro,Dumbser2008,HidalgoDumbser,DumbserZanotti,GassnerDumbserMunz,Balsara2013934}.
The local space-time Galerkin method is only used for 
the construction of an element-local predictor solution of the PDE neglecting the influence 
of the neighbors. This predictor solution is further inserted into a corrector step described in the next section, 
which then provides the appropriate coupling between neighbor elements via a numerical flux function (Riemann solver).  

Let us first transform the PDE \eqref{eqn.pde.nc} into a space-time reference coordinate system 
$(\boldsymbol{\xi}, \tau)$ of the space-time reference element $[0;1]^{d+1}$, 
with $\nabla_\xi = \partial \boldsymbol{\xi} / \partial \x \cdot \nabla$. The spatial reference 
elements are denoted by $T_E$ and are defined as $T_E=[0;1]^d$, i.e. the unit square  
in two space dimensions and the unit cube in three space dimensions, respectively. Time is
transformed according to $t=t^n + \Delta t \, \tau$. 
This yields 
\begin{equation}
\label{eqn.pde.nc.2d}
 \frac{\partial \Q}{\partial \tau}
    + \nabla_\xi \cdot \F^* \left( \Q \right) = 0,
\end{equation}
with the modified flux
\begin{equation}
  \F^* := \Delta t \left( \partial \boldsymbol{\xi} / \partial \x  \right)^{T} \cdot \F(\Q).
\end{equation}

To simplify notation, let us define the following two operators  
\begin{equation}
 \label{eqn.operators1}
  \left<f,g\right> =
      \int \limits_{0}^{1} \int \limits_{T_E} \left( f(\boldsymbol{\xi}, \tau) \cdot g(\boldsymbol{\xi}, \tau) \right) d \boldsymbol{\xi} \, d \tau,
\end{equation}
\begin{equation}
  \left[f,g\right]^{\tau} =
      \int \limits_{T_E} \left( f(\boldsymbol{\xi}, \tau) \cdot g(\boldsymbol{\xi}, \tau) \right) d \boldsymbol{\xi},
\end{equation}
which denote the scalar products of $f$ and $g$ over the space-time reference element $T_E \times \left[0;1\right]$ and  
over the spatial reference element $T_E$ at time $\tau$, respectively. 

Now we multiply \eqref{eqn.pde.nc.2d} with a space-time test function 
$\theta_k=\theta_k(\boldsymbol{\xi},\tau)$ and subsequently integrate over 
the space-time reference control volume $T_E \times [0;1]$ to obtain the weak formulation 
\begin{equation}
\label{eqn.pde.nc.weak1}
 \left< \theta_k, \frac{\partial \q_h}{\partial \tau}  \right>
    + \left< \theta_k, \nabla_\xi \cdot \F_h^* \left(\q_h\right) \right>
    = 0.
\end{equation}
The discrete solution of equation \eqref{eqn.pde.nc.weak1} in space and time, denoted by $\q_h$, 
as well as the discrete space-time representation of the flux tensor $\F_h^*$ are assumed to have 
the following form 
\begin{equation}
\label{eqn.st.state}
 \q_h = \q_h(\boldsymbol{\xi},\tau) =
 \sum \limits_l \theta_l(\boldsymbol{\xi},\tau) \hat{\q}_l := \theta_l \hat{\q}_l,
\end{equation}
\begin{equation}
\label{eqn.st.flux}
 \F^*_h = \F^*_h(\boldsymbol{\xi},\tau) =
 \sum \limits_l \theta_l(\boldsymbol{\xi},\tau) \hat{\F}^*_l := \theta_l \hat{\F}^*_l,
\end{equation}
where we have used the Einstein summation convention over two repeated indices. 
If one employs a nodal basis as in \cite{Dumbser2008}, one simply has 
\begin{equation}
  \hat{\F}^*_l = \F^*(\hat{\q}_l).
\end{equation}
Throughout this paper we use space-time basis functions that are the Lagrange interpolation polynomials
passing through the tensor-product Gauss-Legendre quadrature points on the space-time reference
element $[0;1]^{d+1}$, see \cite{AMR3DCL}. In that way, all the resulting matrices have a sparse block structure and the 
computations can be done efficiently in a dimension-by-dimension fashion.  
Integration by parts in time of the first term in \eqref{eqn.pde.nc.weak1} yields
\begin{equation}
\label{eqn.pde.nc.dg1}
   \left[ \theta_k, \q_h \right]^1 - \left[ \theta_k, \u_h \right]^0
   - \left< \frac{\partial}{\partial \tau} \theta_k, \q_h  \right>
    + \left< \theta_k, \nabla_\xi \cdot \F^*_h \right>
    = 0.
\end{equation}
Note that the piecewise high order polynomial representation $\u_h$ is taken into account 
as initial condition of the Cauchy problem in the small in a weak sense
by the term $\left[ \theta_k, \u_h \right]^0$. \\
Eqn. \eqref{eqn.st.state} is then substituted into \eqref{eqn.pde.nc.dg1} and yields the following 
iterative scheme, see \cite{Dumbser2008,HidalgoDumbser,DumbserZanotti} for more details, 
\begin{equation}
\label{eqn.pde.nc.dg2}
   \left( \left[ \theta_k, \theta_l \right]^1 
   - \left< \frac{\partial}{\partial \tau} \theta_k, \theta_l \right> \right)
    \hat{\q}_l^{r+1}    
    =  \left[ \theta_k, \Phi_l \right]^0 \hat \u_l^n - \left< \theta_k, \nabla_\xi \theta_l \right> \cdot \F^*(\hat{\q}_l^r),
\end{equation}
where $r$ is the iteration index. The iterative method \eqref{eqn.pde.nc.dg2} converges very efficiently to the unknown expansion 
coefficients $\hat{\q}_l$ of the local space-time predictor solution, see \cite{Dumbser2008,HidalgoDumbser,DumbserZanotti}. 
Once the $\hat{\q}_l$ are obtained the local space-time predictor $\q_h$ is known inside each cell (\ref{eqn.st.state}). 
The above iterative procedure has replaced the cumbersome Cauchy-Kovalewski procedure that has been initially employed in the 
original version of the ADER finite volume and ADER discontinuous Galerkin schemes 
\cite{schwartzkopff,toro3,toro4,titarevtoro,dumbser_jsc,taube_jsc,DumbserKaeser07}. 

\subsection{Fully discrete one-step ADER-DG scheme}
\label{sec.ADERNC}

The fully discrete one-step ADER-DG scheme \cite{dumbser_jsc,taube_jsc,QiuDumbserShu} is obtained after multiplication of the governing PDE \eqref{eqn.pde.nc} by a test 
function $\Phi_k \in \mathcal{U}_h$, which is identical with the spatial basis functions, and subsequent integration over the space-time control volume 
$T_i \times [t^n;t^{n+1}]$. 
The flux divergence term is then integrated by parts and one obtains the weak formulation 
\begin{equation}
\label{eqn.pde.nc.gw1}
\int \limits_{t^n}^{t^{n+1}} \int \limits_{T_i} \Phi_k \frac{\partial \u_h}{\partial t} d\x dt
+\int \limits_{t^n}^{t^{n+1}} \int \limits_{\partial T_i} \Phi_k \, \F\left(\u_h \right)\cdot\mathbf{n} \, dS dt 
-\int \limits_{t^n}^{t^{n+1}} \int \limits_{T_i} \nabla \Phi_k \cdot \F\left(\u_h \right) d\x dt 
= 0, 
\end{equation}
where $\mathbf{n}$ is the outward pointing unit normal vector on the surface $\partial T_i$ of element $T_i$. 
As usual in the DG finite element framework \cite{cbs0,cbs1,cbs2,cbs3,cbs4} the boundary flux term in \eqref{eqn.pde.nc.gw1} 
is then replaced by a numerical flux function (Riemann solver) in normal direction, 
$\mathcal{G}\left(\q_h^-, \q_h^+ \right) \cdot\mathbf{n}$, which is a function of the left and right boundary-extrapolated data, 
$\q_h^-$ and $\q_h^+$, respectively. 
Inserting the local space-time predictor 
$\q_h$ into \eqref{eqn.pde.nc.gw1} then yields the following arbitrary high order accurate one-step Discontinuous Galerkin (ADER-DG) scheme: 
\begin{equation}
\label{eqn.pde.nc.gw2}
\left( \int \limits_{T_i} \Phi_k \Phi_l d\x \right) \left( \hat{\u}_l^{n+1} -  \hat{\u}_l^{n} \right) +
\int\limits_{t^n}^{t^{n+1}} \int_{\partial T_i} \Phi_k \, \mathcal{G}\left(\q_h^-, \q_h^+ \right)\cdot\mathbf{n} \, dS dt 
-\int\limits_{t^n}^{t^{n+1}} \int_{T_i} \nabla \Phi_k \cdot \F\left(\q_h \right) d\x dt  = 0. 
\end{equation}
In general we use the simple Rusanov (local Lax Friedrichs) flux \cite{Rusanov:1961a}, or the 
Osher-type flux recently proposed in \cite{OsherUniversal} 
as numerical flux function at the element boundaries, although any other kind
of Riemann solver could be also considered, see \cite{toro-book} for an overview of state-of-the-art Riemann solvers. \\

\subsection{Time step restriction} \label{sec.dt}
The DG method applied to convective problems enjoys a Courant-Friedrichs-Lewy (CFL) number that decreases with the 
approximation order $N$, roughly it follows $\frac{1}{2N+1}$, \cite{CBS-book,CBS-convection-dominated}. This rather restrictive 
condition is caused by the growth of the spectrum of the spatial discretization operator of the semi-discrete scheme  
\cite{krivodonova2013analysis}, hence the time step of an explicit DG scheme in multiple space dimensions has to satisfy 
\begin{equation}
\Delta t \leq \frac{1}{d} \frac{1}{(2N+1)} \frac{h}{|\lambda_{\max}|},
\label{eqn.cfl} 
\end{equation}
where $h$ and $|\lambda_{\max}|$ are a characteristic mesh size and the maximum signal velocity, respectively.  
This rather restrictive condition on the time step is to be put in perspective with the 
subcell resolution capability of any DG method, which allows the use of \textit{coarse} or even 
\textit{very coarse} grids.  
Recall that the representation polynomials are in principle of arbitrary degree $N \geq 0$.
Consequently, within one cell multiple subcell features can be captured. We refer the reader to the next 
section for a detailed description about subcell resolution. \\

This closes the brief description of the family of one-step ADER Discontinuous Galerkin schemes that is 
further used in this paper. Note that this family of schemes is the \textit{unlimited} version 
which, as such, can not prevent numerical oscillations and/or Gibbs phenomenon from occurring  
in the presence of steep gradients or shock waves. 
For sufficiently smooth solutions the unlimited ADER-DG scheme has shown an effective order of accuracy of
$N+1$, linked to the representation space $\mathcal{U}_h$. See \cite{Dumbser2008} for convergence tests in 
2D and 3D. 
Nonetheless, adding some sort of limiting is of paramount importance to ensure the numerical stability of  
the method, despite the nonlinear stability properties that can be proven for the DG scheme in $L_2$ norm and 
the associated local cell entropy inequality \cite{jiangshu}.  
It is the purpose of the next two sections to discuss the issue of limiters of DG schemes and to propose a new 
\aposteriori subcell-based finite volume limiting strategy.



\section{Nonlinear stability via \apriori limiting}
\label{sec:apriori_limiter}

In this section we briefly review existing limiters for DG that have been designed
upon paradigms such as slope limiters, artificial viscosity or essentially 
non oscillatory reconstruction procedures, such as ENO, WENO or HWENO. A discussion  
follows that will give birth to the paradigms and design principles of our 
\aposteriori subcell-based MOOD limiting presented in the next section.

\subsection{A brief and non-exhaustive review of \textit{a priori} limiters for DG schemes} \label{ssec:apriori_limiter}

The main difficulty of solving nonlinear hyperbolic systems of conservation laws arises due to the 
fact that solutions of the system may become discontinuous, even if the initial conditions are smooth. 
This important discovery was due to the groundbreaking work of Bernhard Riemann \cite{riemann1,riemann2}. 
Despite their provable nonlinear $L_2$ stability \cite{jiangshu}, even DG schemes 
can fail in the presence of strong shock waves or steep gradients and may generate 
strong oscillations that can ultimately lead to a failure of the computation. 
This is a consequence of the Godunov theorem \cite{godunov} that also applies to the DG method. 
Therefore, some sort of non-linear limiting is needed and a vast literature on designing such techniques 
for DG methods does exist. Being exhaustive is almost impossible. As a consequence we only recall the underlying 
basic principles and some of the most important limiters used. Philosophically speaking, any limiter 
procedure for DG answers the following questions:
\begin{itemize} 
\item 
  Q1: \textit{Where are the locations where limiting is needed?} 
  This is the so-called \textit{troubled zone} indicator. We emphasize that the number of these 
	locations may be very small, since for high order DG schemes many features can be captured 
	within one cell length, see Fig.~\ref{fig:P7}. 
\item 
  Q2: \textit{How do we achieve high order of accuracy along with a non-oscillatory 
    property close to these locations?} 
  Ideally one should manipulate or replace the DG polynomials in such 
  a way that additional numerical dissipation is supplemented close to the detected 
	locations, but nowhere else, and preferably without destroying the subcell resolution of the DG method.
\end{itemize}
The design of such \apriori limiters is a difficult task.
However, amongst all DG limiters, three families seem to emerge:
artificial viscosity based limiters such as 
\cite{Hartman_02,Persson_06,Barter_2010,Nguyen11,spacetimedg1,spacetimedg2,Feistauer4,Feistauer5,Feistauer6,Feistauer7},
``Slope'' or moment limiters, for instance 
the total variation bounded (TVB) limiter of Cockburn and Shu \cite{cbs0,cbs1,cbs2,cbs3}, 
the moment-based limiters \cite{Biswas_94,Burbeau_2001,Yang_AIAA_2009,Yang_parameterfree_09}, 
or the hierachical slope limiters \cite{Kuzmin2010,Kuzmin2013,Kuzmin2014},
as well as WENO and HWENO based limiters, such as the ones developed in \cite{QiuShu1,QiuShu2,QiuShu3,Qiu_2005,balsara2007,Luo_2007,Zhu_2008,Zhong_2013,Zhu_3D_12}. 
Let us briefly describe and comment the design principles of these three families.

\subsubsection{Artificial viscosity (AV)}

The idea of using an artificial viscosity (AV) to stabilize shock waves dates back 
to von Neumann and Richtmyer \cite{vonneumann50} during the Manhattan project
in the 1940's at Los Alamos National Laboratory.
This idea has been fruitful especially in the Lagrangian community since then
and a lot of AV models emerged, bulk AV, edge AV, tensor AV, see
pages 17-21 of \cite{Loubere_HDR_2013} and the reference herein for a review. 
Back to the origin, the illuminating idea of von Neumann and Richtmyer 
was to introduce a purely artificial dissipative mechanism
of such a form and strength that the shock transition would be a smooth one,  
extending over a small number of cell lengths, 
and then to include this dissipation into the finite 
difference equations, see \cite{neu}, page 312.
A lot of minor and/or major improvements have been developed since then
but any artificial viscosity technique revolves around the basic ideas: 
(i) define the region of the shock waves, 
(ii) add some dissipative mechanism over a small number of cells. \\
While it is nowadays rarely used in the context of shock capturing 
finite volume schemes, the artificial viscosity concept has become popular
again in the context of DG schemes in order to capture shocks. 
In \cite{Hartman_02} the authors use the magnitude of the residual to determine the amount 
of artificial viscosity added to the shock region. In \cite{Persson_06} the authors have 
introduced a subcell shock-capturing method based on the manipulation of the density variable 
for determining the shock region and also the magnitude of artificial viscosity to be added. 
Piecewise constant artificial viscosity leading to oscillations 
has later been overcome in \cite{Barter_2010} where the AV model is based on the system
of equations solved (one additional equation is solved for the AV). 
In the DG community researchers have improved 
the AV models adding artificial terms to the physical 
viscosity coefficients \cite{Cook:2004}
or using an analytic function of the dilatation 
\cite{Nguyen11} \footnote{ 
This phenomenon also occurred in the Lagrangian community.  
For instance, the function of the dilatation in 
\cite{Nguyen11} is related to
the so-called compression switch used in the Lagrangian community known 
from an unpublished work by Rosenbluth from Los Alamos in the 1950's,  
where he suggested that the artificial viscosity should be zero when the fluid 
is undergoing an expansion. This 'trick' is nowadays known as the 'artificial 
viscosity switch' \cite{Csw98,Caramana-Burton-Shashkov-Whalen-98,Cam01}.}.

\subsubsection{WENO limiting procedure}
Some other authors suggest the use of a (H)WENO limiting procedure for Runge-Kutta DG 
methods, for instance in \cite{QiuShu1,QiuShu2,QiuShu3,balsara2007,Zhong_2013,Zhu_13}. 
Usually these authors adopt the following framework:
\begin{enumerate}
\item Identify invalid cells (called ``troubled cells''), 
  namely, the cells which might need limiting. 
  In \cite{Qiu_2005} a very detailed comparison among different 
  possible troubled-cell indicators is carried out. 
  These indicators are a very important ingredient for triggering the subsequent WENO limiter 
	procedure and are often based on minmod-type slope limiters, such as the modified TVB minmod 
	limiter for the 2D scalar case \cite{cbs3}, the modified TVB minmod limiter in characteristic 
	variables proposed in \cite{cbs2} and \cite{cbs4} for nonlinear systems in one and multiple 
	space dimensions, respectively. 
  Sometimes, other shock detection techniques such as the one suggested in \cite{Krivodonova2004} 
	are used as troubled cell indicator. A review can be found in \cite{wang2011adaptive}, chapter 6. 
	A troubled zone indicator based on subcell information has been forwarded by Balsara et al. in 
	\cite{balsara2007}. 
\item Replace the DG polynomials in these detected cells with \textit{reconstructed} polynomials, 
  that keep the original cell averages, preserve the same high order of accuracy, and
  are less oscillatory through a more or less classical WENO reconstruction procedure, see 
	\cite{Zhong_2013} for the details. We underline that in the case of WENO limiters, 
	the subcell resolution property of the DG method is lost to some extent, since the higher order moments
	are reconstructed from the cell averages defined on the coarse main grid. 
\end{enumerate}
Note that in the same spirit Hermite WENO schemes have also been designed to be used 
as limiters for DG schemes, see e.g. \cite{Qiu_2004,Luo_2007,Zhu_2009,balsara2007}, to avoid the 
non-compactness of the WENO limiters, but they are following a similar idea. \\
Usually improving the identification of the troubled cells is a key point, as the procedure may rely on 
user-defined parameters, which are often problem dependent even if in \cite{Zhong_2013} the authors have 
shown a relative insensibility of the method to these parameters.
Last, the extension to 3D and more general meshes are important key points, see \cite{Zhu_3D_12}. 
More recently the troubled cell indicator has also been used for adaptive methods $h$ (mesh refinement),
$p$ (order enrichment), or $r$ (mesh motion) by refining the troubled cells and/or coarsening the 
others \cite{Zhu_hadap_2013}. 

\subsubsection{"Slope"/moment reduction}
Slope limiting or moment reduction techniques 
may permit to control the jumps of the DG polynomials by constraining 
or nullifying the high-order components in designated cells 
\cite{BarthJespersen,Burbeau_2001,cbs4,Yang_AIAA_2009,Yang_parameterfree_09}.
A well-designed slope limiter must filter out non-physical oscillations 
without sacrificing the order of accuracy at smooth extrema. 
To do so some authors rely on monotonicity-preserving limiters 
frequently combined with \textit{ad hoc} smoothness or oscillation detection procedures. 
Originally, slope limiters were developed to constrain piecewise linear polynomial 
reconstructions based on a discrete maximum principle, therefore constructing ``slope'' limiters 
for higher order polynomials are difficult to build into this paradigm \cite{Krivodonova2004,Michoski2011}. 
Many classical limiters have been already used, for instance the minmod-based TVB limiter \cite{cbs1},  
the family of moment limiters \cite{Biswas_94,Burbeau_2001}, as well as 
monotonicity preserving limiters \cite{SureshHuynh,Rider_2001}. 
Most of these limiters succeed in controlling spurious numerical oscillations, but they may also 
present a tendency to degrade accuracy when mistakenly used in certain smooth regions 
of the solution.
In \cite{Kuzmin2010,Kuzmin2013,Kuzmin2014} the authors have developed a 
so-called hierarchical slope limiter. Briefly this limiter constrains the derivatives of 
polynomials, written in the Taylor polynomial form, in order to eliminate under/overshoots measured 
at the vertices of the cell, see \cite{Kuzmin2014}. This vertex-based hierarchical slope limiter has 
some common features with the moment limiter from \cite{Krivodonova2004} and seems to preserve smooth 
extrema without using any troubled cell indicator. Recently this hierarchical slope reduction 
problem has also been recast into an optimization problem \cite{Kuzmin2014}.

\subsection{Discussion}



Most of the aforementioned shock capturing procedures for DG schemes have a common feature: they rely on the fact 
that spurious numerical oscillations can be detected and corrected in the discrete solution by looking at only one 
time step $t^n$ and usually \textit{without using the PDE}. The \textit{a priori} character is evident for artificial viscosity 
based approaches, where the shock detector and the magnitude of the artificial viscosity are chosen based on the current 
solution $\u_h(\x,t^n)$. The TVB, (H)WENO and moment limiters can be to some extent considered as \textit{a posteriori} 
limiter techniques, since they modify the higher order moments in troubled cells at the end of each time step 
(or Runge-Kutta stage) after having used an unlimited version of the DG scheme.
Nonetheless, the troubled cell indicator as well as the limiter usually only consider the \textit{current} discrete solution 
and higher order moments are replaced by some sort of nonlinear \textit{data reconstruction} based on the degrees of 
freedom of the current solution, \textit{without using the PDE}. 

In a DG method many features can be captured within one characteristic cell length. Certain limiters may therefore dissipate 
a lot of these subscale features, and only subtle ones would maintain at the same time the high subcell accuracy and assure  
stability of the numerical solution at shock waves. 
According to numerical evidence provided in \cite{CasoniHuerta1} for the one-dimensional case, the artificial 
viscosity approach \cite{Persson_06} seems to be more appropriate for capturing subcell features compared 
to other DG limiters, especially for very high polynomial degrees $N$. 
Ultimately, one may demand that a DG limiter acts only on the smallest length scale within one cell to 
avoid that excessive numerical dissipation impacts all features that are represented by $\u_h$. This smallest 
length scale is related to the size of the cells and to the polynomial degree $N$. \\ 

Recent developments have been made concerning the construction of sub-cell limiters, which use either a finite volume method 
on subcells \cite{Sonntag} or a smooth switch between a high order DG scheme and a first order finite volume subgrid method,  
\cite{CasoniHuerta1,CasoniHuerta2}, but none of these is based on the a posteriori concept proposed in this paper, nor do 
they rely on the use of higher order finite volume schemes on the subgrid. 


Following the ideas of a novel \aposteriori detection approach (MOOD), that has been introduced for the first time in 
the context of very high order accurate finite volume schemes in \cite{CDL1,CDL2,CDL3,ADER_MOOD_14}, we propose to 
extend the MOOD concept to the context of very high order accurate DG finite element schemes in the following.


\section{Nonlinear stability via \textit{a posteriori} sub-cell limiting (SCL)}
\label{sec:aposteriori_limiter}

\subsection{General design principles}

The subcell resolution is also evident from the amount of data that are stored per cell in the DG context,
namely from the number of the degrees of freedom, which is a function of the polynomial degree $N$.
For our tensor-product basis functions, the number of spatial degrees of freedom per cell is $(N+1)^d$. 

As already discussed before, limiting acts in two steps: First the troubled cell indicator detects which 
regions of the computational domain need limiting and, second, the limiting effectively adds some sort of 
numerical dissipation in these regions, either directly via artificial viscosity or via a nonlinear data 
reconstruction or slope limiting procedure. However, an appropriate DG limiter should ideally detect and 
correct problematic situations \textit{on a subscale level}. 

Starting with piecewise high order polynomial data $\u_h(\x,t^n)$ it seems difficult to detect 
\apriori problematic situations which will occur between time $t^n$ and time $t^{n+1}$.  
In \cite{CDL1,CDL2,CDL3,ADER_MOOD_14} in the context of high order finite volume (FV) 
schemes the authors have adopted a different strategy, called MOOD. This latter consists 
in testing \aposteriori a so-called \textit{candidate solution} $\u_h^*(\x,t^{n+1})$ that has been obtained 
using an \textit{unlimited} high order scheme.  
If this candidate solution does not fulfill a set of properties called \textit{detection criteria}, 
then the cells for which it has failed are \textit{recomputed} using a more robust and more dissipative lower 
order scheme (decrementing of the order of accuracy of the scheme). For an invalid cell, the iterative MOOD loop 
ends either with a valid solution that has passed the detection criteria, or, with a solution that is updated 
with the lowest order scheme that is supposed to be monotone and positivity preserving under CFL condition. 
The detection criteria developed in \cite{CDL1,CDL2,CDL3,ADER_MOOD_14} for different hyperbolic systems of 
conservation laws have proven to be sensitive enough to avoid excessive numerical diffusion while maintaining 
an effective high order of accuracy for smooth problems. More important, they are sufficient to dissipate 
numerical oscillations and ensure stability. \\

\subsection{Subcell data representation, projection and reconstruction}

\begin{figure}
  \begin{center} 
    \begin{tabular}{c} 
      \includegraphics[width=0.85\textwidth]{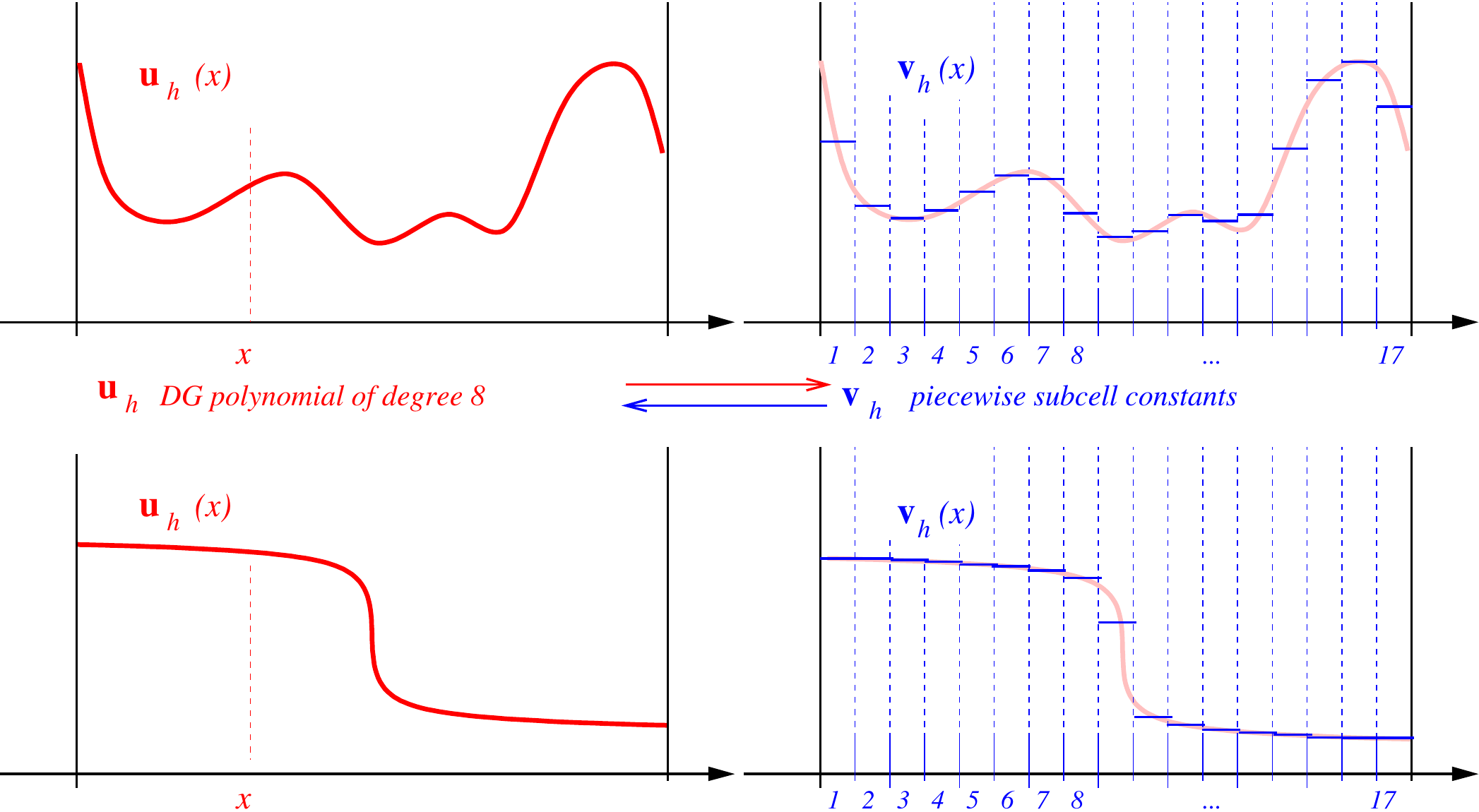} 
    \end{tabular}
    \caption{ \label{fig:P7}
      Examples of DG polynomials $\u_h$ on a cell (red) and associated projection
      $\v_{h}=\mathcal{P}(\u_h)$ onto subcell averages (blue). The information 
			contained in $\u_h$ can be recovered from $\v_{h}$ via the subcell reconstruction 
			operator for $N_s \geq N+1$. Throughout this paper we use $N_s=2N+1$ subcells.}
  \end{center}
\end{figure}

If $\u_h(\x,t^n)$ is the data representation of the DG scheme within cell $T_i$ at time $t^n$ and we consider 
a fine \textit{subgrid} of $T_i$, denoted by $\mathcal{S}_i = \bigcup_j S_{i,j}$ made of $\left(N_s\right)^d$ 
subcells called $S_{i,j}$, $j=1,\ldots, (N_s)^d$, with $N_s \geq N+1$, then we introduce an 
\textit{alternative data representation} denoted by $\v_{h}(\x,t^n)$, which is defined by a set of 
\textit{piecewise constant subcell averages} $\v_{i,j}^n$. These subcell averages 
are directly computed from $\u_h(\x,t^n)$ by $L_2$ projection, which in this case simply means the computation
of the integral average of $\u_h(\x,t^n)$ since $\v_{h}(\x,t^n)$ is piecewise constant on the subcells $S_{i,j}$: 
\begin{equation}
\label{eqn.subcell.p} 
  \v_{i,j}^n = \frac{1}{|S_{i,j}|} \int \limits_{S_{i,j}} \u_h(\x,t^n) d\x = \frac{1}{|S_{i,j}|} \int \limits_{S_{i,j}} \phi_l(\x) d\x \, \hat \u_l^n, \qquad \forall S_{i,j} \in \mathcal{S}_i. 
\end{equation} 
The above equation \eqref{eqn.subcell.p} is in the following also called \textit{projection operator} $\mathcal{P}$, denoted by
$\v_{h}(\x,t^n) = \mathcal{P}\left(\u_h(\x,t^n)\right)$. Throughout this paper, the subcells are chosen to be
\textit{equidistant} Cartesian subcells \cite{AMR3DCL}. 
To gather back the piecewise constant subcell data into a high order DG polynomial, we apply the following 
\textit{reconstruction operator} $\mathcal{R}$: find $\u_h(\x,t^n)$ and therefore $\hat \u_l^n$ so that  
\begin{equation}
\label{eqn.subcell.r} 
  \int \limits_{S_{i,j}} \u_h(\x,t^n) d\x = \int \limits_{S_{i,j}} \v_h(\x,t^n) d\x, \qquad \forall S_{i,j} \in \mathcal{S}_i.  
\end{equation} 
This is a classical reconstruction or recovery problem of a higher order polynomial from known cell averages, which typically 
arises within the finite volume context and also in spectral finite volume methods \cite{spectralfv2d,spectralfv3d,spectralfv.bnd}. 
The approach presented in this paper explicitly admits $N_s > N+1$, hence the resulting system may be \textit{overdetermined}. This 
overdetermined system of reconstruction equations is then solved using a \textit{constrained} least-squares approach, see 
\cite{kaeserjcp,DumbserKaeser06b}, where the constraint is the integral conservation of the cell average over the big cell 
$T_i$, i.e.  
\begin{equation}
\label{eqn.subcell.r.constraint} 
  \int \limits_{T_{i}} \u_h(\x,t^n) d\x = \int \limits_{T_{i}} \v_h(\x,t^n) d\x. 
\end{equation} 
The operator given by the solution of \eqref{eqn.subcell.r} and \eqref{eqn.subcell.r.constraint} is in the following denoted by 
$\u_{h}(\x,t^n) = \mathcal{R}\left(\v_h(\x,t^n)\right)$. It is obvious that finding $\u_h(\x,t^n)$ from \eqref{eqn.subcell.r} 
corresponds to the computation of the (pseudo-) inverse of the matrix associated with the projection operator \eqref{eqn.subcell.p}. 
Hence, the two operators $\mathcal{R}$ and $\mathcal{P}$ satisfy the property 
\begin{eqnarray}
\label{eqn.rid} 
 \mathcal{R} \circ \mathcal{P} = \mathcal{I},
\end{eqnarray}
where $\mathcal{I}$ is the identity operator. 
If the subgrid is large enough ($N_s \geq N+1$) then there is no loss of accuracy when applying the projection operator,  
since $\u_{h}(\x,t^n)$ can always be recovered identically from $\v_{h}(\x,t^n)=\mathcal{P}\left( \u_{h}(\x,t^n) \right)$ 
using the recovery operator $\mathcal{R}$ due to relation \eqref{eqn.rid}, thus $\v_{h}(\x,t^n)$ is able to represent all the 
information contained in $\u_h(\x,t^n)$, see Fig. \ref{fig:P7}. 
If $N_s=N+1$, like in \cite{Sonntag}, then the number of subcells corresponds exactly to the number of degrees of freedom associated 
with the space $\mathcal{U}_h$. 

Consequently, it is equivalent in terms of nominal accuracy to represent data either with a DG scheme of polynomial degree $N$ on the cell $T_i$, 
or, representing the data by piecewise constant cell averages on the subgrid $\mathcal{S}_{i,j}$. 

The use of subgrid information has been used very successfully also in the context of semi-implicit finite volume methods for free surface 
flows, see \cite{Casulli2009,CasulliStelling2011}, where it has led to a significant improvements in terms of numerical accuracy and 
computational efficiency.

\subsection{Extension of MOOD to DG schemes}


The a posteriori MOOD concept is now extended to the DG context as follows. First, a \textit{candidate solution} 
$\u_h^*(\x,t^{n+1})$ is computed from $\u_h(\x,t^n)$ by the \textit{unlimited} DG scheme \eqref{eqn.pde.nc.gw2}. 
Next, we apply the following \textit{detection criteria}. 

\paragraph*{Physical admissibility detection (PAD)} 
The detection criteria must contain physics-based admissibility properties, and, as such can not be uncorrelated with the hyperbolic system 
of conservation laws which is solved.  
Hence, a candidate solution $\u_h^*(\x,t^{n+1})$ is said to be physically valid inside  
cell $T_i$ for this system if 
\begin{equation}
  \pi_k \left( \u_h^*(\x,t^{n+1}) \right) > 0, \qquad \forall \x \in T_i,  \, \forall k, 
	\label{eqn.positivity}
\end{equation} 
where $\pi_k(\Q)$ is the physical quantity that must satisfy the positivity constraint and which is a function of the state vector $\Q$.  
For the Euler equations of compressible gas dynamics, for example, the mass density $\rho$ and the fluid pressure $p$ must be 
positive everywhere and for all times, hence $\pi_1(\Q)=\rho$ and $\pi_2(\Q)=p$.  
There is a growing interest in designing high order accurate and positivity preserving finite volume and DG methods. For the most recent 
developments in the field of DG schemes see, for example, \cite{ShuPositivity1,ShuPositivity2,ShuPositivity3}. 

\paragraph*{Numerical admissibility detection (NAD)} In the past, the discrete maximum principle (DMP) was a very successful guideline for the 
construction of high resolution shock capturing schemes. In this paper, we therefore use the following \textit{relaxed} version of a 
discrete maximum principle, which takes into account the data representation of the DG method under the form of piecewise polynomials.
It is a natural extension of the DMP for cell averages used in \cite{CDL1,CDL2,CDL3,ADER_MOOD_14}.  
The DMP is applied in an \textit{a posteriori} manner as follows. A candidate solution $\u_h^*(\x,t^{n+1})$ is said to fulfill the numerical 
admissibility detection criterion in cell $T_i$ if the following relation is fulfilled componentwise for all conserved variables:  
\begin{equation}
  \min_{\y \in \mathcal{V}_i} \left( \u_h(\y,t^{n}) \right) - \boldsymbol{\delta} \leq \u_h^*(\x,t^{n+1}) \leq  \max_{\y \in \mathcal{V}_i} \left( \u_h(\y,t^{n}) \right) + \boldsymbol{\delta}, \qquad \forall \x \in T_i,   
\label{eqn.dmp.org} 
\end{equation} 
where $\mathcal{V}_i$ is a set containing element $T_i$ together with its Voronoi neighbor cells that share a common node with $T_i$. 
We see from \eqref{eqn.dmp.org} that the discrete maximum principle is now applied in the sense of polynomials. The polynomial that represents the  
candidate solution on element $T_i$ must remain between the minimum and the maximum of the polynomials that have represented the discrete solution 
at the old time step in the neighborhood $\mathcal{V}_i$ of cell $T_i$. 

The quantity $\boldsymbol{\delta}$ is used to relax the strict maximum principle in order to allow some very small overshoots and 
undershoots and to avoid problems with roundoff errors that would occur when applying \eqref{eqn.dmp} in a strict way, without using 
$\boldsymbol{\delta}$. Throughout this paper we set    
\begin{equation}
 \boldsymbol{\delta} = \epsilon \cdot \left( \max_{\y \in \mathcal{V}_i} \left( \u_h(\y,t^{n}) \right) - \min_{\y \in \mathcal{V}_i} \left( \u_h(\y,t^{n}) \right) \right), 
\label{eqn.relax} 
\end{equation}
with $\epsilon = 10^{-3}$. We underline that the slight relaxation of the maximum principle has no influence on the positivity of the solution, because the positivity is detected 
separately under the PAD above. Since it is not very practical from a computational point of view to calculate the maximum and the minimum of the discrete solution $\u_h(\x,t^{n})$ 
within the neighborhood $\mathcal{V}_i$ exactly, 
we use the following discrete version of \eqref{eqn.dmp.org} on the subscale level instead, which can be easily evaluated on the basis of subcell averages: 
\begin{equation}
  \min_{\y \in \mathcal{V}_i} \left( \v_h(\y,t^{n}) \right) - \boldsymbol{\delta} \leq \v_h^*(\x,t^{n+1}) \leq  \max_{\y \in \mathcal{V}_i} \left( \v_h(\y,t^{n}) \right) + \boldsymbol{\delta}, \qquad \forall \x \in T_i,   
\label{eqn.dmp} 
\end{equation} 
with $\v_h^*(\x,t^{n+1}) = \mathcal{P}\left( \u_h^*(\x,t^{n+1}) \right)$ and $\v_h(\y,t^{n})$ given by \eqref{eqn.ader.icbc}. 
A candidate solution is said to be \textit{valid} inside cell $T_i$, if it has passed both the physical and the numerical admissibility detection criteria, i.e. if \eqref{eqn.positivity} 
and \eqref{eqn.dmp} are satisfied. In this case, we set a cell-based indicator function $\beta$ to $\beta_i^{n+1}=0$. If a cell does \textit{not} fulfill the above \aposteriori MOOD 
detection criteria, a cell is marked or flagged like in a typical troubled zones indicator by setting $\beta_i^{n+1}=1$. However, we emphasize again that in our approach the detector uses 
information from \textit{two} different time levels, $t^n$ and $t^{n+1}$, while a classical troubled zone indicator would only look at the discrete solution at \textit{one} time level.   
An alternative self-adjusting \textit{a priori} strategy to trigger the subcell finite volume limiter presented in this paper could be the flattener algorithm 
proposed in \cite{PositivityWENO2}. 

After the detection phase given by \eqref{eqn.positivity} and \eqref{eqn.dmp}, the next operation of our subcell limiter consists in updating the discrete 
solution in invalid cells using a more robust scheme on the subgrid and based on the alternative data representation $\v_h(\x,t^n)$. For this task, one could 
in principle choose a simple and cheap second order TVD finite volume scheme, as used in \cite{Sonntag}, but we prefer to use a higher order one-step ADER-WENO 
finite volume method to avoid the clipping of local extrema on the subgrid. This is particularly important if local extrema on the subgrid are associated 
with physical phenomena, such as sound waves. The WENO scheme on the subcells is able to resolve the subscale features without sacrificing neither the high 
resolution of smooth features, nor the robustness at shock waves and other discontinuities. 

For invalid cells the discrete solution is then \textit{recomputed} starting from the alternative data representation $\v_h(\x,t^n)$, given by piecewise constant 
subcell averages. The update is carried out using a third order ADER-WENO finite volume scheme on the Cartesian subgrid. For details see \cite{titarevtoro,Balsara2013934,AMR3DCL}, 
where all implementation details for reconstruction and one-step time update are given. Here, we just abbreviate the entire ADER-WENO scheme acting on the subcell 
averages by 
\begin{equation}
  \v_h(\x,t^{n+1}) = \mathcal{A} \left( \v_h(\x,t^n) \right). 
	\label{eqn.aderweno} 
\end{equation} 
If the cell $T_i$ has been flagged also in the previous time step as a troubled cell, the initial data for $\v_h(\x,t^n)$ are directly available on the subgrid from the 
ADER-WENO finite volume scheme of the previous time step, otherwise, the initial data are given by the projection operator applied to the DG polynomials on the main grid.  
Appropriate boundary conditions are also needed for the WENO reconstruction on the subgrid of cell $T_i$, like in the context of high order cell-by-cell AMR schemes 
\cite{Mulet1,AMR3DCL,GrothAMR}, hence the neighbor cell data are also scattered onto virtual subcells in the same way. To summarize, the initial and boundary data for the 
subcell ADER-WENO finite volume scheme are given by 
 \begin{equation}
 \v_h(\x,t^n) = \left\{ \begin{array}{ccc} 
											  \mathcal{P}\left( \u_h(\x,t^{n}) \right) & \textnormal{ if } & \beta_i^n = 0, \\
															\mathcal{A}(\v_h(\x,t^{n-1}))      & \textnormal{ if } & \beta_i^n = 1. 
                          \end{array} \right. \qquad \x \in T_j \qquad \forall T_j \in \mathcal{V}_i. 
\label{eqn.ader.icbc} 
\end{equation} 
Note that this operation is \textit{local} and involves only the cell and its direct neighborhood $\mathcal{V}_i$, hence our subcell limiter \eqref{eqn.aderweno} together 
with the necessary initial and boundary conditions \eqref{eqn.ader.icbc} fits well in the general philosophy of DG schemes. Note also, that both, the ADER-DG scheme   
\eqref{eqn.pde.nc.gw2} as well as the ADER-WENO subcell limiter \eqref{eqn.aderweno} are \textit{one-step} schemes, hence the limiter is applied only once per time 
step, without the need of carrying out the same procedure in each Runge-Kutta substage of a classical RK-DG scheme again.  

Finally, for troubled cells ($\beta_i^{n+1} = 1$) we gather back the subgrid data representation $\v_h(\x,t^{n+1})$ produced by the subcell limiter 
using the reconstruction operator $\mathcal{R}$, which computes the final representation of the high order DG polynomial of degree $N$ on the main grid, 
i.e. $\u_h(\x,t^{n+1})=\mathcal{R}\left(\v_h(\x,t^{n+1})\right)$. This concludes the description of our extension of the MOOD paradigm to DG schemes. 

We stress that in our \aposteriori subcell limiter approach the new discrete solution is \textit{recomputed} by using a different and more robust numerical 
scheme, 
\begin{equation} 
 \u_h(\x,t^{n+1})=\mathcal{R}\left(\mathcal{A}(\v_h(\x,t^n))\right), 
\end{equation} 
while the TVB and (H)WENO limiters \textit{post-process} the unlimited candidate solution $\u^*_h(\x,t^{n+1})$ by a nonlinear \textit{reconstruction operator}, 
acting on selected degrees of freedom of the cell and its neighbors. Our limiter approach presented here is also very different from the one proposed by Sonntag and Munz 
\cite{Sonntag}, who use an \textit{a priori} switch from a DG formulation on the main grid to a subcell TVD finite volume scheme based on some \textit{a priori} indicator function. 
It is also very different from the subcell approach forwarded in \cite{CasoniHuerta2}, which smoothly switches between a high order DG scheme and a first order finite volume 
subcell method based again on an \textit{a priori} indicator function. 


\begin{figure}
  \begin{center} 
    \begin{tabular}{c} 
      \includegraphics[width=0.9\textwidth]{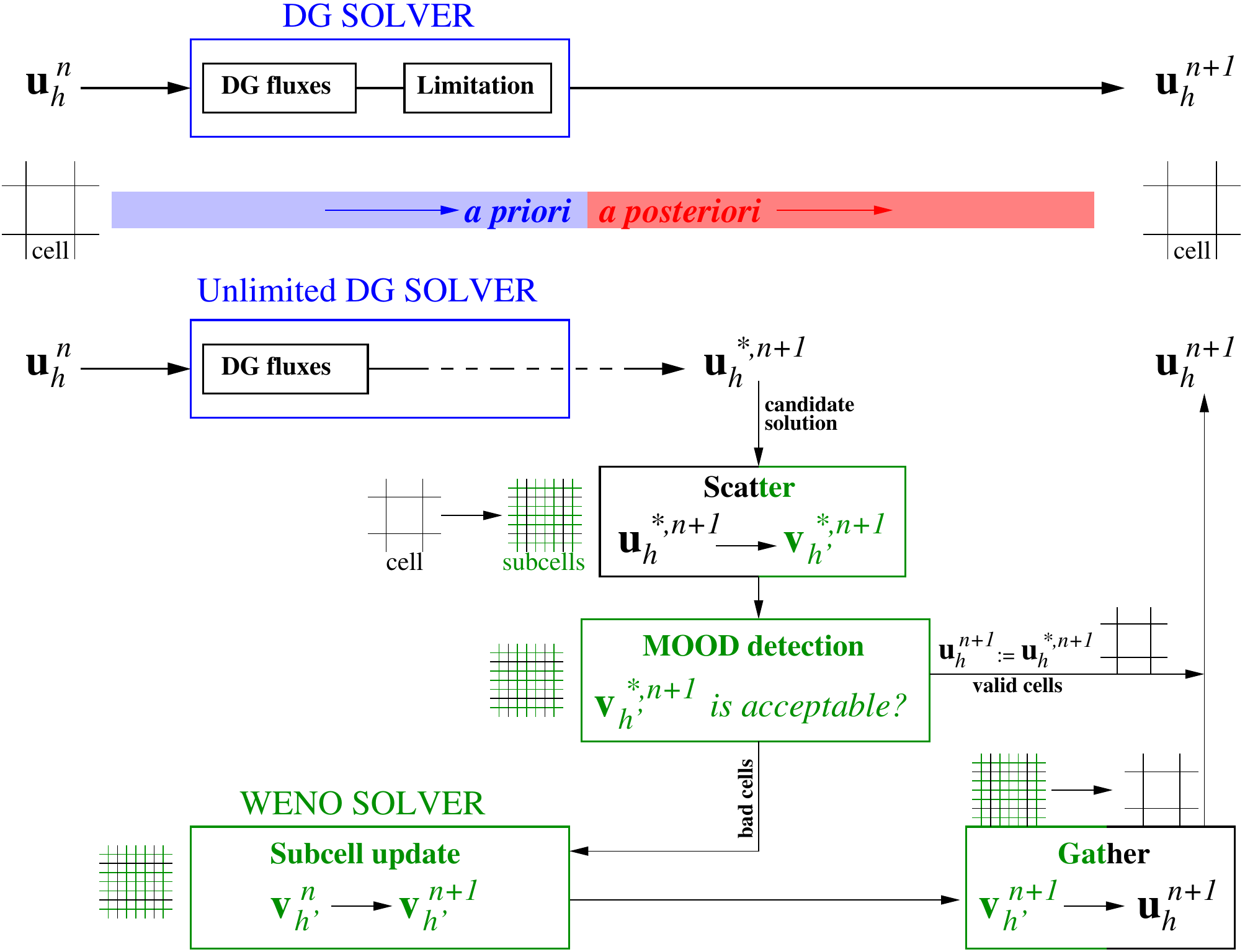} 
    \end{tabular}
    \caption{ \label{fig:algo}
      \textbf{Top:} classical DG algorithm with \apriori limiter embedded in the solver.
      After the fluxes and the limiter have been computed no more action is taken and the solution
      can not be corrected if needed anymore ---
      \textbf{Bottom:} sketch of our approach, where an unlimited candidate 
      solution $\u^*_h(\x,t^{n+1})$ is provided by the DG solver. Subsequently, the scatter step projects the 
      candidate solution onto the subgrid to get $\v^*_h(\x,t^{n+1}) = \mathcal{P}\left( \u^*_h(\x,t^{n+1}) \right)$ for 
			applying the subgrid detection criteria. 
      An \aposteriori MOOD-type procedure is further applied based on the relaxed discrete maximum principle \eqref{eqn.dmp} 
			and some physical admissibility criteria like positivity \eqref{eqn.positivity}. 
			Valid cells are left unmodified, whereas bad cells are \textit{recomputed} using an ADER-WENO finite volume scheme on the 
			subgrid. Last, for these troubled cells, a gathering step reconstructs the piecewise polynomials 
			of degree $N$ on the main grid. The green color corresponds to operations made on the subgrid, 
      black and blue colors correspond to cell based operations on the main grid.
    }
  \end{center}
\end{figure}

\subsection{Summary of the \aposteriori subcell based MOOD limiting}

The \aposteriori subcell based MOOD limiter is a five-act play, see Fig.~\ref{fig:algo}. 
We start with the cell centered piecewise polynomial representation of degree $N$  
on the main grid $\u_h(\x,t^{n}) \in \mathcal{U}_h$ in each cell and perform the following steps: 
\begin{enumerate}
\item \textit{High-order cell-based ADER-DG update.} 
  From $\u_h(\x,t^{n})$ the unlimited ADER-DG scheme \eqref{eqn.pde.nc.gw2} is used to compute a candidate 
	solution $\u^*_h(\x,t^{n+1})$. This DG scheme is meant to be the highest order accurate scheme that one 
	wishes to use.
\item \textit{Cell-to-subcell scattering (projection $\mathcal{P}$).} Project $\u^*_h(\x,t^{n+1})$ onto the subcells, 
  $\v^*_h(\x,t^{n+1})=\mathcal{P} \left( \u^*_h(\x,t^{n+1}) \right)$. 
  The subgrid of any cell must be dense enough to recover the original degrees of freedom of $\mathcal{U}_h$
	via the subcell reconstruction operator $\mathcal{R}$, hence $N_s \geq N + 1$. 
\item \textit{A posteriori MOOD detection procedure.} 
  Given a set of detection criteria like positivity \eqref{eqn.positivity} and DMP \eqref{eqn.dmp}, determine on a subcell basis 
	which cells of the main grid contain an invalid candidate solution $\u^*_h(\x,t^{n+1})$. Mark the associated cell $T_i$ as 
	invalid by setting $\beta_i^{n+1}=1$. Otherwise, if all detection criteria are satisfied on a subcell basis, the cell $T_i$ does
	not require limiting and we set the flag $\beta_i^{n+1}=0$. In this case the algorithm can exit with the candidate solution 
	$\u_h(\x,t^{n+1})$:=$\u^*_h(\x,t^{n+1})$ without further treatment. 
\item \textit{Subcell-based ADER-WENO update.} 
  For each invalid cell project the initial and boundary data onto the subcells using \eqref{eqn.ader.icbc}. 
  Then update the subgrid data using a higher order ADER-WENO finite volume scheme on the subgrid \eqref{eqn.aderweno} 
  to obtain $\v_h(\x,t^{n+1}) = \mathcal{A} \left( \v_h(\x,t^n) \right)$. Note that the initial and boundary conditions 
	for the ADER-WENO scheme are provided in each grid cell and its neighbors either by the projection operator $\mathcal{P}$ 
	(for unlimited cells), or by the ADER-WENO scheme applied in the previous time step (for limited cells). 
\item \textit{Subcell-to-cell gathering (reconstruction $\mathcal{R}$).}  
  For any troubled cell with $\beta_i^{n+1}=1$ gather the new subgrid information into a cell-centered DG polynomial of degree $N$ on the 
	main grid by applying the subcell reconstruction operator \eqref{eqn.subcell.r}, $\u_h(\x,t^{n+1})=\mathcal{R}\left( \v_h(\x,t^{n+1}) \right)$.
\end{enumerate}
Note that in this paper we do \textit{not} consider the full successive order decrementing loop proposed in the original references on MOOD schemes, 
\cite{CDL1,CDL2,CDL3,ADER_MOOD_14}, but we consider only two possible schemes: either the unlimited high order ADER-DG scheme on the coarse grid, or the 
more robust ADER-WENO finite volume scheme on the subgrid. This means that each grid zone is decremented at most once. Note also that, of course, 
the numerical fluxes in unlimited DG cells adjacent to limited cells must be recomputed in order to maintain conservation and consistency of the scheme. 
Within the DG finite element framework and the finite volume framework used here, however, it is no problem to deal with such \textit{hanging nodes} at nonconforming 
grid interfaces properly, see \cite{AMR3DCL,AMR3DNC} for details in the context of high order finite volume schemes on space-time adaptive meshes with hanging nodes in space and time.  
For the DG method, the flux integral over the element boundary in Eqn. \eqref{eqn.pde.nc.gw2} is simply written as the sum of integrals over the sub-edges between  
the main grid and the subgrid. On both grids and for both methods (ADER-DG and ADER-WENO finite volume scheme), a space-time predictor $\q_h(\x,t)$ needed for 
one-step numerical flux integration is available. 

To summarize: handling non-conforming grids with different mesh size and different orders of accuracy in each zone (so-called $hp$-refinement) is possible in a very 
natural manner in the combined DG and finite volume framework used here, see also \cite{CBS-book}. However, it means that if a cell is flagged for 
\aposteriori limiting, those direct neighbors of the cell which are not marked for limiting also need to be recomputed to ensure conservation and consistency. 
Note, however, that only direct neighbors are affected and that the same treatment  
is also necessary in the original MOOD method, \cite{CDL1,CDL2,CDL3,ADER_MOOD_14}. Furthermore, the recomputation affects only those edges of unlimited cells that are 
adjacent to a limited element, hence the simplest way of practical implementation is a flux-correction approach that subtracts the original unlimited DG fluxes and 
adds the new subcell ADER-WENO fluxes across the respective subedges to all moments of the unlimited DG cell. Of course, this correction step of unlimited zones 
adjacent to troubled zones needs additional MPI communication within a parallel implementation of the scheme, but only direct neighbors are involved, hence the scheme 
is still fully local. 

Concerning the additional memory requirements of our algorithm we would like to emphasize that the alternative data representation needs to be stored  
\textit{only} in the troubled cells, while for the detection criterion \eqref{eqn.dmp} it is sufficient to store for each cell $T_i$ the maximum and minimum 
of each conserved variable found in the subcells of the neighborhood $\mathcal{V}_i$. If the limiter is acting only in a restricted number of zones of the 
computational domain, which is usually the case, then the memory overhead produced by our approach is very small. 
Last but not least, the same a posteriori detection approach is also applied to the $L_2$ projection of the initial condition, i.e. if the DMP and the PAD are 
not satisfied for $\u_h(\x,0)$ with respect to the exact initial condition $\Q(\x,0)$, we immediately activate the subcell representation in troubled zones at 
the initial time $t=0$. 

\subsection{On the optimal choice of the subgrid size}

At a first glance, the natural choice for the number of subgrid cells seems to be $N_s=N+1$, in order to represent exactly the same amount of information 
within the subcell averages as originally contained in the high order DG polynomials on the main grid. This choice has been made for example in \cite{Sonntag}. 
However, we are convinced that this choice is \textit{not} the best one in terms of accuracy and local truncation error of the finite volume scheme on the subgrid. 
For a very detailed modified equation analysis (differential approximation) of high order ADER finite volume schemes, together with their dissipation and dispersion 
properties, see \cite{dumbser_diffapprox}. In \cite{dumbser_diffapprox} it was shown for the linear scalar advection equation in 1D that the error terms of ADER finite 
volume schemes of up to order 16 in space and time contain the factor $(1-\textnormal{CFL})$, which means that the schemes are the more accurate the larger the CFL number. 
Looking at the typical (severe) time step restriction of explicit DG schemes given by \eqref{eqn.cfl}, we note that by using $N+1$ subcells the CFL number for the finite 
volume method on the subgrid is only about \textit{half} of the maximum admissible CFL number of an explicit Godunov-type finite volume scheme, since the finite volume method 
on the subgrid must satisfy the stability condition 
\begin{equation}
  \Delta t \leq \frac{1}{d} \frac{1}{N_s} \frac{h}{|\lambda_{\max}|}.
	\label{eqn.cfl.fv} 
\end{equation} 
Note that $\frac{h}{N_s}$ is the size of the subgrid cells. If we want to use the \textit{optimal time step} on the subgrid with a Courant number close to the 
maximal one, the subgrid must satisfy $N_s = 2N+1$, which clearly follows from \eqref{eqn.cfl} and \eqref{eqn.cfl.fv}. A coarser subgrid would lead to 
\textit{small} local CFL numbers for the finite volume scheme on the subgrid and thus to \textit{more numerical dissipation and dispersion}, 
in addition to the reduced sub-cell resolution due to the coarser mesh itself! In other words, the optimal value of $N_s = 2N+1$ yields lower errors not only due to a 
smaller mesh size, but also in terms of smaller constants in front of the error terms in the local truncation error analysis. 
Of course, also the choice $N_s > 2N+1$ is suboptimal, since in this case the subgrid finite volume scheme would limit the (already small) time step of the DG 
scheme on the main grid. At the end, the choice of $N_s$ is left to the user, but the aim of this section was to discuss the choice of the \textit{optimal} value 
of $N_s$ in terms of accuracy and computational efficiency.

\section{Numerical results} \label{sec:numerics}

In this paper we focus on the the Euler equations of compressible gas dynamics
\begin{equation}
 \frac{\partial}{\partial t} \left( \begin{array}{c} \rho \\ \rho \mathbf{v} \\ \rho E \end{array} \right) +
 \nabla \cdot \left( \begin{array}{c} \rho \mathbf{v} \\ \rho \mathbf{v} \mathbf{v} + p \mathbf{I} \\ \mathbf{v} (\rho E + p )\end{array} \right) = 0,
\end{equation}
where $\rho$ denotes the mass density, $\mathbf{v}=(u,v,w)$ is the velocity vector, $p$ is the fluid pressure, $E$ is the total energy density and $\mathbf{I}$ denotes 
the $d \times d$ identity matrix. With the notation $\mathbf{v} \mathbf{v}$ we intend the dyadic product of the velocity vector with itself. 
To close the system the equation of state (EOS) of a perfect gas with adiabatic index $\gamma$ is used:
\begin{equation}
 p = (\gamma-1)\left(\rho E - \frac{1}{2}\rho \mathbf{v}^2 \right).
\end{equation}
The \aposteriori MOOD-type subcell finite volume limiter has been implemented within an MPI parallel 3D code devoted to hyperbolic system 
of conservation laws on Cartesian grids, see \cite{AMR3DCL} for the general framework.  
For most test cases presented in this section we have employed the ADER Discontinuous 
Galerkin method with piecewise polynomials of degree $N=5$ or $N=9$,  
referred to as ADER-DG-$\mathbb{P}_5$ and ADER-DG-$\mathbb{P}_9$, in the following. 
This DG method is then supplemented by the new \aposteriori sub-cell limiter (SCL). 
The scheme acting at the subcell level is a third order ADER-WENO finite volume method \cite{titarevtoro,AMR3DCL} 
with $\mathbb{P}_2$ reconstruction, denoted by WENO3. Hence the full scheme is called 
DG-$\mathbb{P}_N$+WENO3 SCL. For comparison purposes we will also employ low order 
DG schemes with $N=1,2$ in one test case. \\

In the Discontinuous Galerkin finite element framework, each variable inside a computational cell is 
represented by a polynomial of degree $N>0$. In what follows, we will therefore always use the following visualization
technique for 2D and 3D contour plots. We plot the numerical solution $\u_h(\x,t^n)$ as pointvalues interpolated 
onto a subgrid with $N_s=N+1$ if the cell is unlimited ($\beta_i^{n+1}=0$) and we will plot the alternative 
representation of the solution $\v_h(\x,t^n)$ on the subgrid used for the limiter with $N_s=2N+1$ if the cell has 
been detected as a troubled cell ($\beta_i^{n+1}=1$). For 1D cuts we usually take equidistant samples at
subgrid level of the solution representations $\u_h(\x,t^n)$ and $\v_h(\x,t^n)$, respectively. 
This is extremely important in order to verify that the subscale structure of the DG polynomials 
really represents a physically valid state within one large cell. Moreover, an acceptable DG method should 
maintain smooth solutions and, more importantly, it should produce genuinely discontinuous profiles  
at shock waves without spurious oscillations.
In order to visualize which cells have been limited, in the rest of this section we will systematically represent
in blue the unlimited cells ($\beta_i^{n+1}=0$) and in red ($\beta_i^{n+1}=1$) the limited cells,  
see for instance Figure~\ref{fig:sodlax3D} below. \\
A well known difficulty of a high order accurate method dealing with discontinuous solutions is to deposit entropy 
(and dissipation) on a length scale which is much smaller than the characteristics length of the coarse cell. 
Indeed, if the coarse cell is very large, spreading a shock wave over one or two of such cell(s) typically generates 
an excessively large numerical dissipation. Our \aposteriori subcell limiter is supposed to act differently and  
in this section we will provide the numerical evidence of this property.\\
The methodology of validation is based on a sequence of classical test cases, namely:
\begin{itemize}
\item Sod and Lax shock tube \cite{toro-book} - these shock tubes are classical tests to
  assess the ability of a numerical method to deal with simple waves (rarefaction,
  contact discontinuity and shock wave).
\item Smooth vortex - this test is designed to observe the high-order of 
  accuracy of a numerical method when the exact solution is smooth. Note that,
  in theory, in this case any limiter should not be activated at all if the mesh is fine enough.
  A detailed numerical convergence table for polynomial degrees ranging from $N=1$ to $N=9$ is provided. 
\item Shu-Osher oscillatory shock tube \cite{shuosher2} -
  this test is designed to show 
  the difficulty of capturing smooth small-scale features and shock waves.
\item Double Mach reflection \cite{woodwardcol84} - 
  this classical test is meant to measure visually
  the ability of a numerical method to capture complex patterns created after
  the interaction of shock waves. We use this test to show the behavior of the 
  method when the polynomial degree of the DG scheme is increased.
\item Forward facing step \cite{woodwardcol84} - 
  this classical test simulates a supersonic 
  flow over a forward facing step. The solution approaches a
  steady-state solution composed of multiple interacting shock waves and vortex like 
  structures which are often dissipated with low order accurate schemes.
\item 2D Riemann problems \cite{kurganovtadmor} - 
  this classical suite of test problems is meant to assess 
  the ability of a numerical method to solve genuinely two-dimensional Riemann problems emerging
  from four piecewise constant states joining at the origin.
\item Shock-vortex interaction \cite{Shock_vortex_95,caishu93} - 
  this test is designed to observe the 
  interaction of a planar shock wave with a cylindrical vortex and the subsequent complex
  structures of primary and secondary waves.
\item 3D explosion problem \cite{toro-book} - this test is designed to show the
  ability of a scheme to deal with separate spherical waves in 3D and subsequently validates
  the approach in 3D. We use this test also to prove that our approach is well suited for being
	used within a massively parallel MPI framework running on 8000 CPU cores and capable of 
	dealing with 10 billion space-time degrees of freedom per time step and conserved variable.
\end{itemize}

\subsection{Sod and Lax shock tube}  \label{ssec:Sod_Lax}
\begin{figure}
  \begin{center} 
  \begin{tabular}{cc} 
    \includegraphics[width=0.47\textwidth]{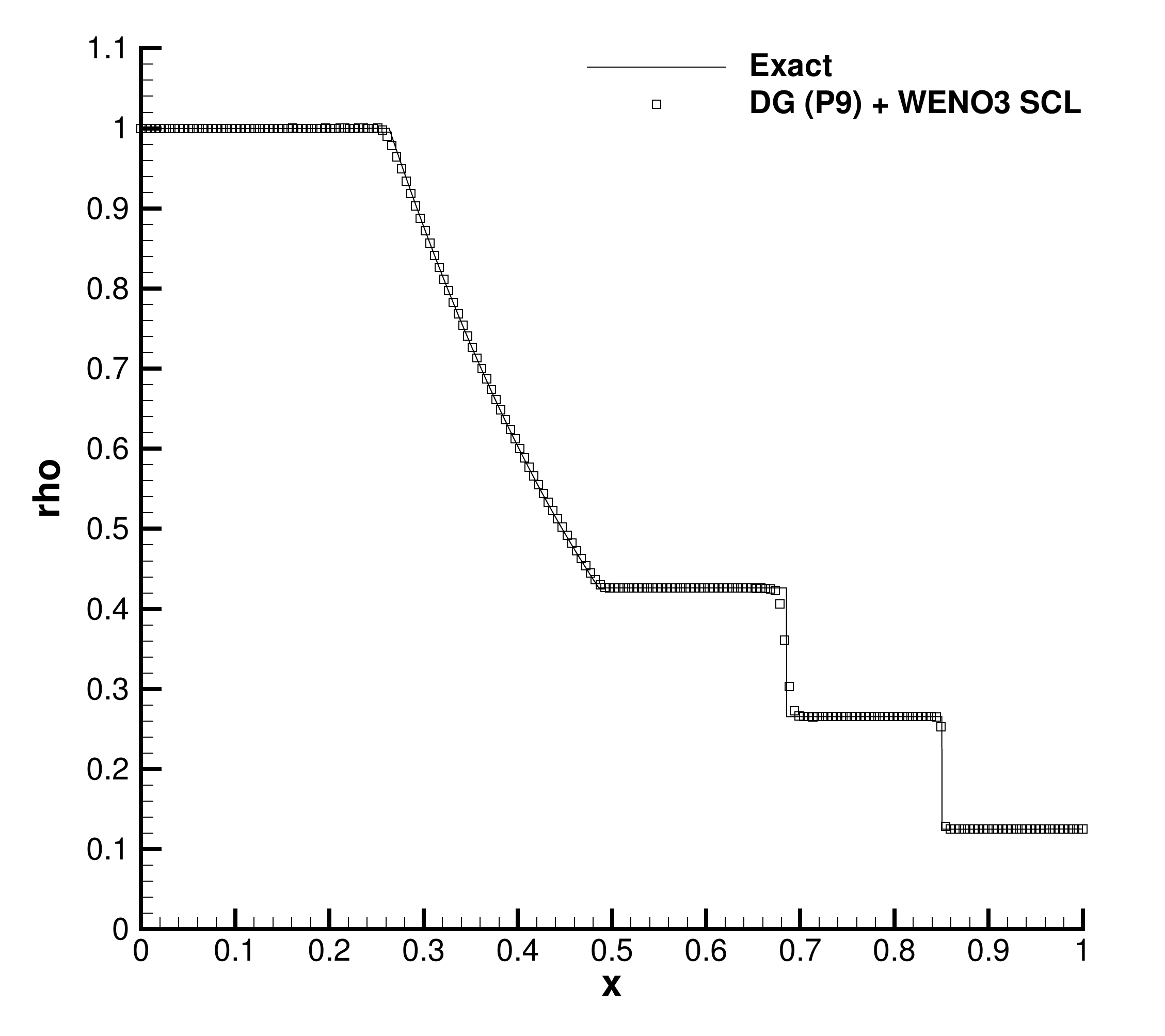}&
    \includegraphics[width=0.47\textwidth]{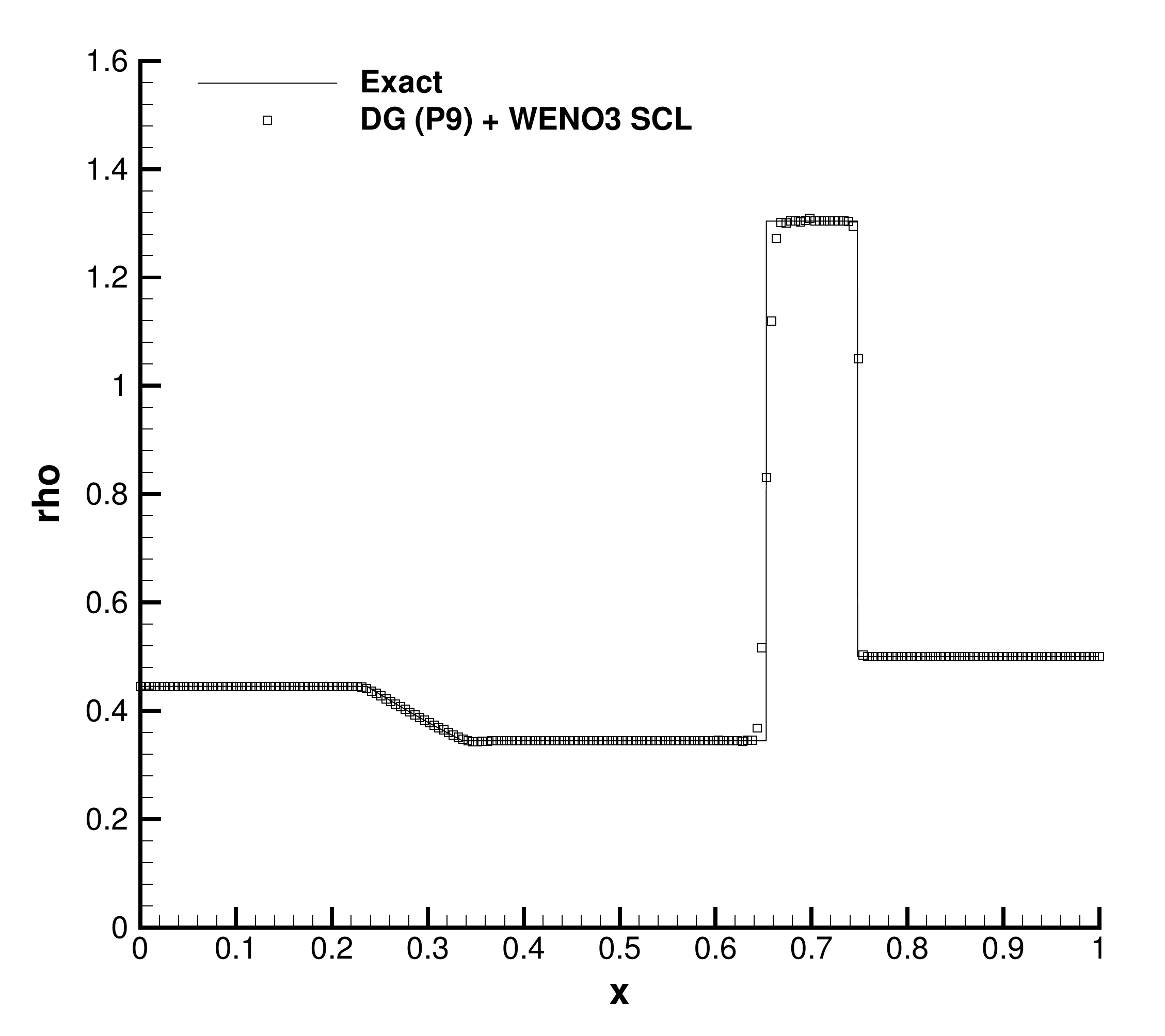}
    \\
    \includegraphics[width=0.47\textwidth]{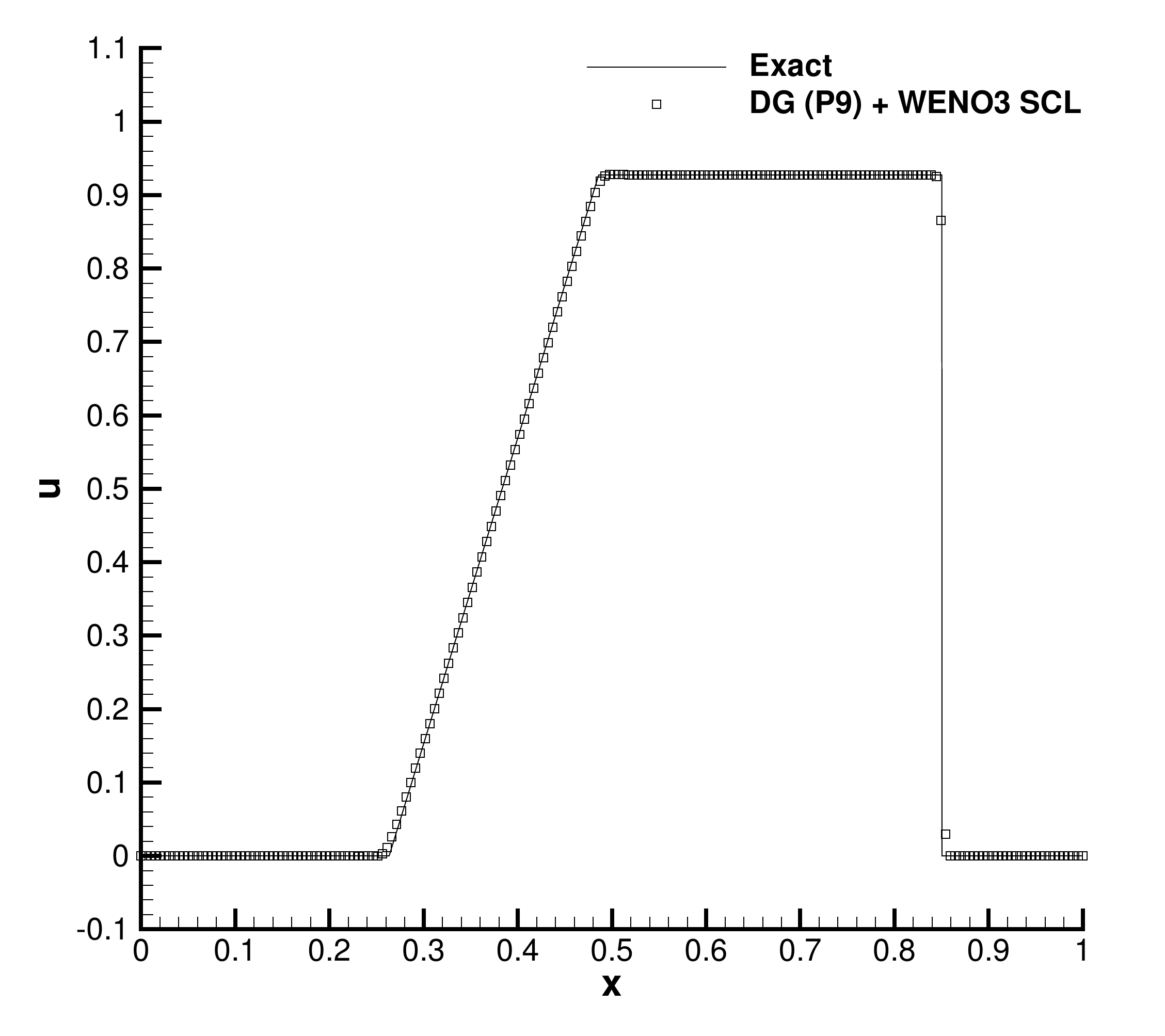}  &
    \includegraphics[width=0.47\textwidth]{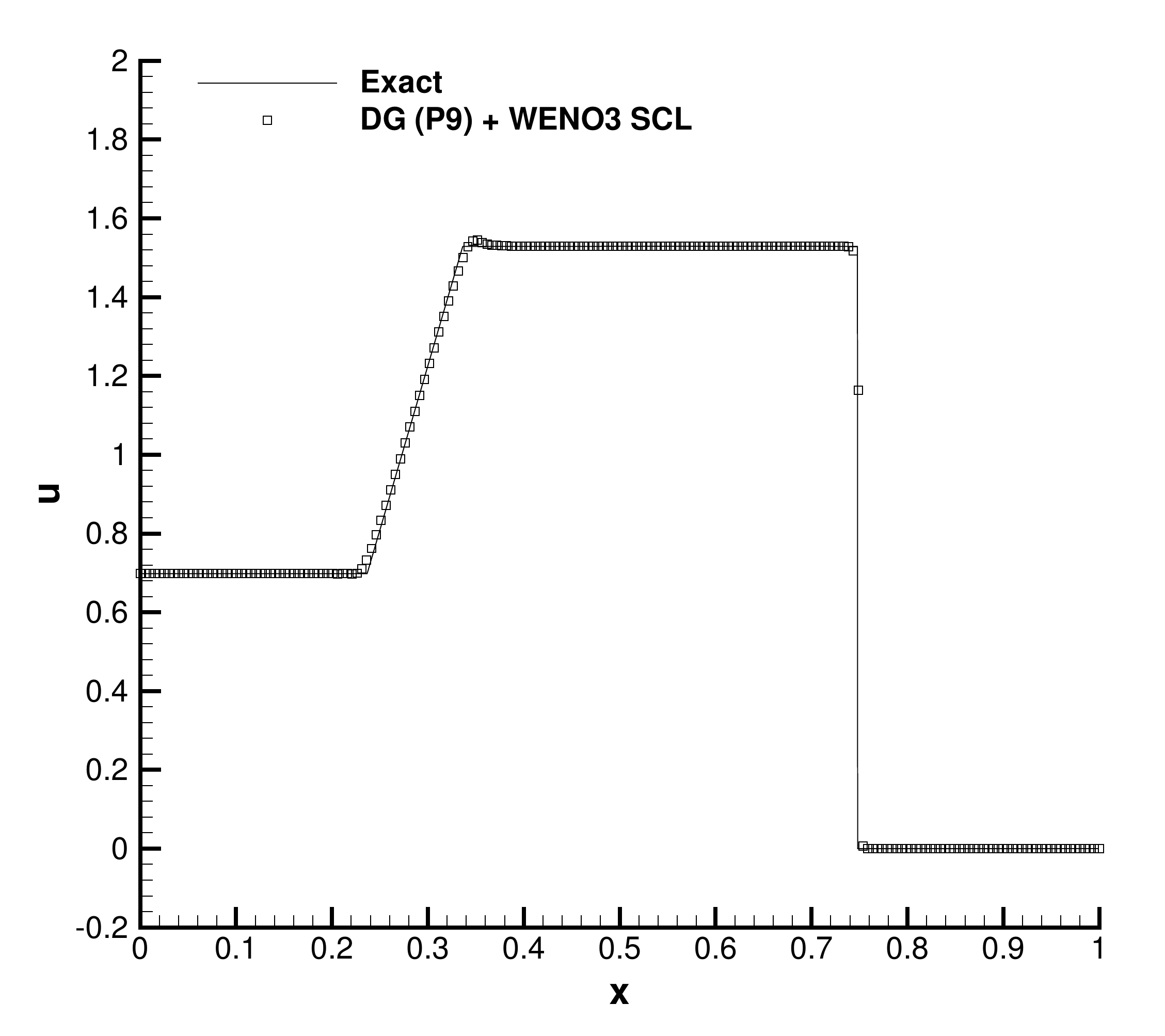} 
    \\
    \includegraphics[width=0.47\textwidth]{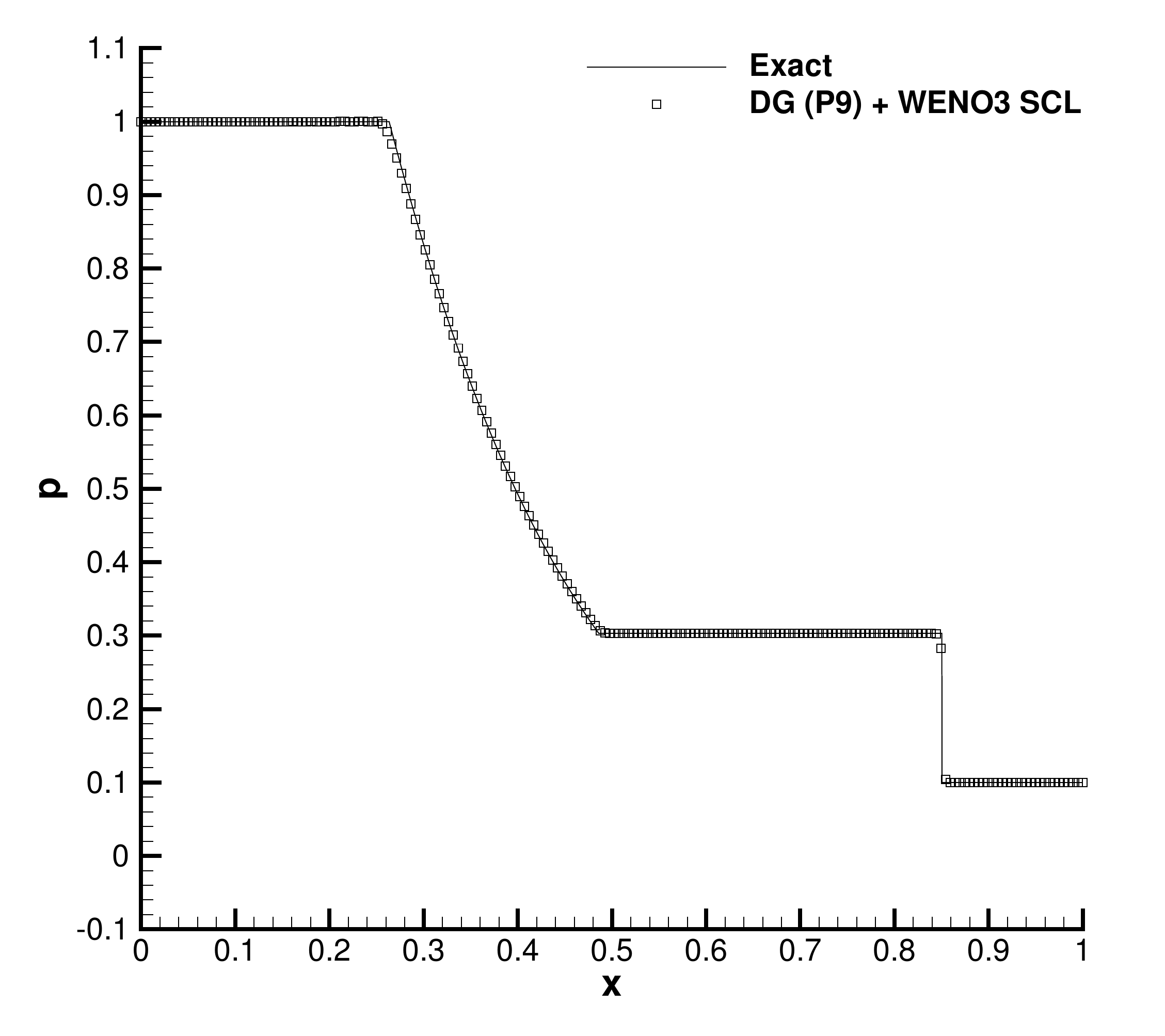}  &
    \includegraphics[width=0.47\textwidth]{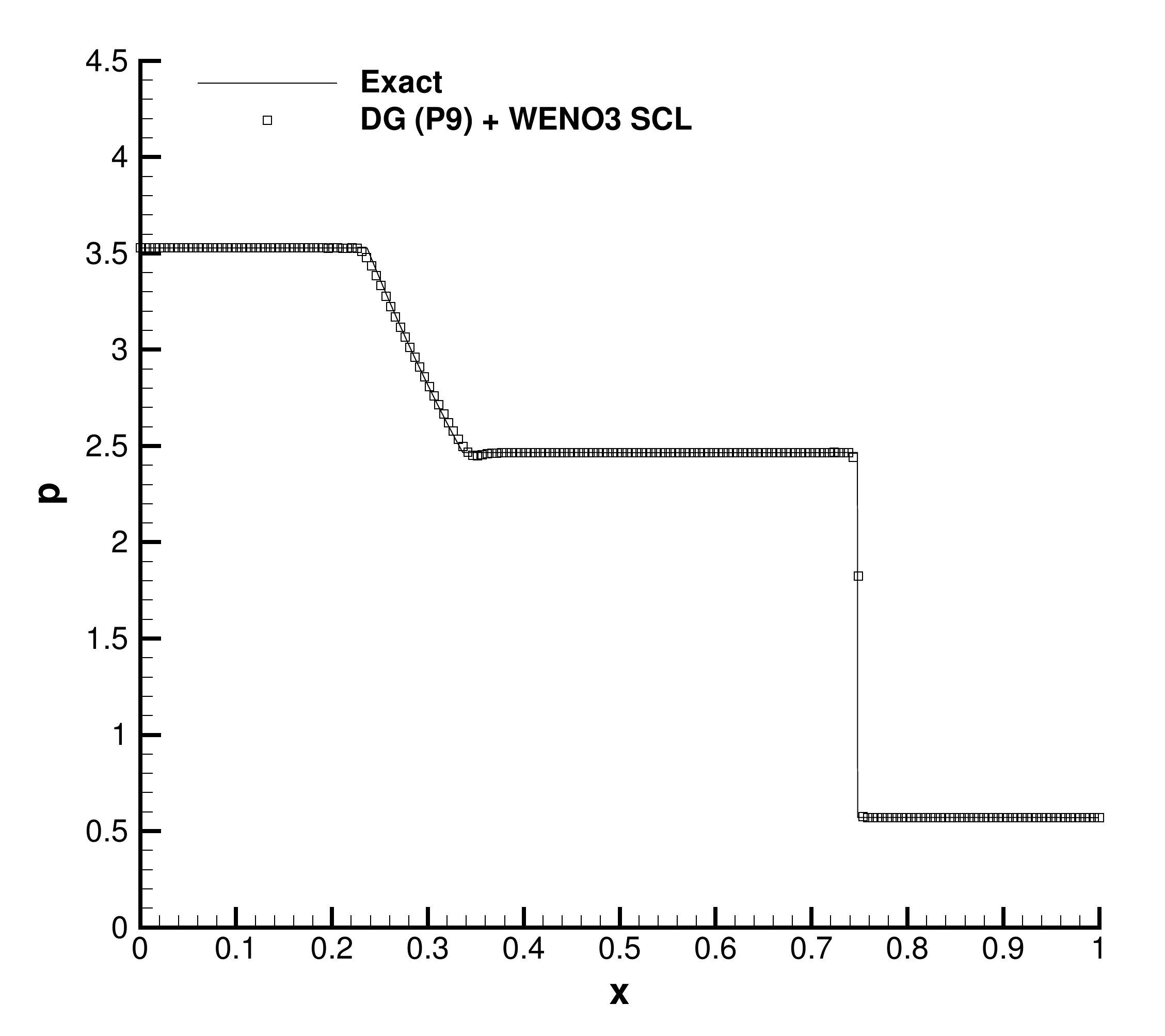}    
  \end{tabular}
   \caption{ \label{fig:sodlax1D}
			Sod shock tube problem (left panels) at $t_{\text{final}}=0.2$ 
      and Lax problem (right panels) at $t_{\text{final}}=0.14$. In both cases a very coarse mesh of only 
			$20 \times 5$ cells on the main grid has been used. An ADER-DG-$\mathbb{P}_9$ scheme supplemented with 
			\aposteriori ADER-WENO3 subcell limiter has been used --- 1D cut on 200 equidistant sample points through the 
			numerical solution (symbols) \textit{vs} exact solution for density (top), velocity component $u$ (middle) 
			and pressure (bottom).  
     }
  \end{center}
\end{figure}

\begin{figure}
  \vspace{-2cm}
  \begin{center} 
    \includegraphics[width=0.8\textwidth]{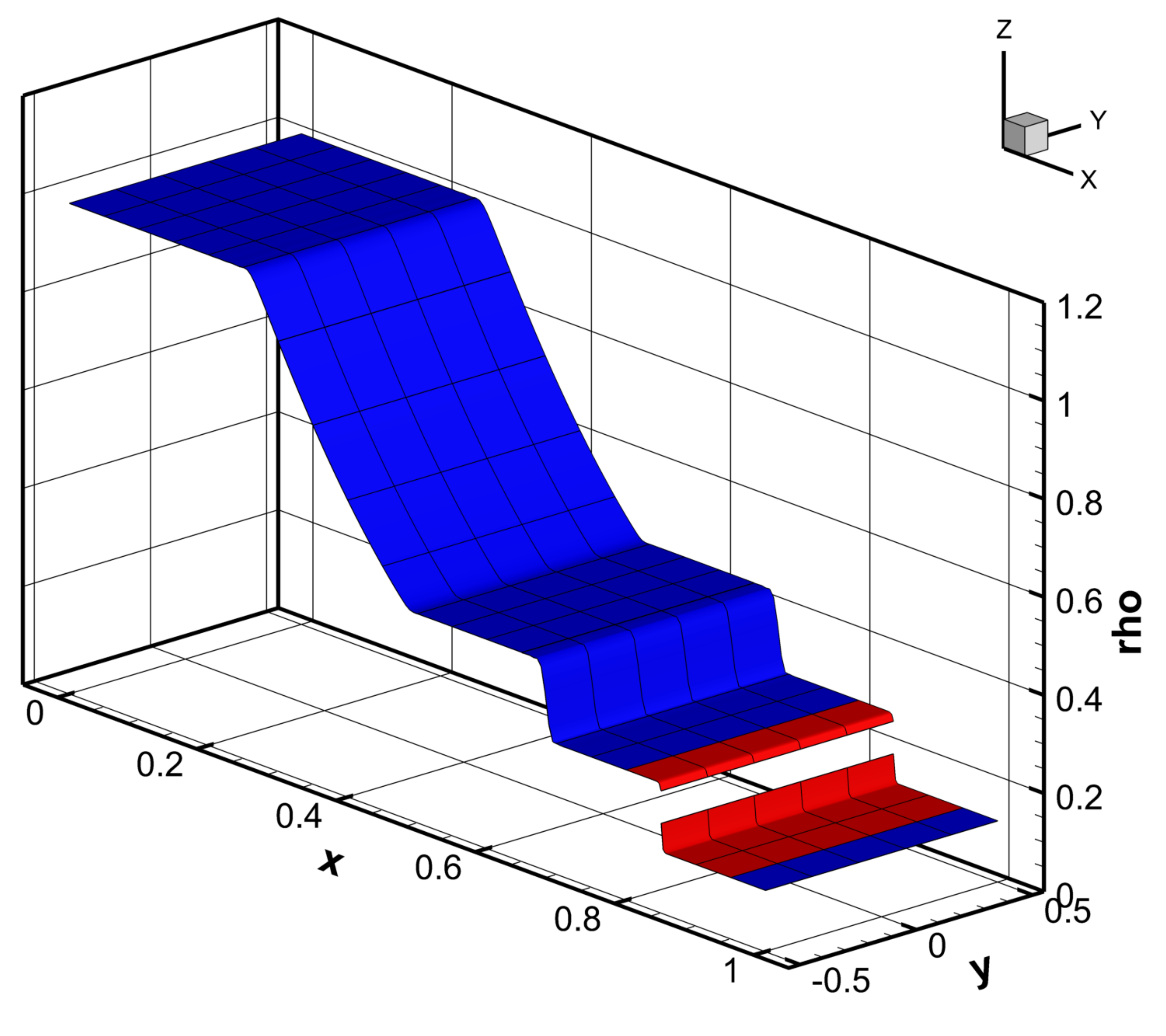}   \\
    \vspace{-0.4cm}
    \includegraphics[width=0.8\textwidth]{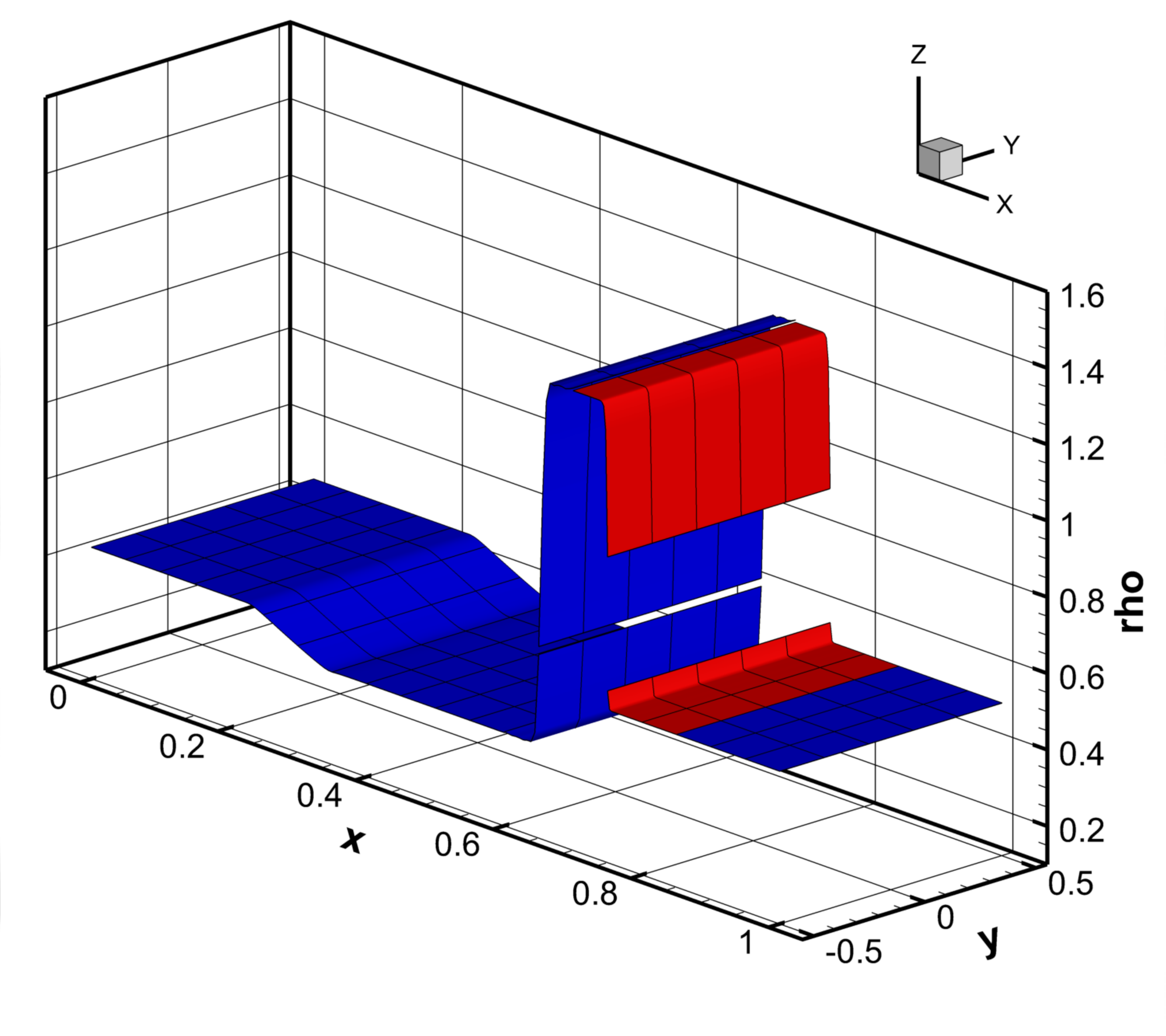} 
    \vspace{-0.4cm} 
   \caption{ \label{fig:sodlax3D}
      Sod problem (top panel $t_{\text{final}}=0.2$) 
      and Lax problem (bottom panel $t_{\text{final}}=0.14$)  
      solved on a $20\times 5$ element mesh using ADER-DG-$\mathbb{P}_9$ with \aposteriori WENO3 subcell limiter --- 
      The density variable is displayed. Troubled cells are shown in red, while blue cells are updated with the unlimited 
			ADER-DG-$\mathbb{P}_9$ on the main grid. }
  \end{center}
\end{figure}

Here, we run the planar Sod shock tube problem and the classical Lax shock 
tube problem
on a 2D structured mesh to assess the ability of the methods to capture 
one-dimensional
simple waves. The initial conditions for density, velocity component $u$ and 
pressure are listed
in Table~\ref{tab.rp.ic}. The other velocity component is initialized with $v=0$.  
\begin{table}[!htbp] 
 \numerikNine
  \begin{center}
    \begin{tabular}{|l||ccc|ccc|c|}
      \hline      
      \multirow{2}{*}{\textbf{\textsl{Problem}}} & \multicolumn{3}{c|}{\textbf{Left state}} & 
      \multicolumn{3}{c|}{\textbf{Right state}} & \textbf{Final time}\\
       & $\rho_L$ & $u_L$ & $p_L$ & $\rho_R$ & $u_R$ & $p_R$ & $t_{\textnormal{final}}$ \\
      \hline
      \hline
      {Sod}     &  1.0     & 0.0   & 1.0   & 0.125    & 0.0   & 0.1   & 0.2 \\
      \hline
      {Lax}     &  0.445   & 0.698 & 3.528 & 0.5      & 0.0   & 0.571 & 0.14 \\
      \hline
    \end{tabular}
  \end{center}
  \caption{ \label{tab.rp.ic}
    Initial left and right states for the density $\rho$, velocity $u$
    and the pressure $p$ for the Sod and Lax shock tube problems. Final simulation times
    $t_{\text{final}}$ are also provided.  }
\end{table}
The ratio of specific heats is $\gamma=1.4$ and for both problems
the initial discontinuity is located in $x=0.5$. 
The exact solution for these one-dimensional Riemann problems can be found in 
\cite{toro-book}.
The computational domain $\Omega=[0;1]\times [-0.5;0.5]$ is
paved with a very coarse structured mesh made of $N_x\times N_y$ cells with $N_x=20$ and $N_y=5$,
 see Figure~\ref{fig:sodlax3D}.
Dirichlet boundary conditions are imposed in $x$-direction while periodic boundaries are applied in $y$-direction. \\
In Figure~\ref{fig:sodlax1D} we present the results for 
DG-$\mathbb{P}_9$ where density, velocity in $x$-direction and pressure 
are displayed along the $x$ direction versus the exact solution.
In this figure each DG-$\mathbb{P}_9$ polynomial is represented by $10$ sample points. 
We observe a very good agreement with
the exact solution, the shock being resolved just in one cell. 
In Figure~\ref{fig:sodlax3D} we also represent the unlimited polynomials (blue) 
and the limited ones (red), we can clearly see that the contact discontinuity 
in the Sod problem is resolved \textit{within} one single cell in an almost S-type shape,  
thanks to the use of very high order polynomials of degree $N=9$.  
On the other hand, although the cells embracing the shock region have been limited 
(they are colored in red), the subcell data are first of all non-oscillatory, second, 
the plateaus before and after the shock are well captured and third the shock is properly 
represented by a jump at an element interface. Back to the 1D panels we can see that this jump on 
subcell values is perfectly located at the exact shock wave location. 
The same comments also hold for the Lax problem. In conclusion, although the number of cells is very 
small, the ADER-DG-$\mathbb{P}_9$ scheme with SCL can produce a very accurate solution and also the discontinuities 
are very sharp even on a very coarse mesh. 
The \aposteriori limiter has acted properly to maintain the shock wave within one or two \textit{subcells}.

\subsection{Smooth vortex}  \label{ssec:Vortex}
The isentropic vortex problem was initially introduced for the two-dimensional
compressible Euler equations in \cite{Shu1} to test the accuracy of numerical 
methods, since the exact solution is smooth and has a simple analytical expression. 
Let us consider the computational domain $\Omega = [-5,5]\times [-5,5]$
and an ambient flow characterized by
$\rho_\infty=1.0$, $u_\infty=1.0$, $v_\infty=1.0$, $w_\infty=0.0$, $p_\infty=1.0$,
with a normalized ambient temperature $ T^*_\infty =1.0$
computed with the perfect gas equation of state and $\gamma=1.4$.\\
A vortex is centered on the $z$ axis line at
$(x_{\text{vortex}}, y_{\text{vortex}})=(0,0)$
and supplemented to the ambient gas at the initial time $t=0$ with the following conditions
$u = u_\infty + \delta u$, $ v = v_\infty + \delta v$, $w=w_\infty$, $ T^* = T^*_\infty + \delta T^*$ where
\begin{eqnarray}
  \nonumber
  \delta u   = -y' {\frac {\beta} {2 \pi}} \exp \left( {\frac {1-r^2} {2}} \right), \; \;
  \delta v   = x' {\frac {\beta} {2 \pi}} \exp \left( {\frac {1-r^2} {2}} \right), \; \;
  \delta T^* = - { \frac {(\gamma - 1 ) \beta} {8 \gamma \pi^2}} \exp \left( {1-r^2} \right),
\end{eqnarray}
with $ r = \sqrt{{x'}^2 + {y'}^2}$ and $x' = x - x_{\text{vortex}}, y' = y - y_{\text{vortex}}$.
The vortex strength is given by $\beta=5.0$ and the initial density follows the relation
\begin{eqnarray}
\rho
       = \rho_\infty \left( {\frac{T^*}{T_\infty^*} } \right)^{\frac{1}{\gamma-1} }
       = \left(1 - { \frac {(\gamma - 1 ) \beta} {8 \gamma \pi^2}} \exp \left( {1-r^2} \right)
            \right)^{\frac{1}{\gamma-1} }.
\end{eqnarray}
Periodic boundary conditions are prescribed everywhere, so that at the final time $t_{\text{final}}=10$ 
the vortex is back to its original position. The numerical flux used here was the Osher-type flux presented 
in \cite{OsherUniversal}. This problem has a smooth solution and thus 
should be simulated with effective high-order of accuracy if the limiter behaves properly.
Four successively refined grids made of $N_x\times N_x$ squares are employed to compute these
errors. We compute the discrete $L^1$, $L^2$ and $L^\infty$ error norms between the exact 
solution and the numerical solution for the density at the final time. The computation of these 
error norms is performed using a sufficiently high order accurate Gaussian quadrature rule. 
The errors and the rate of convergence are reported in Table~\ref{tab:Vortex_Error} 
for the ADER-DG-$\mathbb{P}_N$ supplemented with the \aposteriori WENO3 subcell limiter and 
$N$ varying from $2$ to $9$. 
From this table we can observe that the optimal order of convergence is essentially achieved 
by the scheme and that the proposed \aposteriori subcell limiter does not destroy the 
accuracy of the high order DG scheme. 
Beyond $\mathbb{P}_7$ it seems more difficult to get perfect orders due to the fact that the
errors are very small and roundoff starts to play a role. However, even on these ultra coarse 
meshes the results are very accurate, for instance comparing the ADER-DG-$\mathbb{P}_9$ on $8\times 8$ 
cells which is roughly as accurate as ADER-DG-$\mathbb{P}_3$ on $50\times 50$ cells  
or ADER-DG-$\mathbb{P}_2$ on $100\times 100$ elements. 
%
 \begin{table}[!htbp] 
 \centering
 \numerikNine
 \begin{tabular}{|c|c||ccc|ccc|c|}
   \hline
   \multicolumn{9}{|c|}{\textbf{2D isentropic vortex problem --- ADER-DG-$\mathbb{P}_N$ + WENO3 SCL}} \\
   \hline
   \hline
    & $N_x$ & $L^1$ error & $L^2$ error & $L^\infty$ error & $L^1$ order & $L^2$ order & $L^\infty$ order &
   Theor. \\
   \hline
   \hline
   \multirow{4}{*}{\rotatebox{90}{{DG-$\mathbb{P}_2$}}}
& 25	& 9.33E-03	& 2.07E-03	& 2.02E-03	&---	&---	&--- & \multirow{4}{*}{3}\\
& 50	& 6.70E-04	& 1.58E-04	& 1.66E-04	& 3.80	& 3.71	& 3.60 &\\
& 75	& 1.67E-04	& 4.07E-05	& 4.45E-05	& 3.43	& 3.35	& 3.25 &\\
& 100	& 6.74E-05	& 1.64E-05	& 1.82E-05	& 3.15	& 3.15	& 3.10 &\\
   \cline{2-8}
   \hline
   \multirow{4}{*}{\rotatebox{90}{{DG-$\mathbb{P}_3$}}}
& 25	& 5.77E-04	& 9.42E-05	& 7.84E-05	& ---	& ---	& --- & \multirow{4}{*}{4}\\
& 50	& 2.75E-05	& 4.52E-06	& 4.09E-06	& 4.39	& 4.38	& 4.26 &\\
& 75	& 4.36E-06	& 7.89E-07	& 7.55E-07	& 4.55	& 4.30	& 4.17 &\\
& 100	& 1.21E-06	& 2.37E-07	& 2.38E-07	& 4.46	& 4.17	& 4.01 &\\
   \cline{2-8}
   \hline
   \multirow{4}{*}{\rotatebox{90}{{DG-$\mathbb{P}_4$}}}
& 20	& 1.54E-04	& 2.18E-05	& 2.20E-05	& ---	& ---	& --- & \multirow{4}{*}{5}\\		
& 30	& 1.79E-05	& 2.46E-06	& 2.13E-06	& 5.32	& 5.37	& 5.75 &\\
& 40	& 3.79E-06	& 5.35E-07	& 5.18E-07	& 5.39	& 5.31	& 4.92 &\\
& 50	& 1.11E-06	& 1.61E-07	& 1.46E-07	& 5.50	& 5.39	& 5.69 &\\
   \cline{2-8}
   \hline
   \multirow{4}{*}{\rotatebox{90}{{DG-$\mathbb{P}_5$}}}
& 10	& 9.72E-04	& 1.59E-04	& 2.00E-04	& ---	& ---	& --- & \multirow{4}{*}{6}\\		
& 20	& 1.56E-05	& 2.13E-06	& 2.14E-06	& 5.96	& 6.22	& 6.55 &\\
& 30	& 1.14E-06	& 1.64E-07	& 1.91E-07	& 6.45	& 6.33	& 5.96 &\\
& 40	& 2.17E-07	& 2.97E-08	& 3.59E-08	& 5.77	& 5.93	& 5.82 &\\
   \cline{2-8}
   \hline
   \multirow{4}{*}{\rotatebox{90}{{DG-$\mathbb{P}_6$}}}
& 5	& 2.24E-02	& 4.15E-03	& 3.11E-03	& ---	& ---	& --- & \multirow{4}{*}{7}\\		
& 10	& 1.76E-04	& 2.75E-05	& 2.86E-05	& 6.99	& 7.24	& 6.76 &\\
& 20	& 1.67E-06	& 2.28E-07	& 2.26E-07	& 6.72	& 6.91	& 6.98 &\\
& 25	& 3.60E-07	& 4.96E-08	& 6.27E-08	& 6.86	& 6.84	& 5.74 &\\
   \cline{2-8}
   \hline
   \multirow{4}{*}{\rotatebox{90}{{DG-$\mathbb{P}_7$}}}
& 5	& 5.50E-03	& 1.22E-03	& 1.46E-03	& ---	& ---	& --- & \multirow{4}{*}{8}\\		
& 10	& 4.63E-05	& 6.26E-06	& 6.95E-06	& 6.89	& 7.61	& 7.71 &\\
& 15	& 1.62E-06	& 2.20E-07	& 2.29E-07	& 8.28	& 8.26	& 8.42 &\\
& 20	& 2.05E-07	& 2.80E-08	& 2.28E-08	& 7.18	& 7.17	& 8.01 &\\
   \cline{2-8}
   \hline
   \multirow{4}{*}{\rotatebox{90}{{DG-$\mathbb{P}_8$}}}
& 4	& 9.11E-03	& 1.80E-03	& 3.44E-03	& ---	& ---	& --- & \multirow{4}{*}{9}\\		
& 8	& 4.97E-05	& 7.51E-06	& 6.93E-06	& 7.52	& 7.90	& 8.96 &\\
& 10	& 7.50E-06	& 1.05E-06	& 1.18E-06	& 8.47	& 8.81	& 7.95 &\\
& 15	& 2.40E-07	& 3.34E-08	& 3.09E-08	& 8.49	& 8.51	& 8.98 &\\
   \cline{2-8}
   \hline
   \multirow{4}{*}{\rotatebox{90}{{DG-$\mathbb{P}_9$}}}
& 4	& 3.95E-03	& 7.89E-04	& 1.42E-03	& ---	& ---	& --- & \multirow{4}{*}{10}\\		
& 8	& 1.01E-05	& 1.44E-06	& 1.52E-06	& 8.61	& 9.09	& 9.87 &\\
& 10	& 1.44E-06	& 2.00E-07	& 2.27E-07	& 8.74	& 8.85	& 8.51 &\\
& 12	& 2.67E-07	& 3.70E-08	& 3.77E-08	& 9.26	& 9.25	& 9.85 &\\
   \cline{2-8}
   \hline
 \end{tabular}
 \caption{ \label{tab:Vortex_Error} $L^1, L^2$ and $L^\infty$ errors and convergence rates for the 
   2D isentropic vortex problem for the ADER-DG-$\mathbb{P}_N$ scheme supplemented
   with \aposteriori ADER-WENO3 subcell limiter for $N=2$ to $9$ from top to bottom.}
 \end{table}

\subsection{Shu-Osher oscillatory shock tube}  \label{ssec:Shu_Osher}
\begin{figure}
  \begin{center} 
  \begin{tabular}{c} 
    \includegraphics[width=0.7\textwidth]{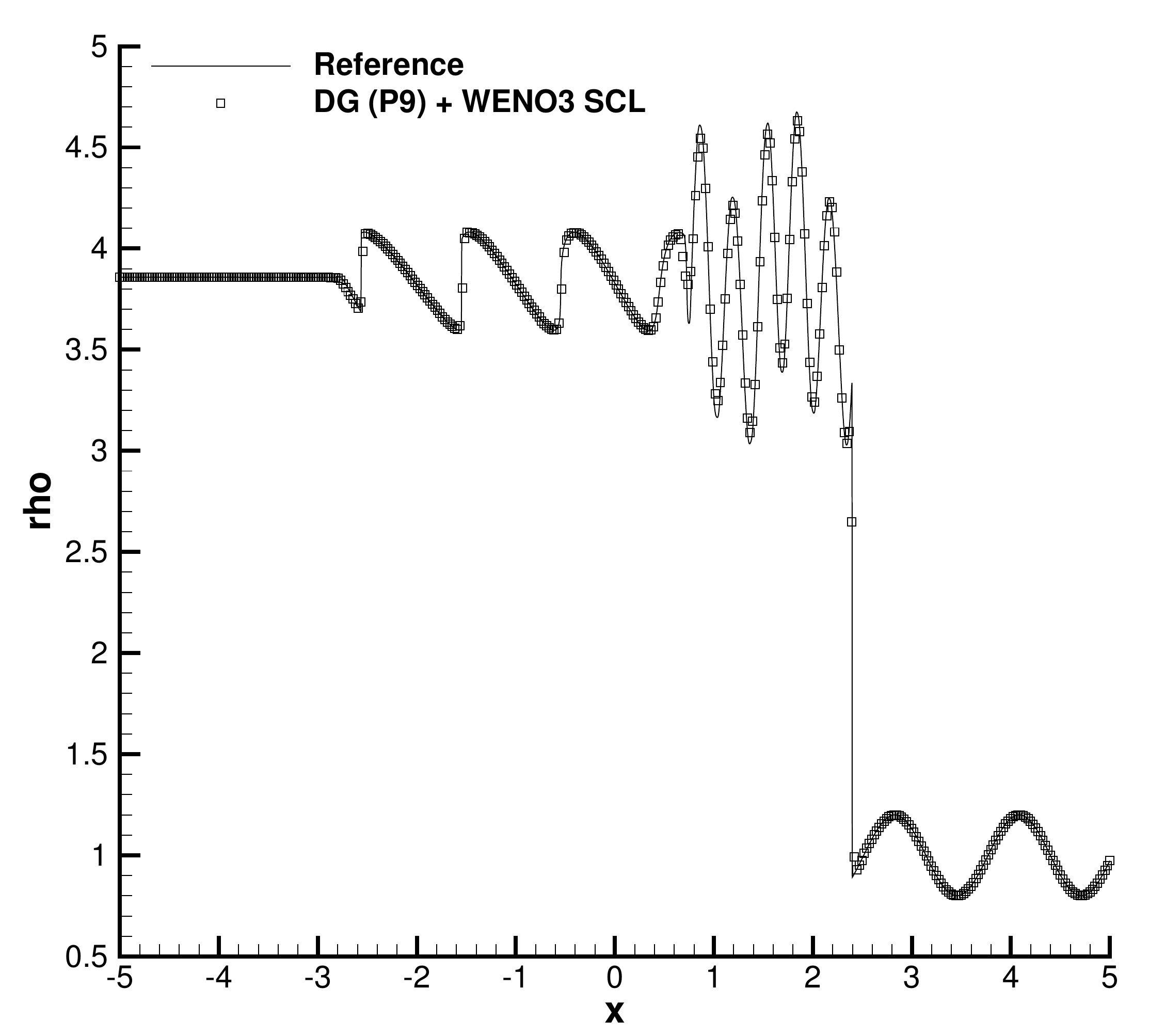}  \\
    \includegraphics[width=0.8\textwidth]{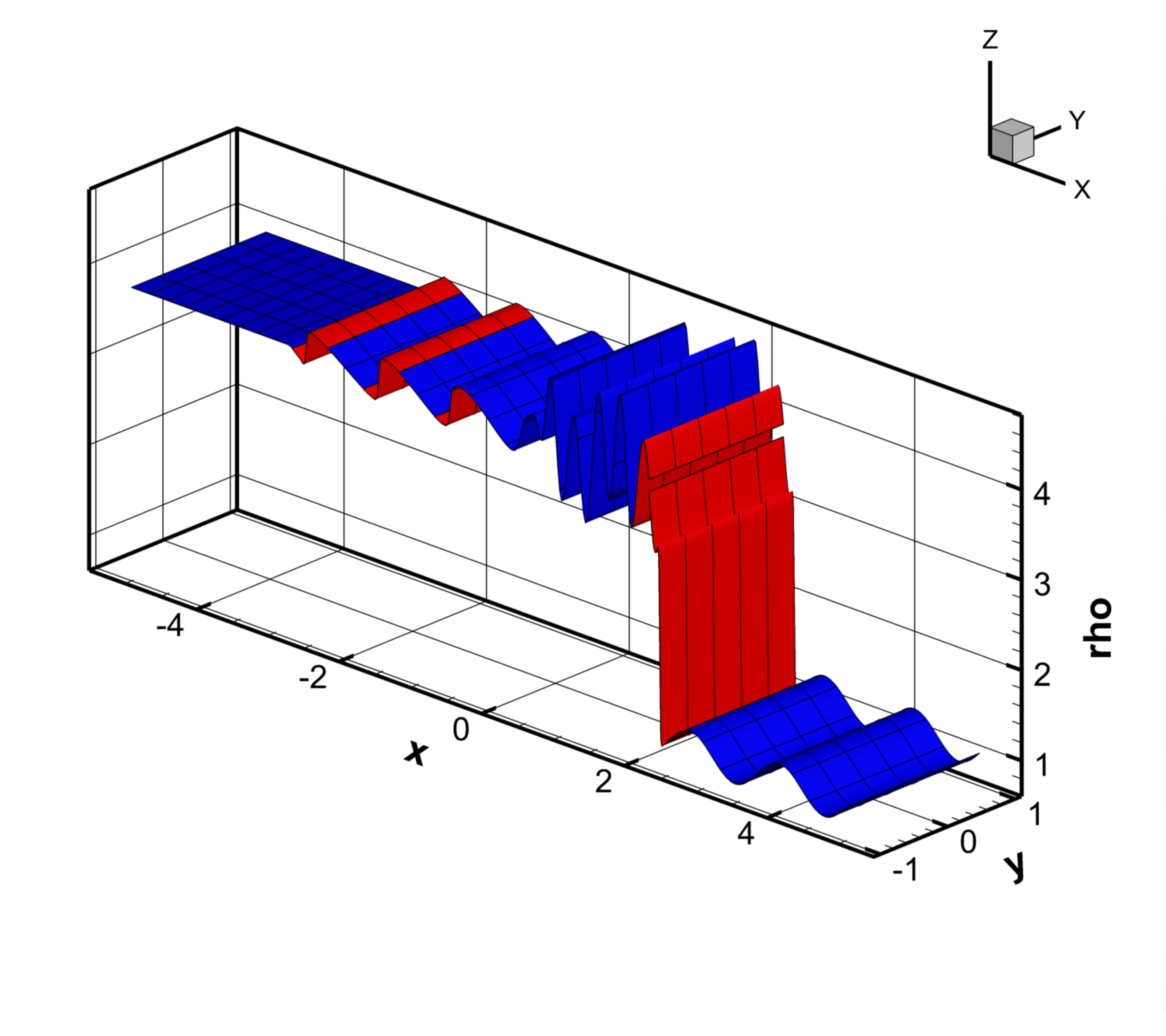}  
  \end{tabular}
  \vspace{-1cm}
   \caption{ \label{fig:ShuOsher}
      Shu-Osher problem at $t_{\text{final}}=0.18$ with ADER-DG-$\mathbb{P}_9$ 
      with \aposteriori ADER-WENO3 subcell limiter on a $40\times 5$ mesh ---
      Top: 1D cut through the numerical solution (symbols) \textit{vs} reference solution 
      (ultra fine ADER-WENO solution in straight line). Any $\mathbb{P}_9$ polynomial 
      is represented by $10$ sample points per cell ---
      Bottom: troubled cells, which have been updated with ADER-WENO3 on the subgrid, are shown in red, blue cells 
      have used the unlimited ADER-DG-$\mathbb{P}_9$ scheme on the main grid.
     }
  \end{center}
\end{figure}

This test \cite{shuosher2} is a 1D hydrodynamic shock tube.
The downstream flow has a sinusoidal density fluctuation 
$\rho=1-\varepsilon \sin(\lambda \pi x)$ with a wave length of $\lambda = 5$ 
and an amplitude of $\varepsilon = 0.2$.  A Mach $3$ shock front 
is initially located at $x=-4$ on domain $[-5;5]$.
The left and the right states are given by
$\rho_L=3.857143$, $u_L= 2.629369$, $p_L = 10.33333$ and
$\rho_R=1 + 0.2 \sin(5\pi x)$, $u_R=0$ and $p_R=1$. The final time is set 
to $t_{\text{final}}=0.18$.
This problem involves small scales after the shock has interacted with the sine
wave that can be captured either with a fine enough mesh or with high order accurate method.
Here a very coarse mesh made of $40$ cells in $x$-direction and $5$ in $y$-direction 
is chosen. The numerical flux was the Osher-type flux presented in \cite{OsherUniversal}. \\
The results 
of ADER-DG-$\mathbb{P}_9$ with \aposteriori WENO3 subcell limiter method are
presented in Figure~\ref{fig:ShuOsher}. A 1D view is first depicted, where each 
$\mathbb{P}_9$ polynomial is represented by $10$ sample points per cell.
A reference solution (straight line) obtained by a classical third order ADER-WENO 
scheme on a very fine mesh is also plotted. The quality of the result is excellent 
for a $40$ cell mesh. 
 The bottom part of the figure displays the cell centered $\mathbb{P}_9$ polynomials.
The red cells have been limited, therefore updated with WENO3 on subcells 
whereas the blue cells have been updated with unlimited ADER-DG-$\mathbb{P}_9$ on the main 
grid. It is important
to note that the DG polynomials are almost continuous across all cell boundaries, except 
for the shock wave, where a real discontinuity is preserved.
In other words the unavoidable limiting has not smeared the shock wave on two cells
of the main grid, but it has been able to keep the subcell description of this discontinuity, 
thus leading to a shock wave spread over only one or two subcells. 
Note that the leftmost waves are actually contained within one single 
cell. These waves are steepening and will later become genuine shock waves\footnote{This is 
the same effect for Burgers equation when a smooth initial profile evolves 
into a shock wave.}.
As expected the \aposteriori detection procedure has flagged 
these waves for subsequent subcell correction, hence their red color.

\subsection{Double Mach reflection problem}  \label{ssec:DM}
\begin{figure}
  \hspace{-2cm}
  \begin{center}
    \begin{tabular}{c}
      \includegraphics[width=1.0\textwidth]{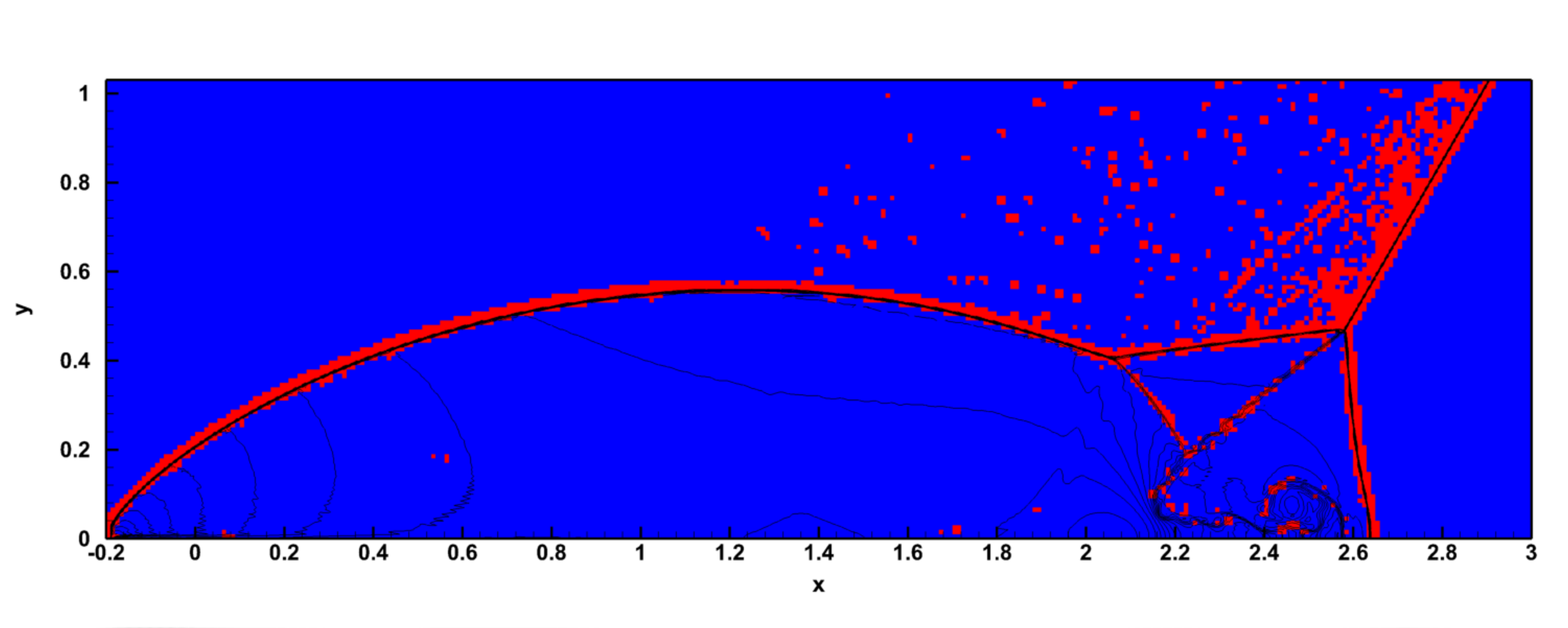}    \\
      \includegraphics[width=1.0\textwidth]{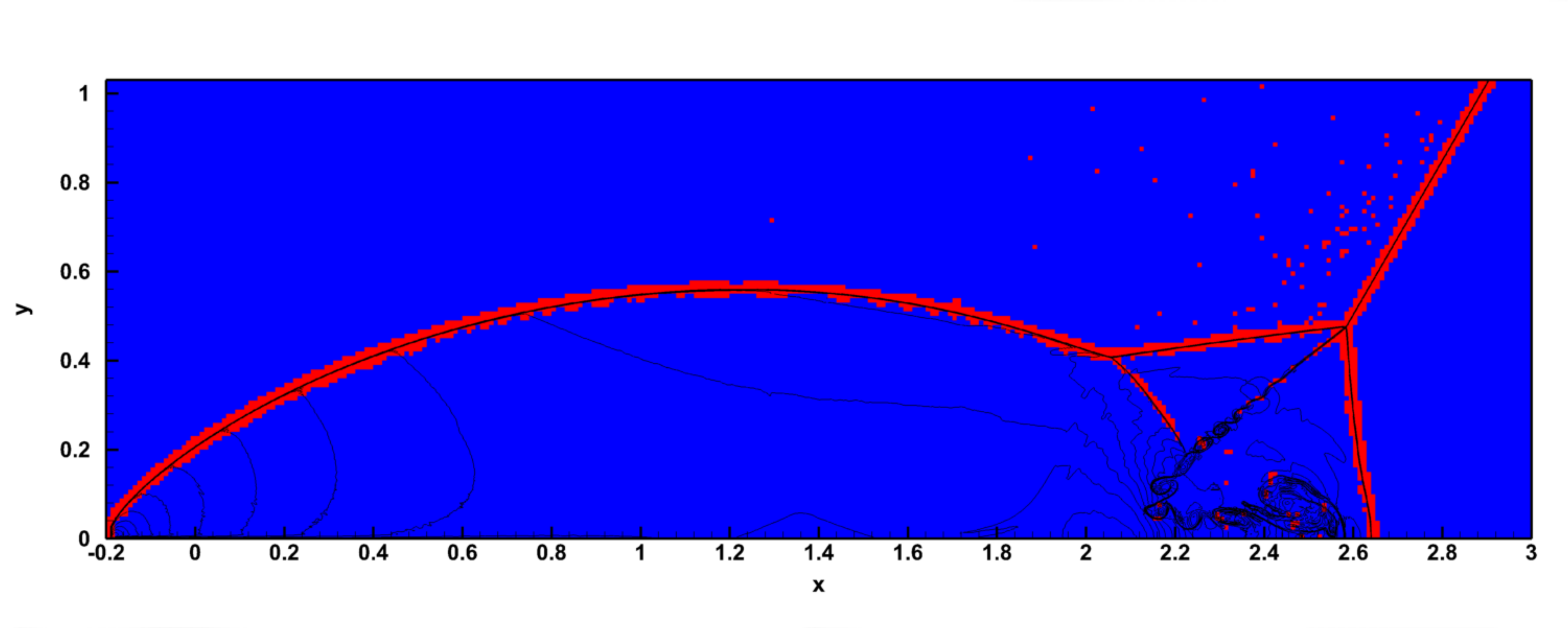}    \\
      \includegraphics[width=1.0\textwidth]{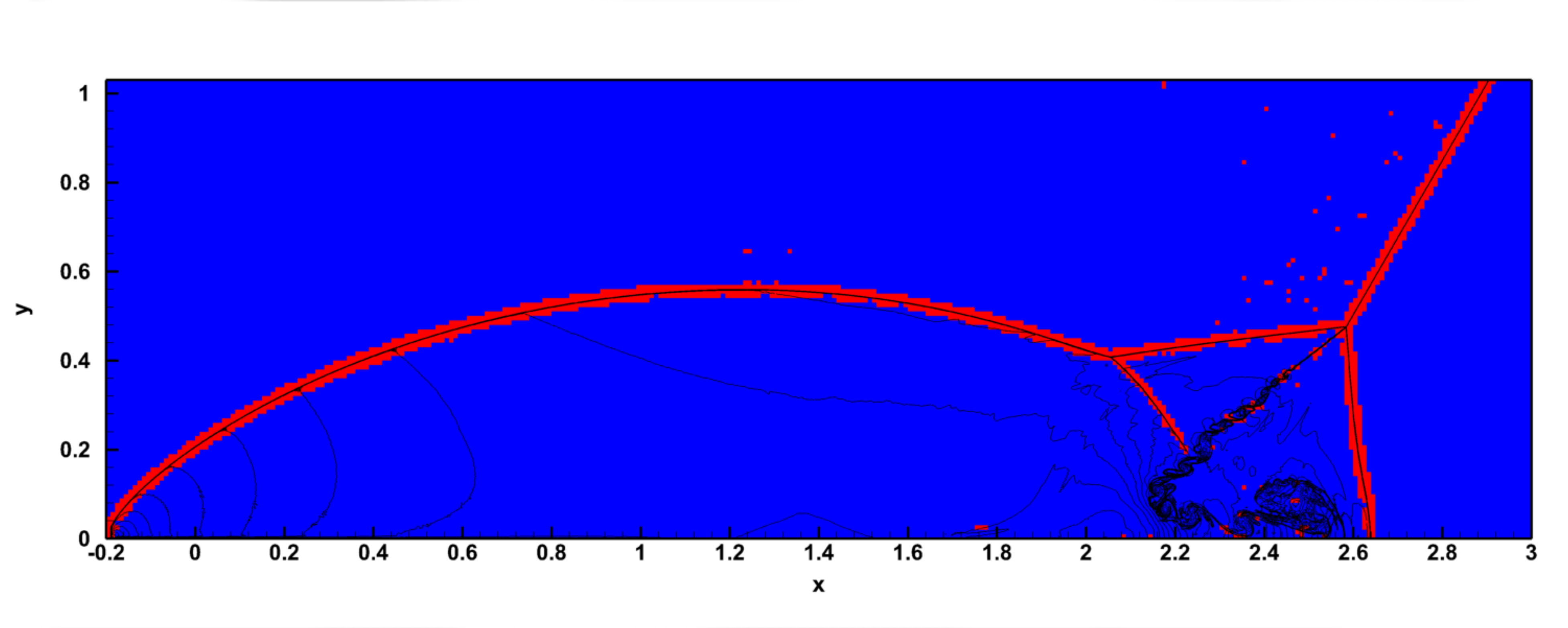}      
    \end{tabular}
    \caption{ \label{fig:DoubleMachLimiter}
      Double Mach reflection problem at $t_{\text{final}}=0.2$ with ADER-DG-$\mathbb{P}_2$, ADER-DG-$\mathbb{P}_5$ 
      and ADER-DG-$\mathbb{P}_9$ (from top to bottom) supplemented 
      with \aposteriori ADER-WENO3 subcell limiter ---
      $33$ isolines ranging from $1.5$ to $17.5$ are displayed for the density variable 
      (black lines) along with the troubled cells (red) and the unlimited cells (blue) 
      at this time level. The characteristic mesh spacing is only $h=1/100$ in all cases. 
    }
  \end{center}
\end{figure}

Next we have run the 2D double Mach reflection problem of a strong shock that was proposed in 
\cite{woodwardcol84}. This test problem involves a Mach $10$ shock in a perfect gas with 
$\gamma=1.4$ which hits a $30^\circ$ ramp with the $x$-axis. 
Using Rankine-Hugoniot conditions we can deduce the initial conditions in front of and after 
the shock wave  
\begin{eqnarray}
(\rho, u, v, p)( \x,t=0) =
\left\{
\begin{array}{cll}
  \frac{1}{\gamma}(8.0, 8.25, 0.0, 116.5), \quad & \text{ if } & \quad x'<0.1, \\
  (1.0, 0.0, 0.0, \frac{1}{\gamma}),       \quad & \text{ if } & \quad x'\geq 0.1, 
\end{array}
\right.
\end{eqnarray}
where $x'$ is the coordinate in a rotated coordinate system. Reflecting wall boundary conditions are prescribed  
on the bottom and the exact solution of an isolated moving oblique shock wave with shock Mach 
number $M_s=10$ is imposed on the upper boundary. A Rusanov (local Lax-Friedrichs) flux has been used for this test problem. 
Inflow and outflow boundary conditions are prescribed on the left side and the right side, respectively. \\
The computational domain is given by $\Omega = [0;3.5] \times [0;1]$ and the main grid is built using a characteristic length of 
$h=1/100$, leading to $350\times 100$ computational cells.  
We solve this problem with four schemes: ADER-DG-$\mathbb{P}_N$
with $N=1$, $2$, $5$ and $9$, all supplemented with the \aposteriori ADER-WENO3 subcell limiter. 
The results are depicted using $33$ isolines ranging from $1.5$ to $17.5$ 
for the density variable at $t_{\text{final}}=0.2$, see Fig~\ref{fig:DoubleMachLimiter}.
We first present the density contour lines along with 
the cells colored by the indicator function $\beta_i^{n+1}$: red cells have been recomputed 
using the subcell ADER-WENO3 scheme ($\beta_i^{n+1}=1$), whereas blue cells ($\beta_i^{n+1}=0$) have 
kept the original unlimited ADER-DG-$\mathbb{P}_N$ method. The shock waves are well detected for 
all DG schemes. Moreover, the interaction zone where a lot of small scale vortex structures develop 
is almost free of limiting for ADER-DG-$\mathbb{P}_5$ and ADER-DG-$\mathbb{P}_9$ and \textit{more}  
troubled cells are detected by our \aposteriori method for a \textit{low} order ADER-DG-$\mathbb{P}_2$ scheme, 
which is quite remarkable. One would have intuitively expected that higher order schemes tend to create 
more oscillations and thus need more limiting. It seems, however, that the higher order DG scheme has 
\textit{better subcell resolution} capabilities and therefore can treat these vortex structures as smooth 
subscale features, without the need of a limiter.  
Also more erroneously detected troubled cells are observed for the ADER-DG-$\mathbb{P}_2$ scheme in the upper part
of the figure, behind the main shock wave. 
Only minor details seem to differ in the indicator function between ADER-DG-$\mathbb{P}_5$ and ADER-DG-$\mathbb{P}_9$. 
Nonetheless, the numerical structures represented by the density isolines seem to be richer with the ADER-DG-$\mathbb{P}_9$ 
scheme.
 
\begin{figure}
  \begin{center}
    \begin{tabular}{cc}
      \includegraphics[width=0.52\textwidth]{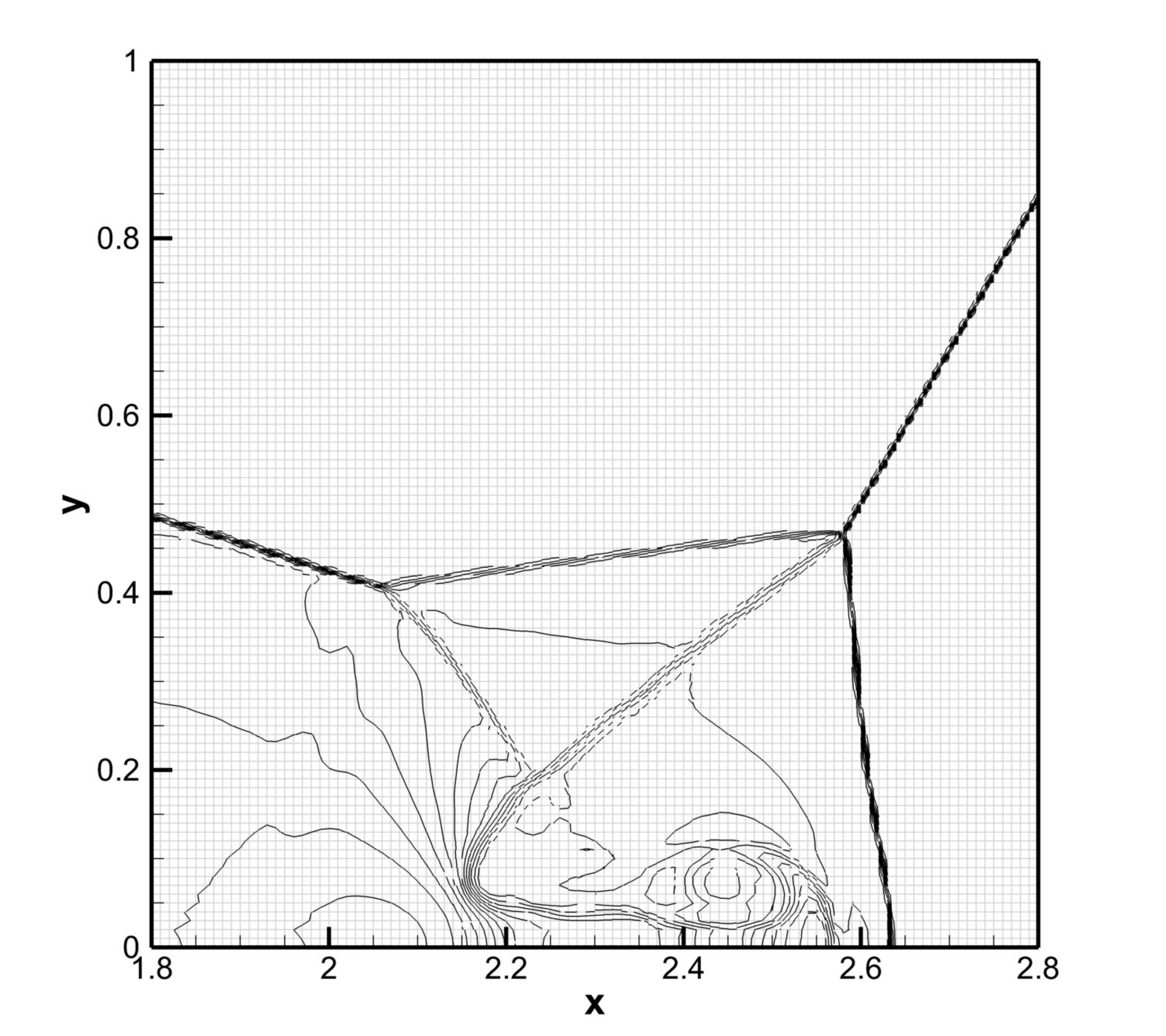}  &
      \hspace{-0.9cm}  
      \includegraphics[width=0.52\textwidth]{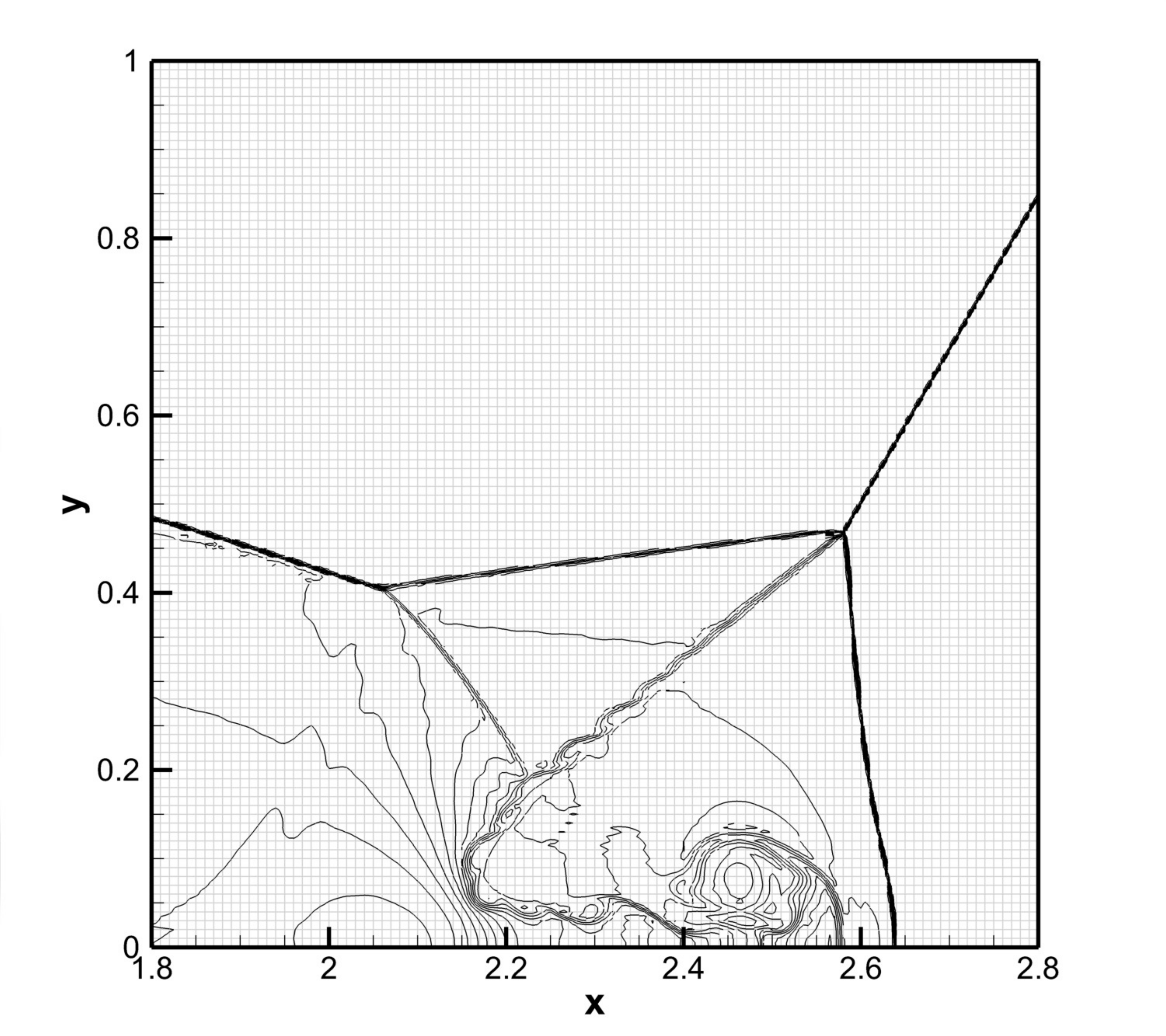}  \\     
      \includegraphics[width=0.52\textwidth]{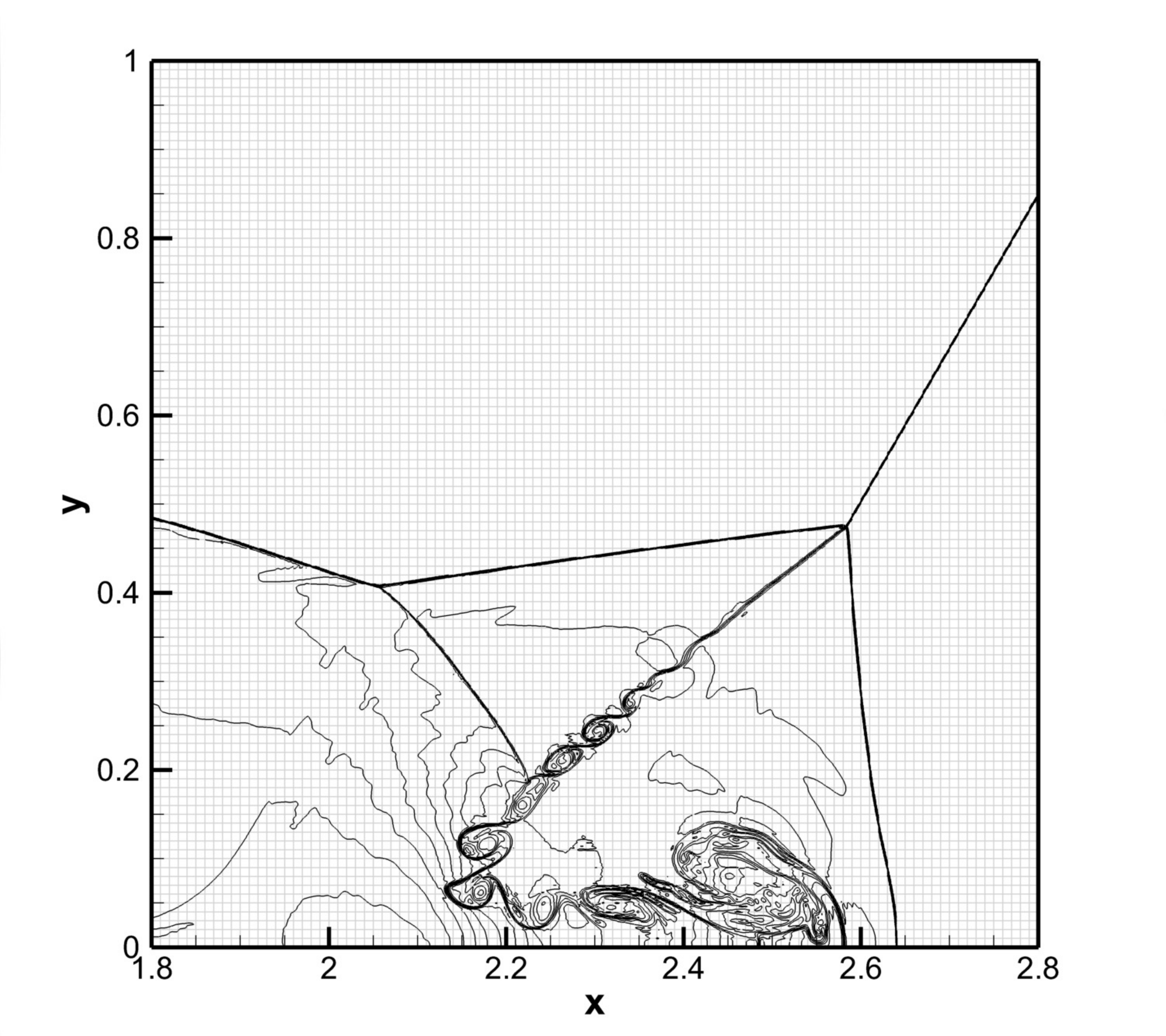}  &
      \hspace{-0.9cm}  
      \includegraphics[width=0.52\textwidth]{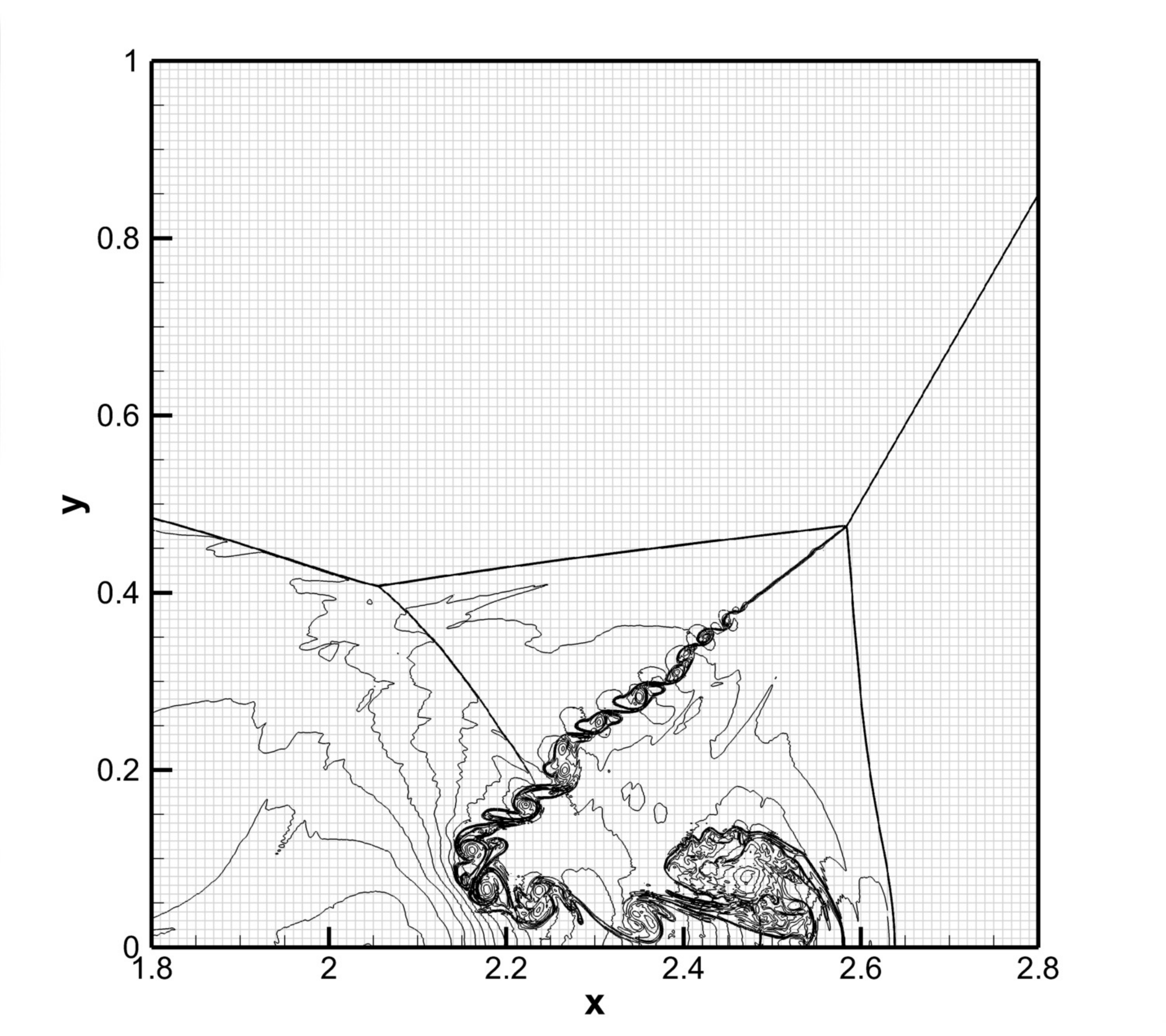}      
    \end{tabular}    
    \caption{ \label{fig:DoubleMachzoom}
      Double Mach reflection problem at $t_{\text{final}}=0.2$ with ADER-DG-$\mathbb{P}_1$, 
      $\mathbb{P}_2$, $\mathbb{P}_5$ and $\mathbb{P}_9$ (from top left to bottom right) 
      supplemented with \aposteriori ADER-WENO3 subcell limiter --- 
      Zooms on the interaction zone.
      $33$ isolines ranging from $1.5$ to $17.5$ are displayed for the density variable 
      (black lines). The mesh is underlined in gray color. 
    }
  \end{center}
\end{figure}

Overall, also for this very difficult problem which exhibits at the same time strong shock waves and smooth flow features, 
our detection procedure seems to identify discontinuities associated with shock waves properly and it ignores complex but 
smooth vortex-type flow structures, as expected. 
In order to confirm this observation we propose in Fig~\ref{fig:DoubleMachzoom}
two zooms on the interaction zone for ADER-DG-$\mathbb{P}_N$ schemes with $N=1$, $2$, $5$, and $9$.
The main grid is also plotted in light gray color to compare the numerical thickness of a 
shocks wave and the size of the vortex-type structures with respect to the main grid size. 
From this figure we confirm that the DG scheme behaves better with increasing $N$, 
the degree of the DG polynomial. ADER-DG-$\mathbb{P}_9$ dissipates less and, as such, is capturing 
more structures than the other schemes. 
Also the smearing of the shocks is less pronounced with ADER-DG-$\mathbb{P}_9$ than with  
ADER-DG-$\mathbb{P}_5$ on the same grid.  
All schemes capture the shock waves properly without spurious oscillations, but  
the numerical dissipation of the lower order schemes destroys the small scale vortex 
structures that are obtained on this rather coarse grid only by the highest order 
DG schemes. The limiter does not seem to destroy this ability. In other words, the 
subcell resolution property of the DG scheme seems to be very well preserved. \\ 
From the results of this test case we can conclude that the \aposteriori subcell 
limiter provides a valid detection strategy of troubled cells and that the 
ADER-WENO3 scheme used on the subgrid is able to preserve the overall accuracy 
of the high order DG scheme on the main grid.

\subsection{Forward facing step}  \label{ssec:FF}
\begin{figure}
  \begin{center}
    \vspace{-2.0cm}
    \begin{tabular}{c}
    \vspace{-0.65cm}
      \includegraphics[width=1.0\textwidth]{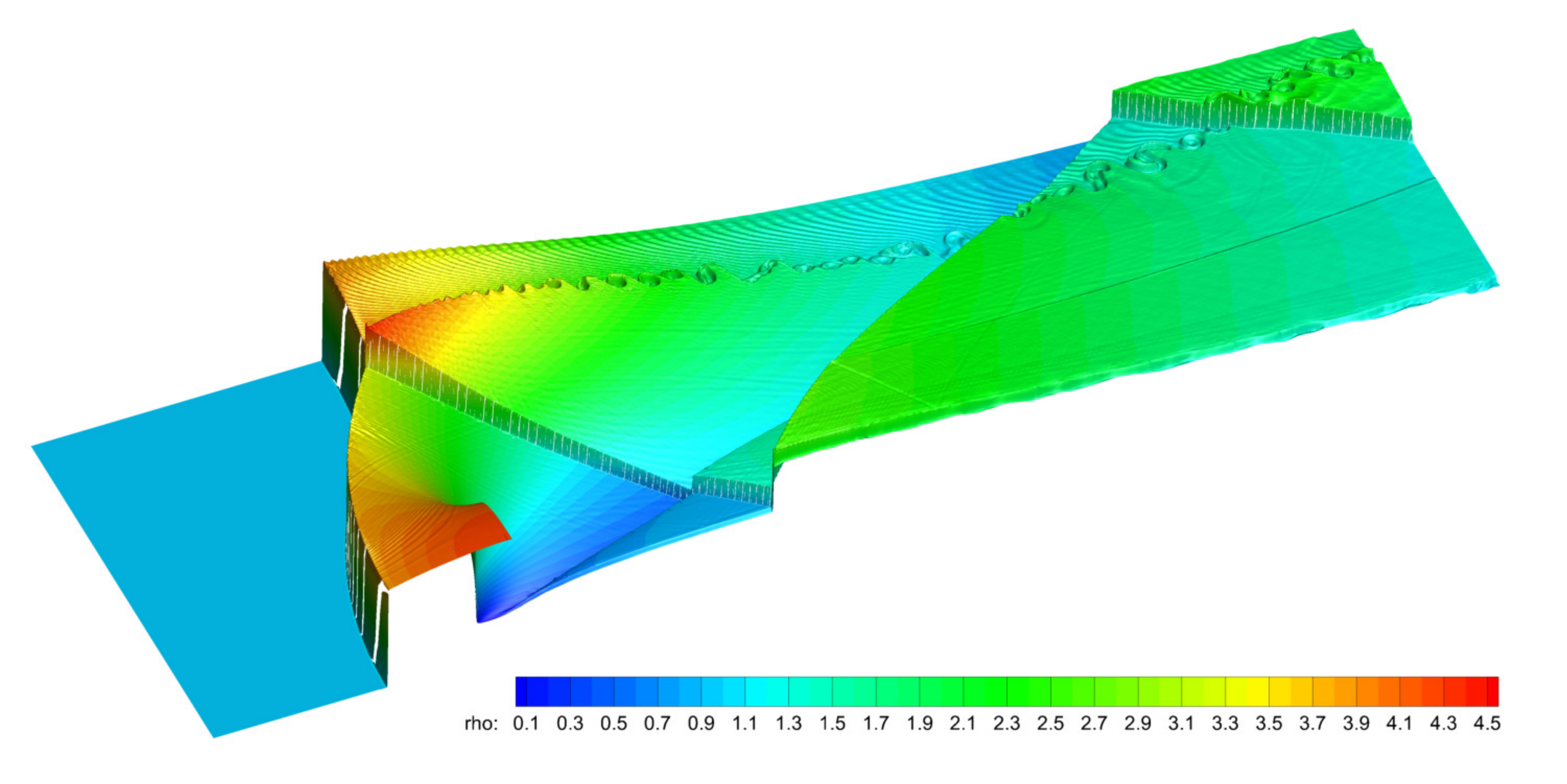}  \\
      \vspace{-1.1cm}
      \includegraphics[width=1.0\textwidth]{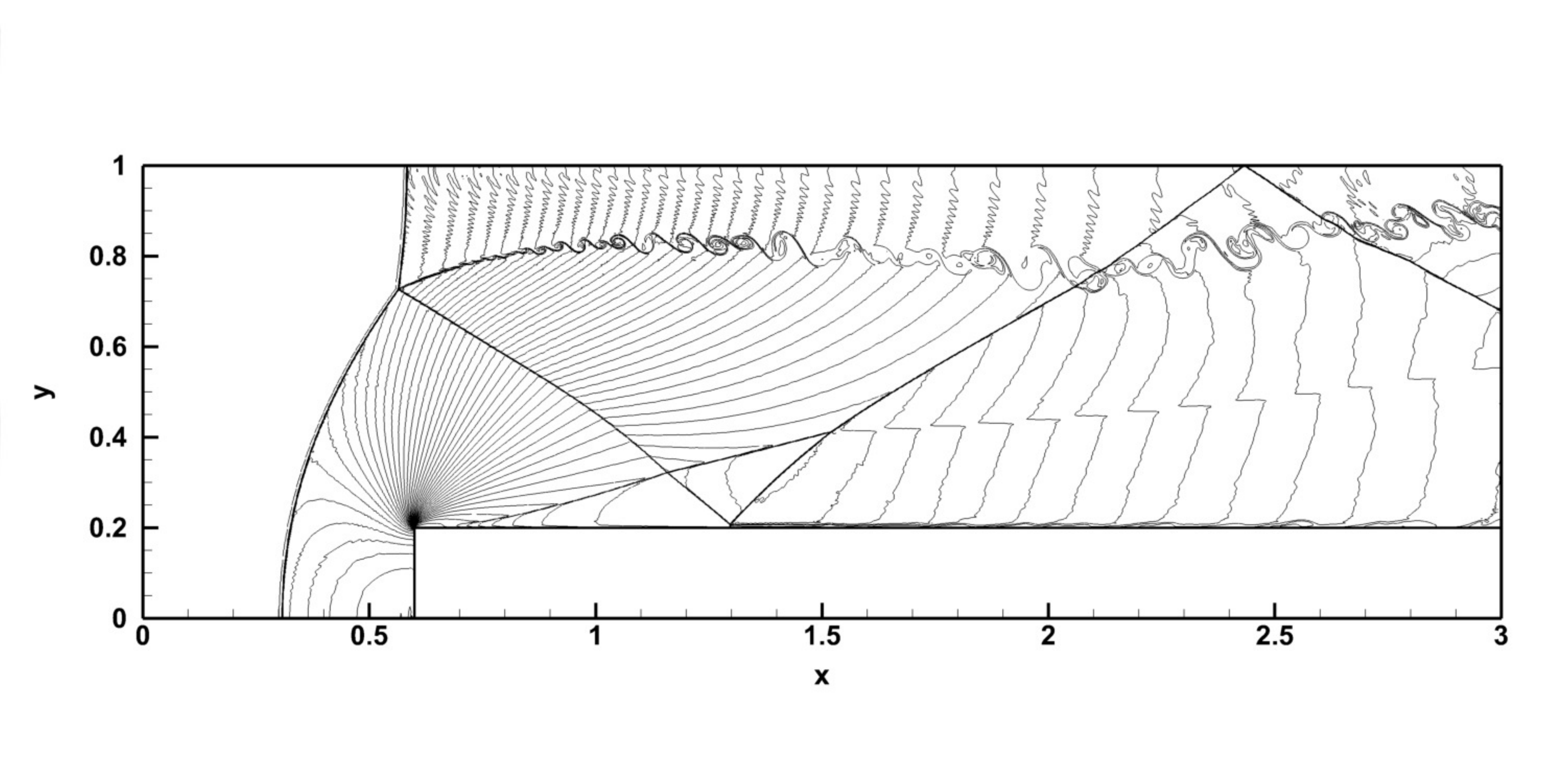}  \\
      \vspace{-1.0cm}
      \includegraphics[width=1.0\textwidth]{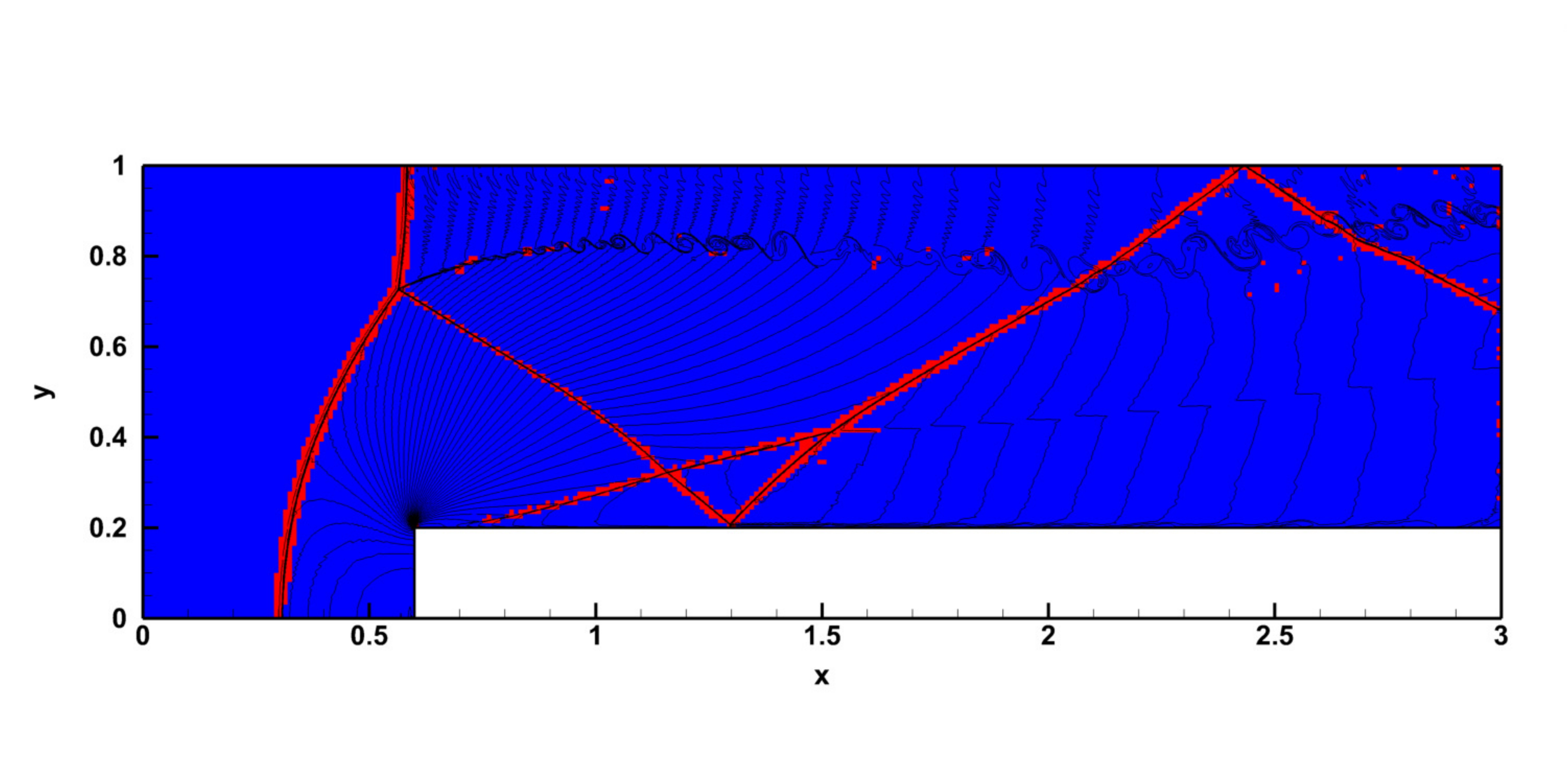}      
    \end{tabular}    
    \caption{ \label{fig:FFS}
      Forward Facing Step problem at $t_{\text{final}}=4$ with ADER-DG-$\mathbb{P}_5$ supplemented
      with \aposteriori ADER-WENO3 subcell limiter ---      
      Top-panel: extruded 3D density (color and azimuth) ---
      Middle-panel: density represented as isolines ---
      Bottom-panel: troubled cells (red) updated using the subcell ADER-WENO3
      method and unlimited ADER-DG-$\mathbb{P}_5$ cells (blue).
    }
  \end{center}
\end{figure}

In this section we simulate the so called forward facing step (FFS) problem, 
also proposed by Woodward and Colella in \cite{woodwardcol84}.
It consists of a Mach 3 wind tunnel with a step. 
The initial condition consists in a uniform gas with density $\rho=\gamma$, 
pressure $p=1$, velocity components $u=3, v=0$ and $\gamma=1.4$. The computational 
domain is given by $\Omega = [0;3] \times [0;1] \backslash [0.6;3] \times [0;0.2]$ 
and reflective wall boundary conditions are applied on the upper and lower boundary 
of the domain, whereas inflow and outflow boundary conditions are applied at the 
entrance and the exit. The solution of this problem 
involves shock waves interacting with the boundaries and it is run up to a final 
time of $t_{\text{final}}=4$. 
The mesh is uniform, using a characteristic length of $h=1/100$, leading to about 
$N=300 \times 100$ cells on the main grid.  
For this test we run the ADER-DG-$\mathbb{P}_5$ scheme with a Rusanov flux, supplemented 
with the \aposteriori ADER-WENO3 subcell limiter (SCL).  
In Fig.~\ref{fig:FFS} we display the density variable at the final time.  
The top panel represents an extruded density (the azimuth and the color
represent the numerical density). From this view we can clearly see 
several smooth zones separated by the shock waves. Furthermore, we can 
distinguish the Kelvin-Helmholtz instabilities developing along the top 
shear wave. These instabilities are maintained when they move across 
the shock waves. Moreover, the unsteady vortices create small amplitude 
acoustic waves, which are expected from such unsteady flows. Yet, they 
are often smeared by more dissipative numerical schemes. 
On the same figure we propose in the middle panel the classical  
black and white isolines in contrast with the previous figure.
It is clear that the perturbed isolines emanating from the vortices are due to 
the acoustic waves, previously seen on the extruded 3D view. 
Finally the bottom panel shows the limited cells updated with the subcell WENO3  
scheme (in red) and the unlimited ADER-DG-$\mathbb{P}_5$ cells (in blue). As 
expected, the great majority of the computational domain is computed with the 
unlimited ADER-DG-$\mathbb{P}_5$ scheme, apart from the shock waves, where the 
subcell limiter mechanism seems to act properly. 

\subsection{2D Riemann problem}  \label{ssec:2DRP}

In this section we consider a set of two-dimensional Riemann problems which has been proposed and extensively studied in 
\cite{schulzrinne,kurganovtadmor} and which has also been recently used to construct genuinely multidimensional Riemann 
solvers of the HLL type, see \cite{balsarahlle2d,balsarahllc2d,BalsaraMultiDRS}. The computational domain is 
$\Omega = [-0.5;0.5] \times [-0.5;0.5]$ and the initial conditions are given by 
\begin{equation}
 \mathbf{u}(x,y,t=0) = \left\{ \begin{array}{ccc} 
 \mathbf{u}_1 & \textnormal{ if } & x > 0 \wedge y > 0,    \\ 
 \mathbf{u}_2 & \textnormal{ if } & x \leq 0 \wedge y > 0, \\ 
 \mathbf{u}_3 & \textnormal{ if } & x \leq 0 \wedge y \leq 0, \\ 
 \mathbf{u}_4 & \textnormal{ if } & x > 0 \wedge y \leq 0.    
  \end{array} \right. 
\end{equation} 
The initial conditions and the final simulation time, $t_{\text{final}}$, for the five configurations presented in this article 
are listed in Table~\ref{tab.rp2d.ic}. For more information about the other configurations the reader is referred to 
\cite{schulzrinne,kurganovtadmor}. For the simulation we have used a main grid of $100 \times 100$ elements with a 
DG-$\mathbb{P}_5$ scheme supplemented with our \aposteriori WENO3 subcell limiter. 
\begin{table}[!b] 
  \numerikNine  
\begin{center} 
\begin{tabular}{|c|c||cccc|cccc|c|} 
\hline
\textbf{\#}    && $\rho$ & $u$ & $v$  & $p$ & $\rho$ & $u$ & $v$  & $p$ & \multirow{2}{*}{$t_{\text{final}}$} \\ 
\cline{3-10}
    && \multicolumn{4}{c|}{$x \leq 0$} & \multicolumn{4}{c|}{$x>0$} & \\
\hline
\hline
\multirow{2}{*}{\rotatebox{90}{\textbf{RP1}}} 
&$y > 0$    & 0.5323 & 1.206   & 0.0     & 0.3   & 1.5    & 0.0 &  0.0    & 1.5 &  \multirow{2}{*}{0.25} \\ 
&$y \leq 0$ & 0.138  & 1.206   & 1.206   & 0.029 & 0.5323 & 0.0 &  1.206  & 0.3 & \\ 
\hline
\hline
\multirow{2}{*}{\rotatebox{90}{\textbf{RP2}}}
&$y > 0$    & 0.5065 &  0.8939 & 0.0     & 0.35 & 1.1    & 0.0 &  0.0    & 1.1  &\multirow{2}{*}{0.25}\\ 
&$y \leq 0$ & 1.1    &  0.8939 & 0.8939  & 1.1  & 0.5065 & 0.0 &  0.8939 & 0.35 &\\ 
\hline
\hline
\multirow{2}{*}{\rotatebox{90}{\textbf{RP3}}}
&$y > 0$    & 2.0    &  0.75  & 0.5   & 1.0  & 1.0  &  0.75 &  -0.5  & 1.0  &\multirow{2}{*}{0.30}\\ 
&$y \leq 0$ & 1.0    & -0.75  & 0.5   & 1.0  & 3.0  & -0.75 &  -0.5  & 1.0  &\\ 
\hline
\hline
\multirow{2}{*}{\rotatebox{90}{\textbf{RP4}}}
&$y > 0$    & 1.0   & -0.6259& 0.1   & 1.0  & 0.5197  &  0.1 & 0.1  & 0.4  &\multirow{2}{*}{0.25}\\ 
&$y \leq 0$ & 0.8   & 0.1    & 0.1   & 1.0  & 1.0  & 0.1 & -0.6259  & 1.0  &\\ 
\hline
\hline
\multirow{2}{*}{\rotatebox{90}{\textbf{RP5}}}
&$y > 0$    & 1.0    &  0.7276 & 0.0     & 1.0  & 0.5313 & 0.0 &  0.0    & 0.4  &\multirow{2}{*}{0.25}\\ 
&$y \leq 0$ & 0.8    &  0.0    & 0.0     & 1.0  & 1.0    & 0.0 &  0.7276 & 1.0  &\\ 
\hline 
\end{tabular} 
\caption{Initial conditions for the 2D Riemann problems numbered from $1$ to $5$. 
  These further correspond to Configurations 3, 4, 6, 8 
  and 12 in \cite{kurganovtadmor}} 
\end{center}
\label{tab.rp2d.ic}
\end{table} 

\begin{figure}
  \vspace{-1cm}
  \begin{center}
    \begin{tabular}{cc}
      \vspace{-0.2cm}
      \includegraphics[width=0.47\textwidth]{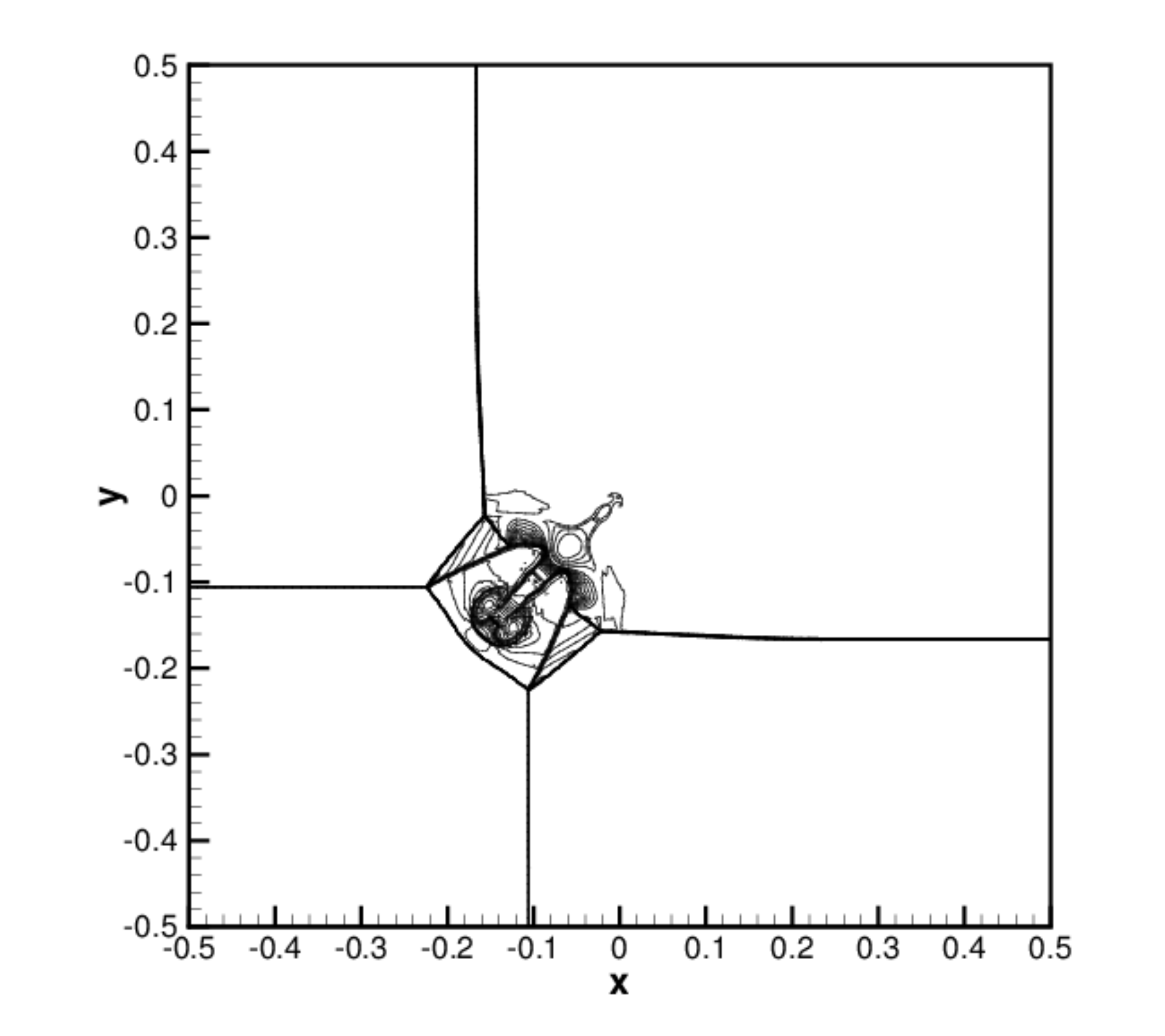}  & 
      \vspace{-0.2cm}
      \includegraphics[width=0.47\textwidth]{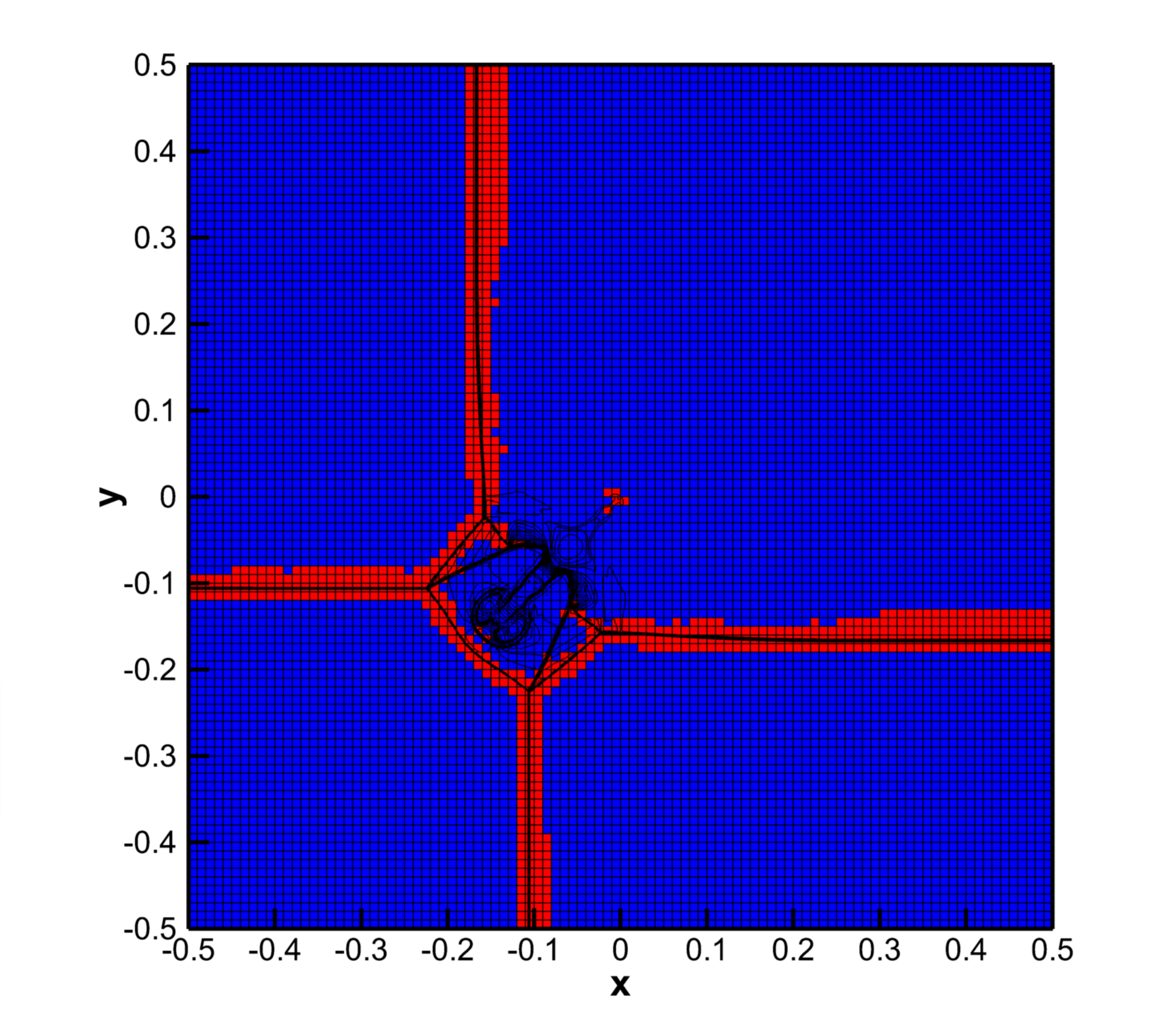}  \\ 
      \vspace{-0.2cm}
      \includegraphics[width=0.47\textwidth]{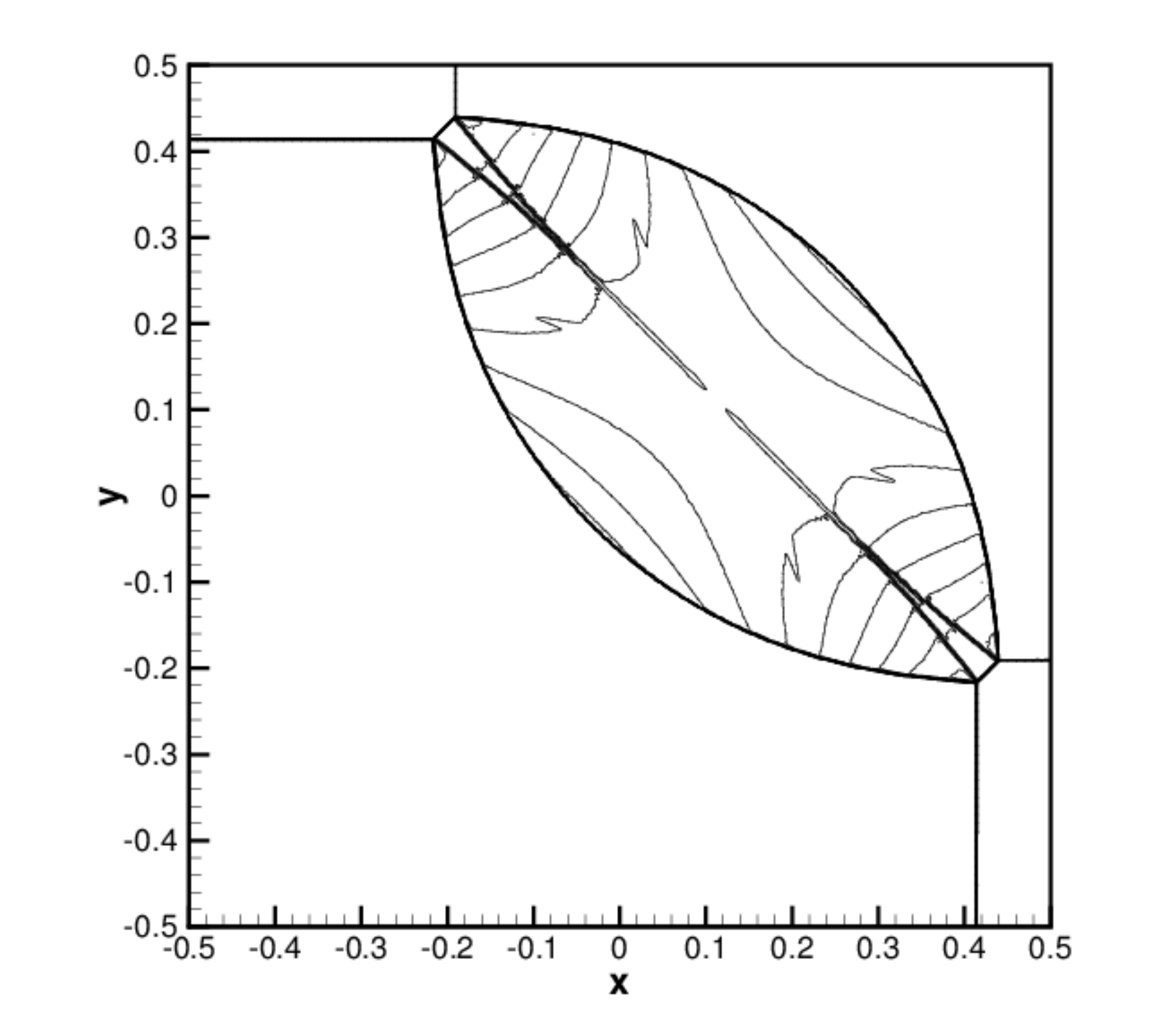}  & 
      \vspace{-0.2cm}
      \includegraphics[width=0.47\textwidth]{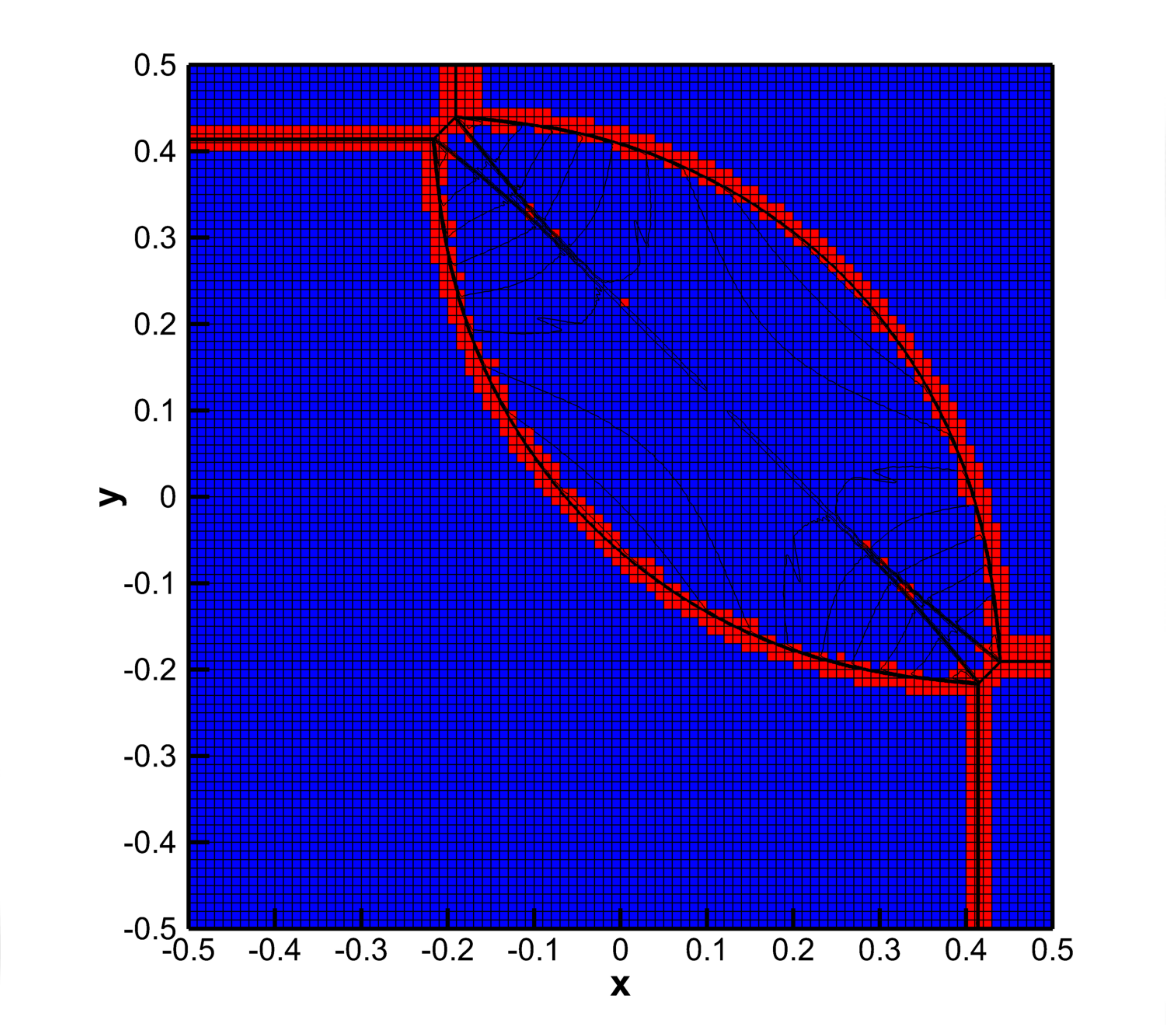}  \\ 
      \vspace{-0.2cm}
      \includegraphics[width=0.47\textwidth]{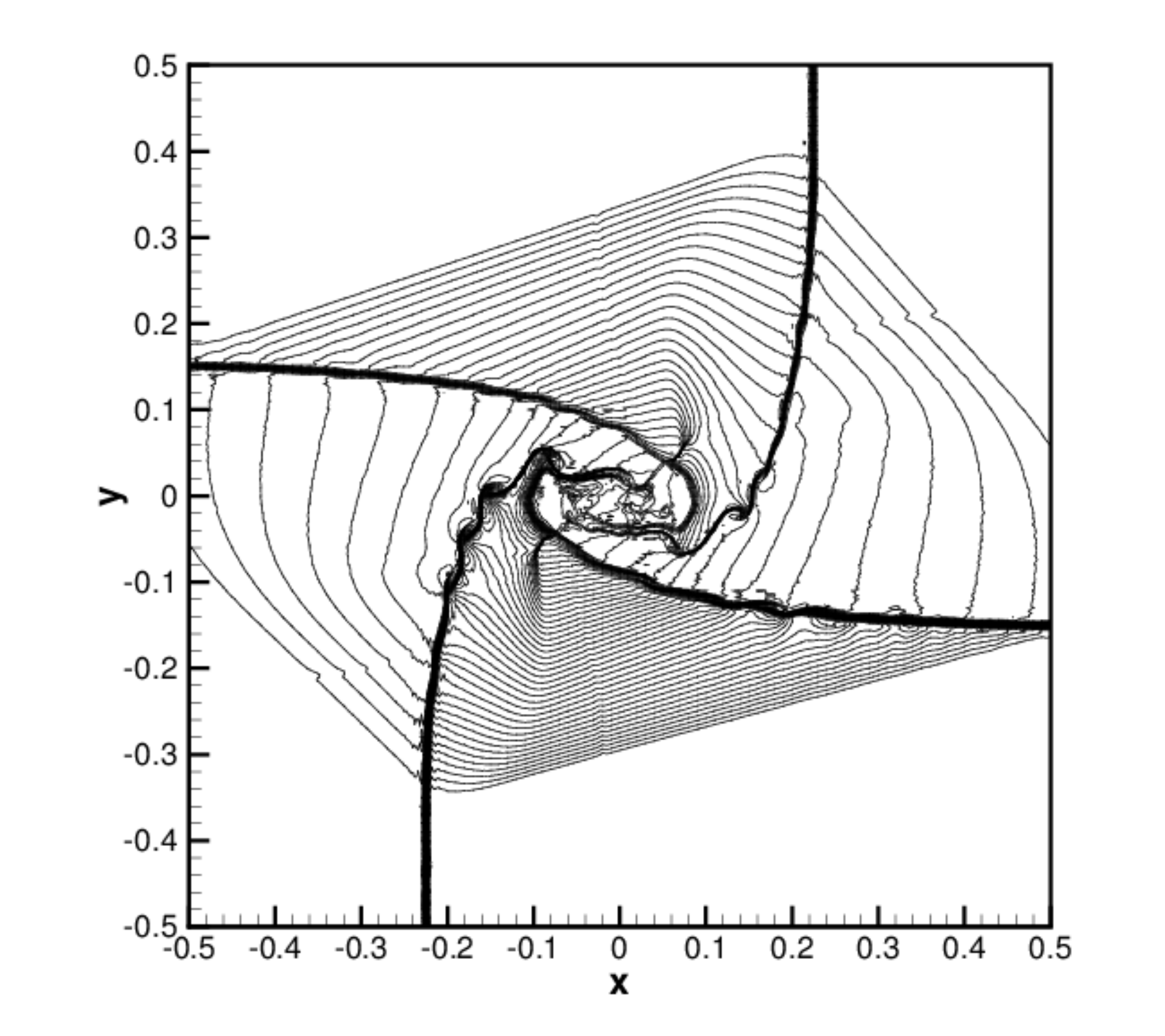}  & 
      \vspace{-0.2cm}
      \includegraphics[width=0.47\textwidth]{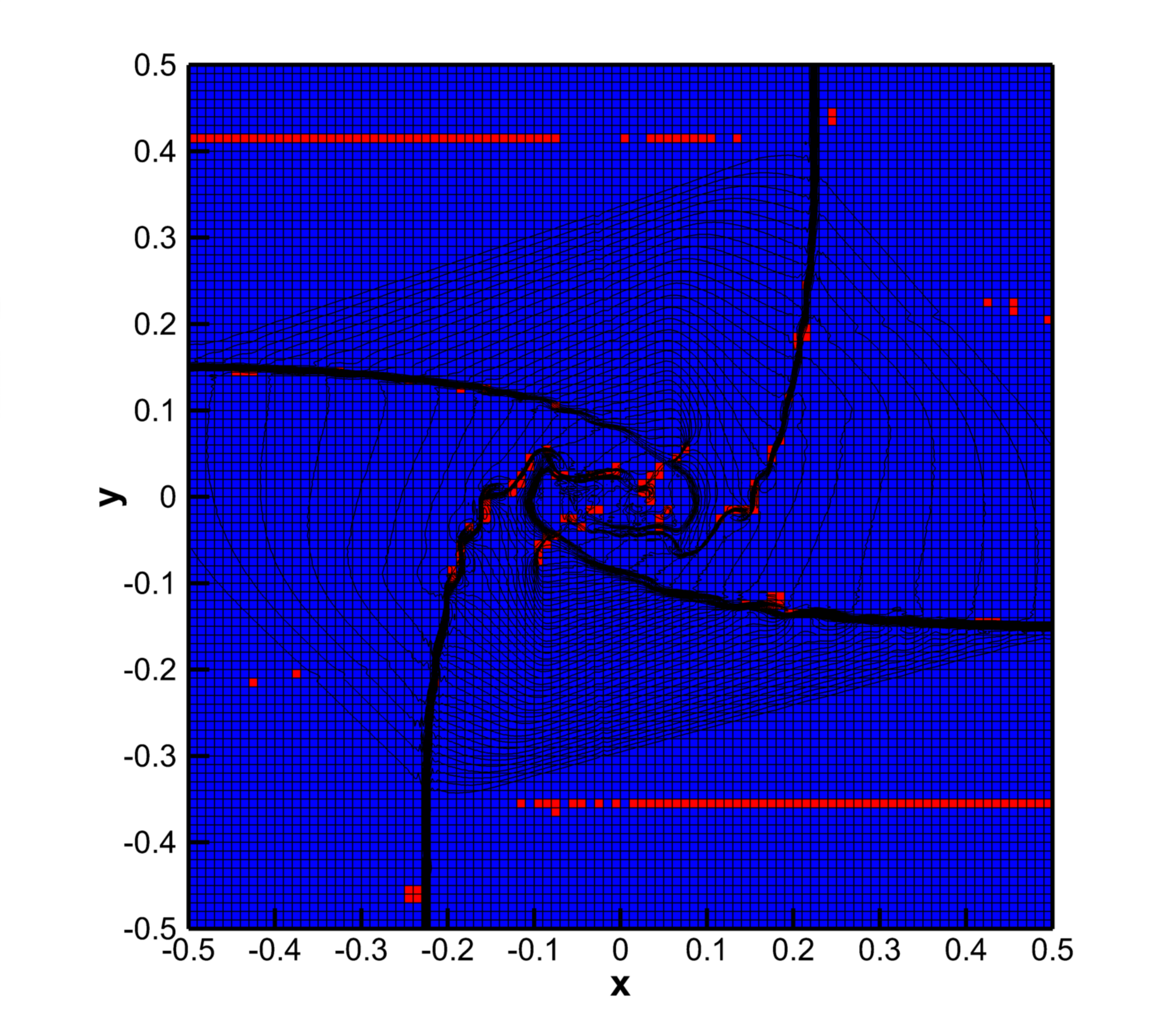} 
    \end{tabular}    
    \caption{ \label{fig:RP2D}
      2D Riemann problems simulated with an ADER-DG-$\mathbb{P}_5$ method supplemented
      with \aposteriori ADER-WENO3 subcell limiter ---  
      Left panels: isolines of the density variable ---
      Right panels: limited cells (red) and unlimited cells (blue). 
      From top to bottom \textbf{RP1}, \textbf{RP2} and \textbf{RP3}.
    }
  \end{center}
\end{figure}

\begin{figure}
  \begin{center}
    \begin{tabular}{cc}
      \vspace{-0.2cm}
      \includegraphics[width=0.47\textwidth]{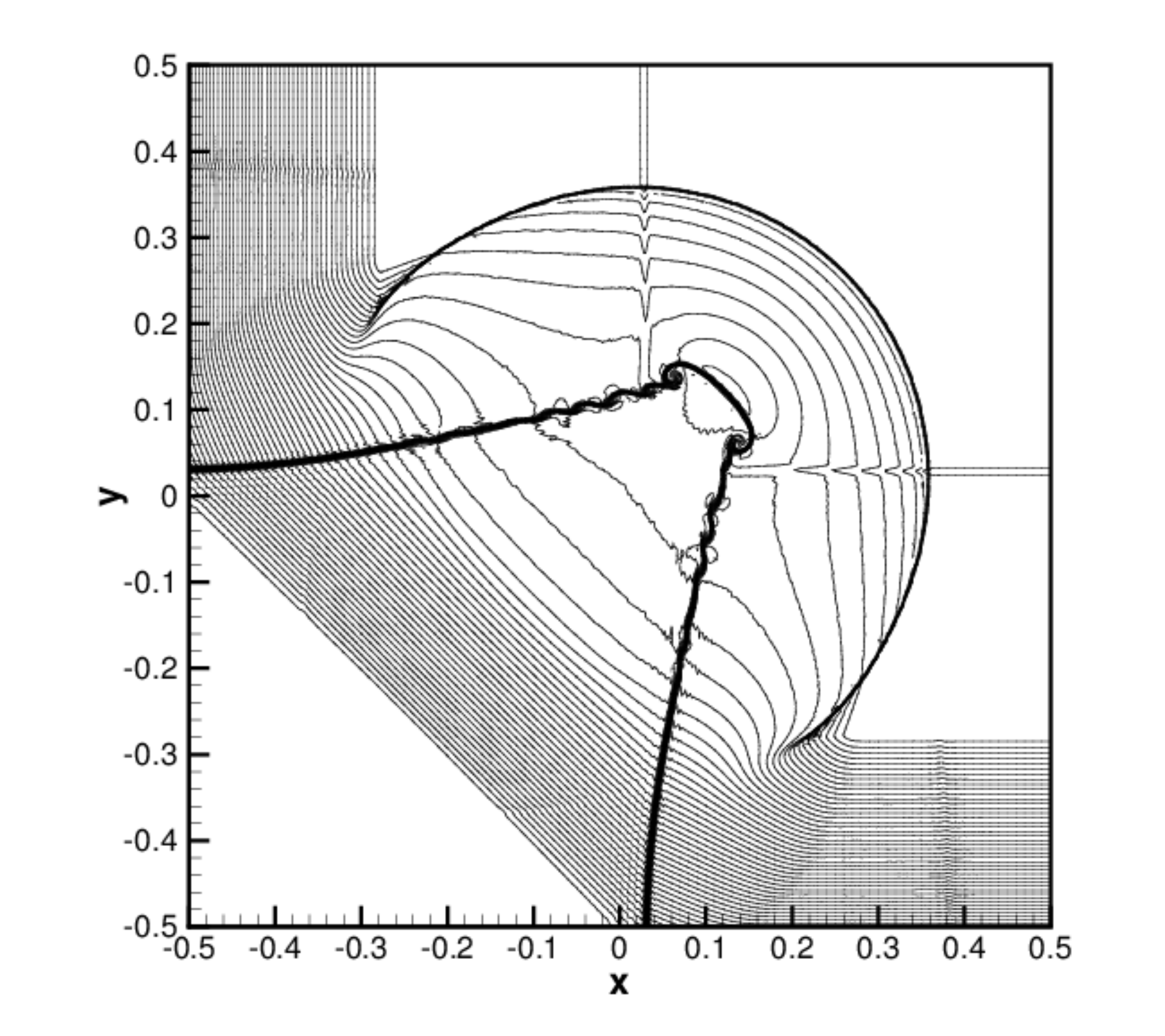}  & 
      \vspace{-0.2cm}
      \includegraphics[width=0.47\textwidth]{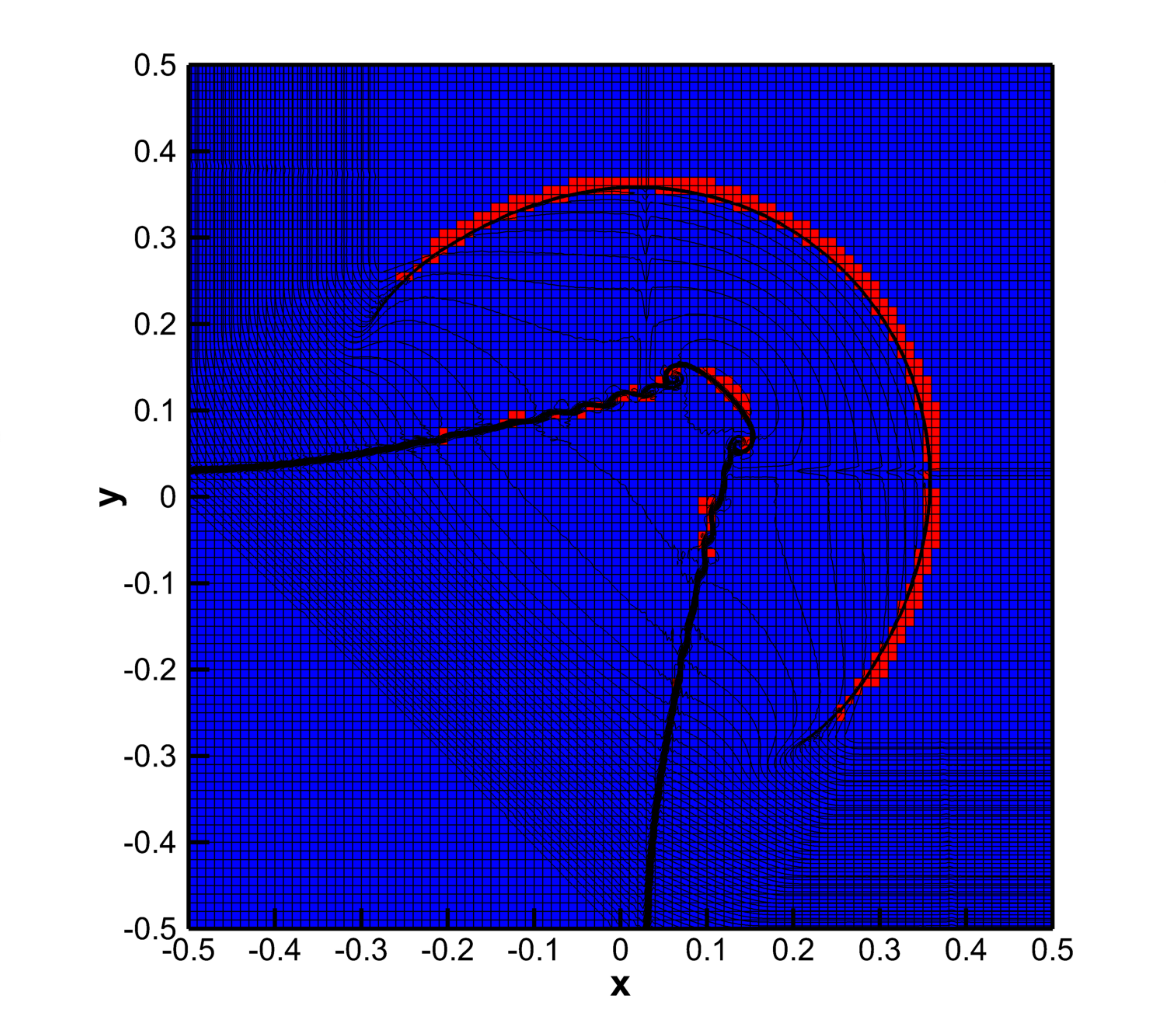}  \\ 
      \vspace{-0.2cm}
      \includegraphics[width=0.47\textwidth]{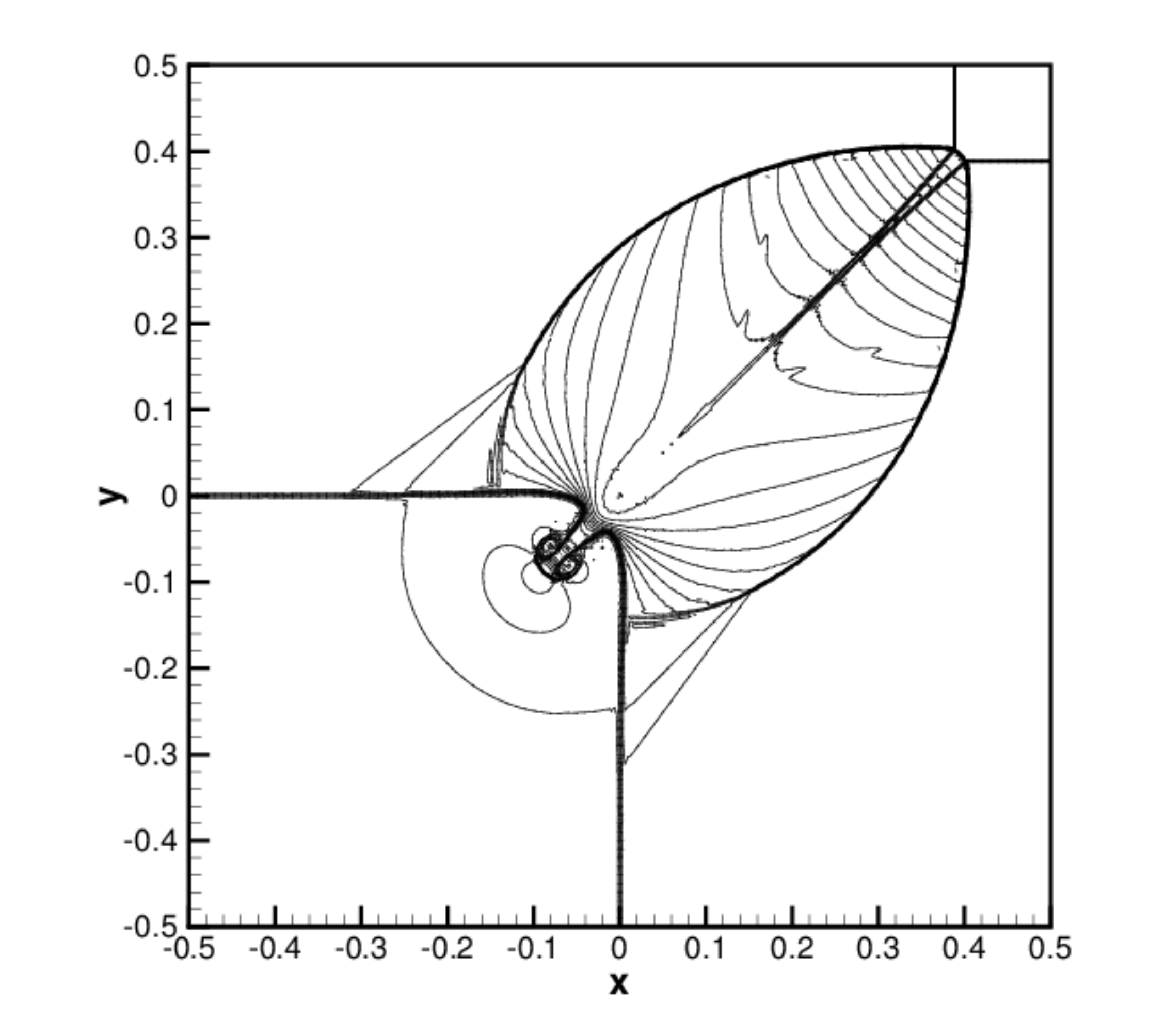}  & 
      \vspace{-0.2cm}
      \includegraphics[width=0.47\textwidth]{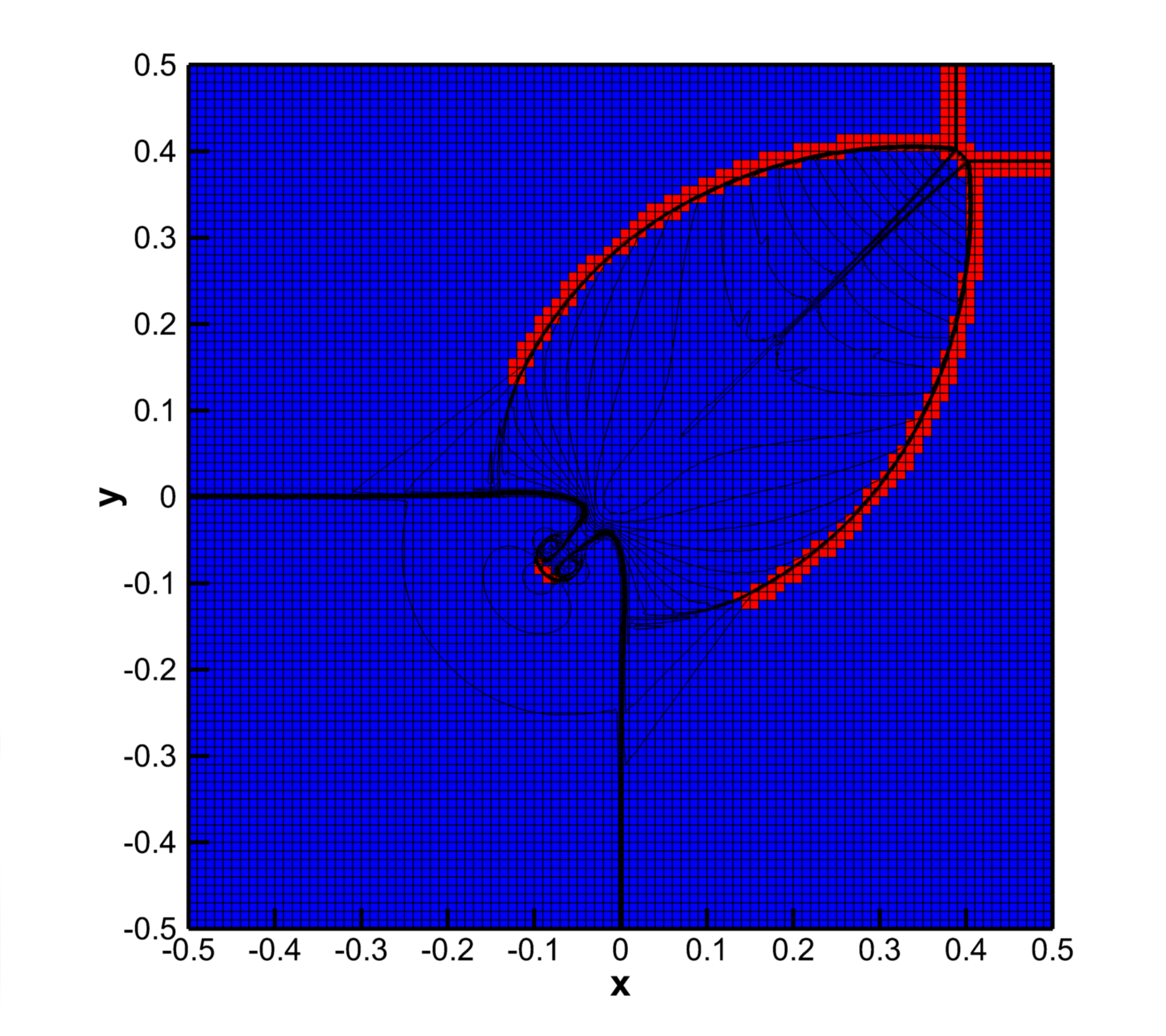}  \\ 
    \end{tabular}    
    \caption{ \label{fig:RP2Dbis}
      Same caption as Fig.\ref{fig:RP2D}, but for
       \textbf{RP4} and \textbf{RP5} (from top to bottom).
    }
  \end{center}
\end{figure}

In Figs.~\ref{fig:RP2D}-\ref{fig:RP2Dbis} the numerical results are presented for the first three and for the last 
two Riemann problems, respectively. In the left panels we show the distribution of the density at the final time, 
with equidistant isolines between the minimum and the maximum value. The number of isolines is $30$ for \textbf{RP1} 
and \textbf{RP2}, $50$ for \textbf{RP3} and \textbf{RP4}, and $40$ for \textbf{RP5}. 
In the right panels we show, as usual, the corresponding mesh and the cells, colored in red, which have been updated with 
the subcell WENO3 scheme, while the unlimited cells are marked in blue. It seems that the limiter is active only along strong  
discontinuities. Apart from these waves (and some parasitical cells along with the start-up error for \textbf{RP3}) the 
limiter is inactive in smooth regions, leading to an optimal precision given by the unlimited ADER-DG-$\mathbb{P}_5$ scheme. \\ 
The computational results, and in particular the generation of the main structures in all these two dimensional Riemann problems, 
are in good agreement with the literature \cite{kurganovtadmor}.  
Nonetheless, in these simulations the numerical dissipation is drastically reduced and, as such, the numerical solution shows much 
more small-scale features than usually reported, which are for example due to the Kelvin-Helmholtz instability along shear waves 
(see \textbf{RP3} and \textbf{RP4} in Figs.~\ref{fig:RP2D}-\ref{fig:RP2Dbis} respectively). 

\begin{figure}
  \begin{center}
    \begin{tabular}{cc}
      \vspace{-0.25cm}
      \includegraphics[width=0.45\textwidth]{./RP2DC6}  & 
      \vspace{-0.25cm}
      \includegraphics[width=0.45\textwidth]{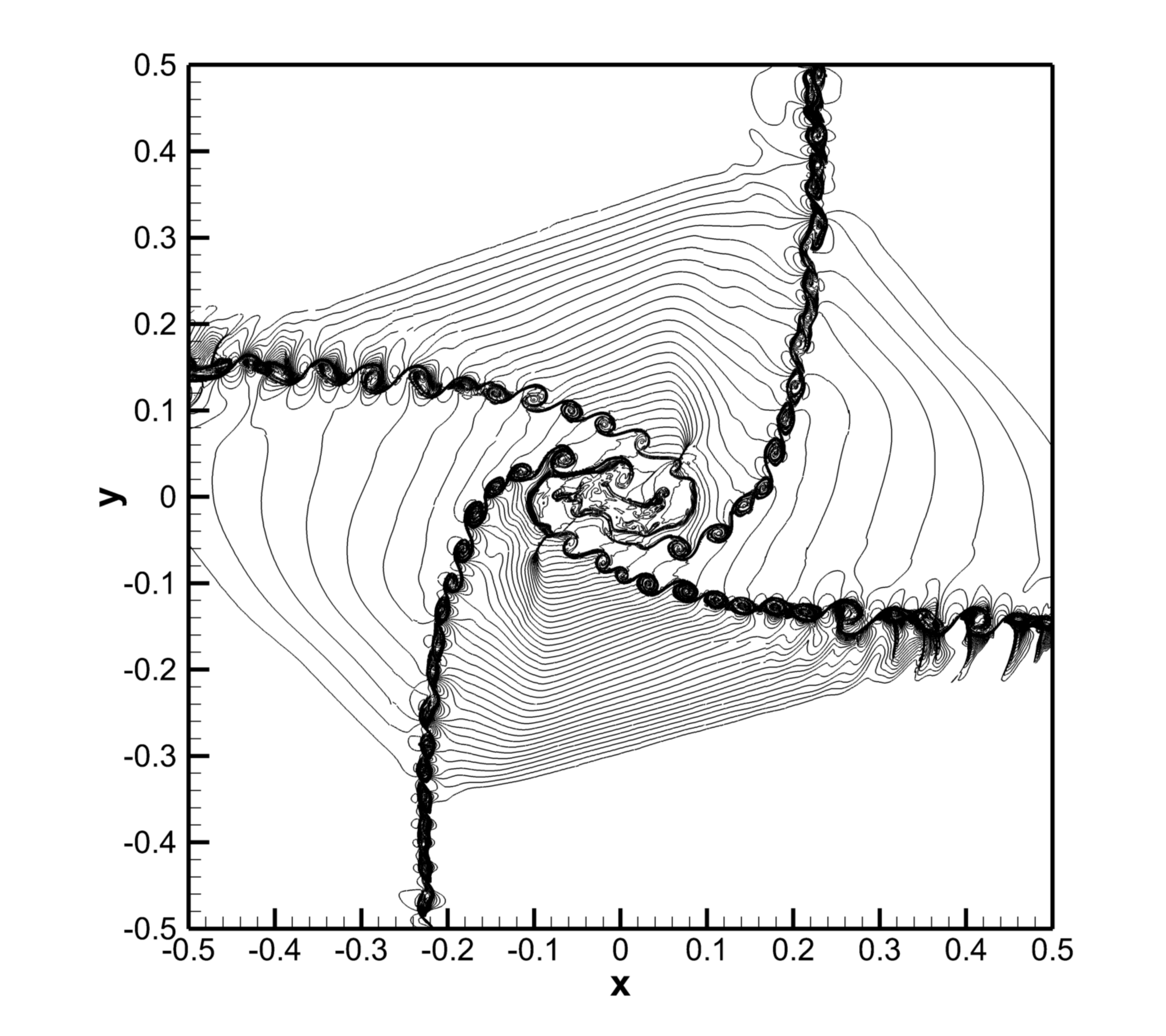}   
    \end{tabular}    
    \caption{ \label{fig:RP2DC6-compare}
      Comparison of \textbf{RP3} solved with two different numerical methods. {\em Left panel}: same  
			as Fig.\ref{fig:RP2D}, obtained with the ADER-DG-$\mathbb{P}_5$ method with subcell limiter on a coarse $100 \times 100$ main grid. 
		  {\em Right panel}: solution obtained with a sixth order ADER-WENO6 finite volume scheme with AMR 
		  (equivalent resolution: $1250 \times 1250$), showing the Kelvin-Helmholtz instability on the
		  shear waves.  
    }
  \end{center}
\end{figure}

To have a reliable comparison, we have repeated the test case \textbf{RP3} using a sixth order finite volume ADER-WENO 
scheme with space-time adaptive mesh refinement (AMR), using the method described in \cite{AMR3DCL}. The results of this comparison 
are illustrated in Fig. \ref{fig:RP2DC6-compare}: in the left panel we have simply reported again the solution of Fig. \ref{fig:RP2D}, 
while in the right panel we can see the solution with the sixth order ADER-WENO scheme using AMR. In the AMR simulations, the level 
0 grid is composed of $50\time 50$ elements which has then been adaptively refined using two levels of refinement and a refinement factor 
of 5 (right panel), which leads to an equivalent resolution of $1250 \times 1250$ elements on the finest level.  
The high order finite volume simulation with an AMR code clearly confirms the onset of a Kelvin-Helmholtz instability, which means that
the results obtained with the high order ADER-DG scheme with subcell limiter are reliable, since this kind of physical instability actually 
should be observed in numerical simulations for sufficiently fine meshes or sufficiently high order accurate schemes. 

For these classical test cases it seems that the coupling of \aposteriori subcell limiter and high order DG schemes is a valid 
tool to capture the discontinuous waves without spurious oscillations, and also the smooth part of the flow without excessive 
numerical dissipation.


\subsection{Shock-vortex interaction}  \label{ssec:SVortex}
\begin{figure}
  \begin{center}
    \includegraphics[width=0.85\textwidth]{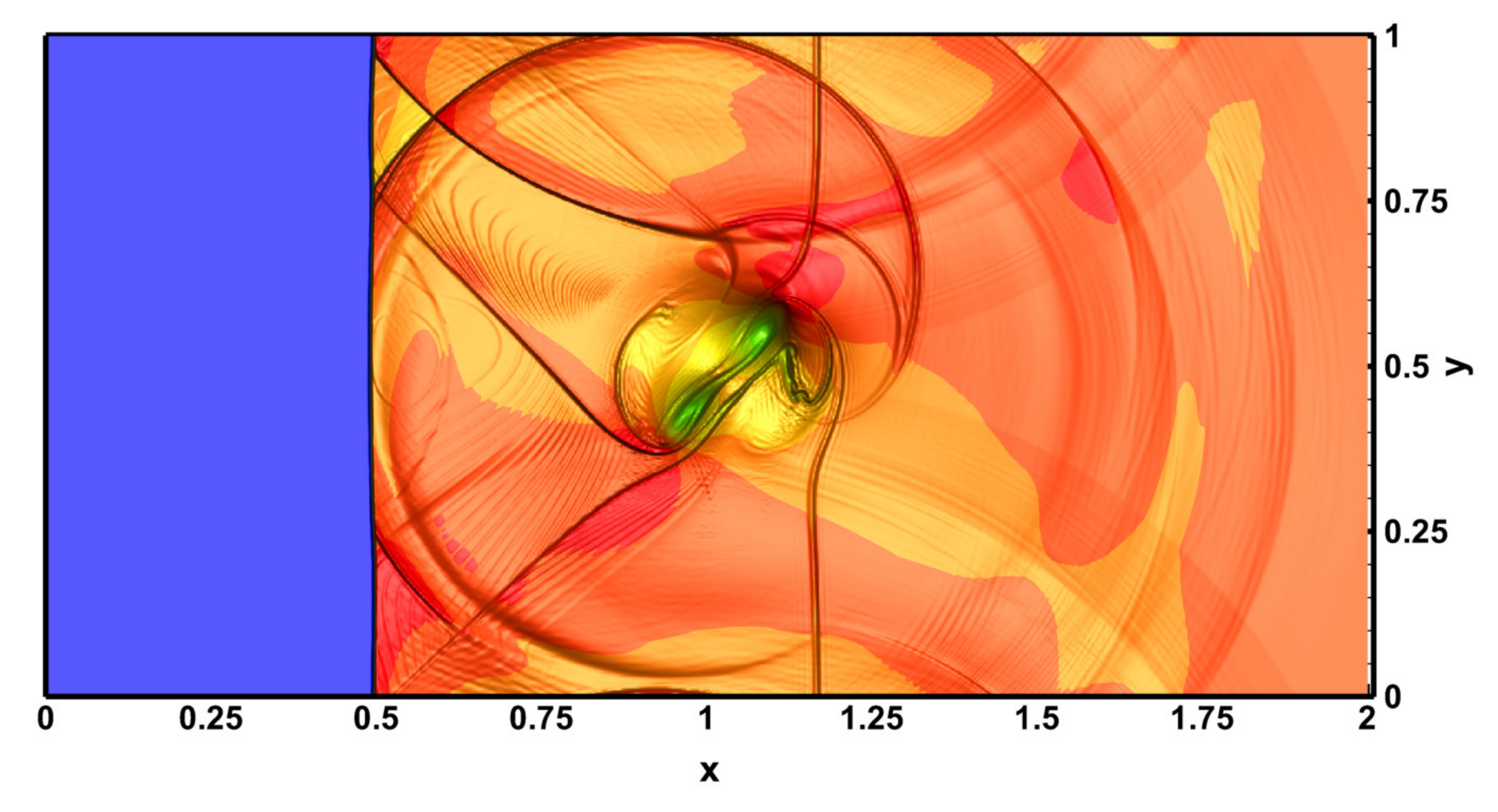}    \\
    \includegraphics[width=0.85\textwidth]{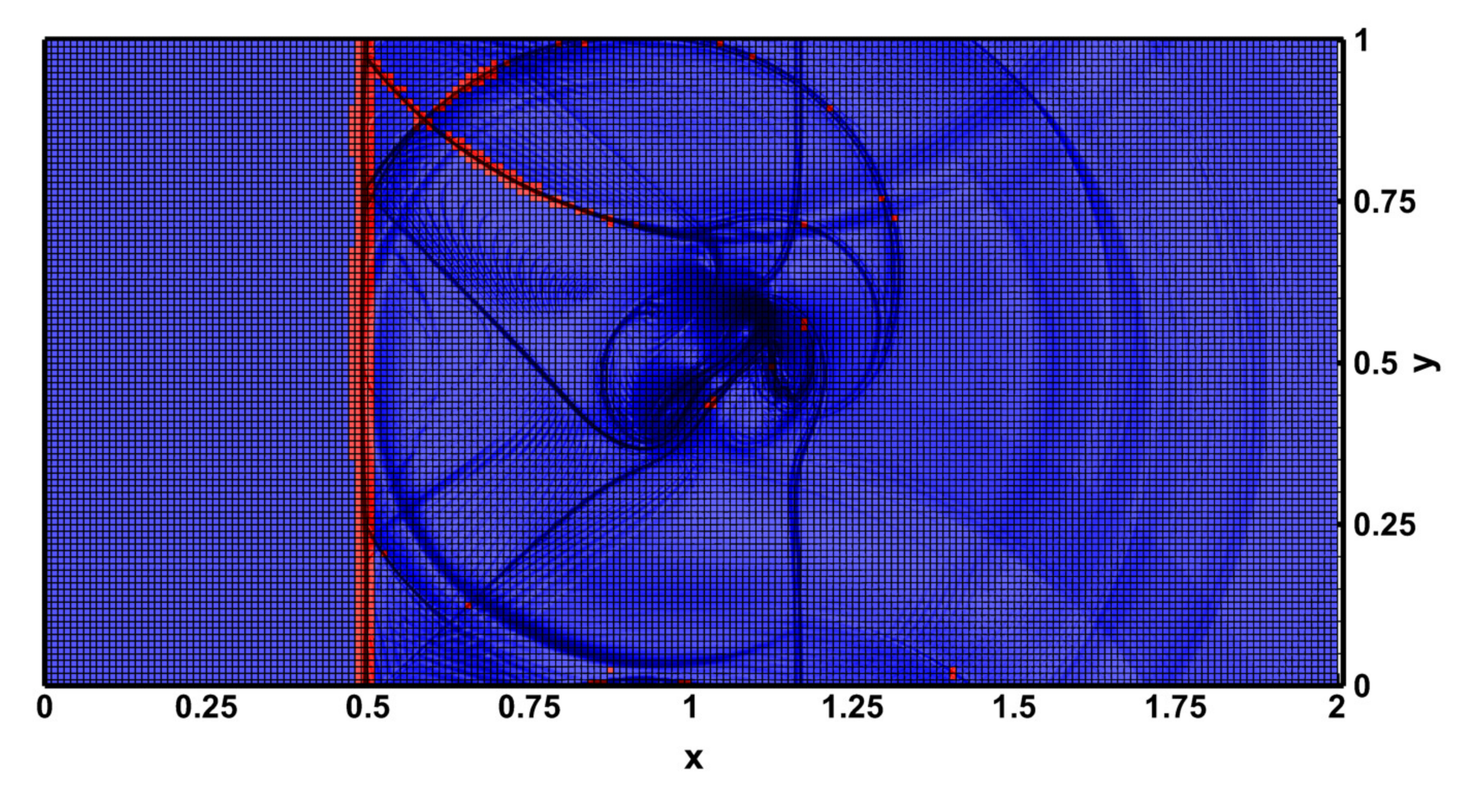}    				
    \caption{ \label{fig:SVI-P5-h1e-2}
      Shock-vortex interaction problem solved with ADER-DG-$\mathbb{P}_5$ supplemented 
      with \aposteriori ADER-WENO3 subcell limiter ---
      Top panel: distribution of the mass density 
      at time $t_{\text{final}}=0.7$ ---
      Bottom panel: troubled cells (red) and unlimited cells (blue).
    }
  \end{center}
\end{figure}

As a final test in two space dimensions, we have considered the interaction of a vortex with a steady shock wave. Originally proposed 
by \cite{donat}, this test is a true benchmark for a high order numerical scheme, as it involves a complex flow pattern with both 
smooth features and discontinuous waves. The initial conditions, defined over the computational domain $\Omega=[0; 2]\times[0,1]$, 
are given by a stationary normal shock wave placed at $x=0.5$ and by a vortex, which is placed at $(x_c,y_c)=(0.25,0.5)$. 
The shock Mach number is denoted by $M_{\rm S}$ and inside the vortex we have the following distribution of the angular velocity:  
\begin{equation}\label{eq:SVI}
  v_\phi = \left\{\begin{array}{ll}
v_m \frac{r}{a} & \;\textrm{for}\quad r \leq a\,, \\
 \noalign{\medskip}
 v_m \frac{a}{a^2-b^2}\left(r-\frac{b^2}{r}\right) & \;\textrm{for}\quad a\leq  r \leq b\,, \\
 \noalign{\medskip}
0 & {\rm otherwise} \,,\\
 \end{array}\right. 
\end{equation} 
with $r^2=(x-x_c)^2 + (y-y_c)^2$. 
The temperature of the vortex is obtained after solving the ordinary differential equation
\begin{equation}
\frac{dT}{dr}=\frac{\gamma-1}{R\gamma}\frac{v_\phi^2(r)}{r}\,,
\end{equation}
from which it is possible to compute the density and the pressure as
\begin{eqnarray}
p=p_0\left(\frac{T}{T_0}\right)^{\frac{\gamma}{\gamma-1}}, \quad \quad
\rho=\rho_0\left(\frac{T}{T_0}\right)^{\frac{1}{\gamma-1}}\,.
\end{eqnarray}
The unperturbed upstream values are of course related through the ideal gas equation of state, i.e. $p_0=R \rho_0 T_0$, where we set the gas constant to $R=1$.  
The strength of the vortex is described in terms of the Mach number $M_{\rm V}=v_m /c_0$, where $c_0=\sqrt{\gamma p_0/\rho_0}$ is the sound speed upstream of the shock. 
In our test, the specific values of these parameters are $\gamma=1.4$, $a=0.075$, $b=0.175$, $M_{\rm S}=1.5$, $M_{\rm V}=0.7$, $p_0=1$, $\rho_0=1$. Finally, the downstream values in the post-shock region
are computed through the classical Rankine--Hugoniot conditions~\cite{Landau-Lifshitz6}. We have solved this problem with the ADER-DG-$\mathbb{P}_5$ version of our scheme over a computational grid with 
characteristic length $h=1/100$. The results of our calculations are reported in Fig.~\ref{fig:SVI-P5-h1e-2}. The top panel shows the distribution of the mass density at time $t_{\text{final}}=0.7$, 
while in the bottom panel the troubled cells are shown in red and the unlimited cells are reported as usual in blue. 
Overall, these results confirm the ability of the scheme in capturing at the same time shock waves as well as smooth vortex features that produce small amplitude acoustic waves. 
Our results can be compared with those obtained through ADER finite volume schemes with fourth order of accuracy (see Fig. 9 in \cite{DumbserKaeser07}) and with the numerical 
solution provided in \cite{donat}. 

\subsection{3D explosion problem}\label{ssec:3Dexplosion}

To validate the new scheme in three spatial dimensions, we have considered an explosion problem on the computational 
domain $\Omega=[-1; 1]^3$. The setup represents a multi-dimensional extension of the classical Sod problem \citep{sod}, 
with initial conditions given by 
\begin{equation}\label{RP3D}
  \big(\rho, u, v, w, p\big) = \left\{\begin{array}{ll}
\big(1, 0, 0, 0,   1\big) & \;\textrm{for}\quad r \leq R\,, \\
 \noalign{\medskip}
 \big(0.125,   0, 0, 0,   0.1\big) & \;\textrm{for}\quad r>R    \,,  
 \end{array}\right.
\end{equation} 
where $r = \sqrt{x^2+y^2+z^2}$ is the radial coordinate, while $R=0.5$ denotes the radius of the initial discontinuity. 
The equation of state is assumed to be that of an ideal-gas, with adiabatic index $\gamma=1.4$. This test is important, 
as it involves the propagation of waves which are not aligned with the Cartesian grid. Since the problem is spherically 
symmetric, the reference solution can be obtained solving an equivalent one dimensional PDE in the radial direction with 
geometric source terms, see \cite{toro-book}. \\ 

The simulation has been performed using an ADER-DG-$\mathbb{P}_9$ scheme together with the new ADER-WENO subcell limiter. We have 
considered two different grid resolutions, with a number of elements given by $(25\times25\times 25)$  and 
$(100\times100\times 100)$, respectively. We emphasize that because of the high degree of the DG polynomial used here ($N=9$), 
which implies $(N+1)^4=10^4$ space-time degrees of freedom to represent the space-time predictor inside each element, the total 
amount of space-time degrees of freedom for the grid $100\times100\times100$ is actually $N_{\rm DOF}=10^{10}$. To the best of 
our knowledge, the explosion problem in three space dimensions has never been solved with a DG scheme on such a big mesh with 
such high order of approximation in space and time. The computation has been performed in parallel using the MPI standard on 
8000 CPU cores of the SuperMUC supercomputer at the \textit{Leibniz Rechenzentrum} (LRZ) in Munich, Germany. Our MPI 
parallelization of the scheme is based on domain decomposition, using the free software packages METIS and parMETIS, 
see \cite{metis}. 

Figure~\ref{fig:EP3D-1D-cuts} shows the profiles of the density along the $x-$ axis at time $t_{\text{final}}=0.2$, 
together with the 1D reference solution \cite{toro-book}. In these 1D cuts the discrete solution has been sampled on 
125 equidistant sample points in order to represent the subcell resolution capabilities of the DG method on the coarse $25^3$ 
grid. 
\begin{figure}[!htbp]
  \begin{center}
    \includegraphics[width=0.85\textwidth]{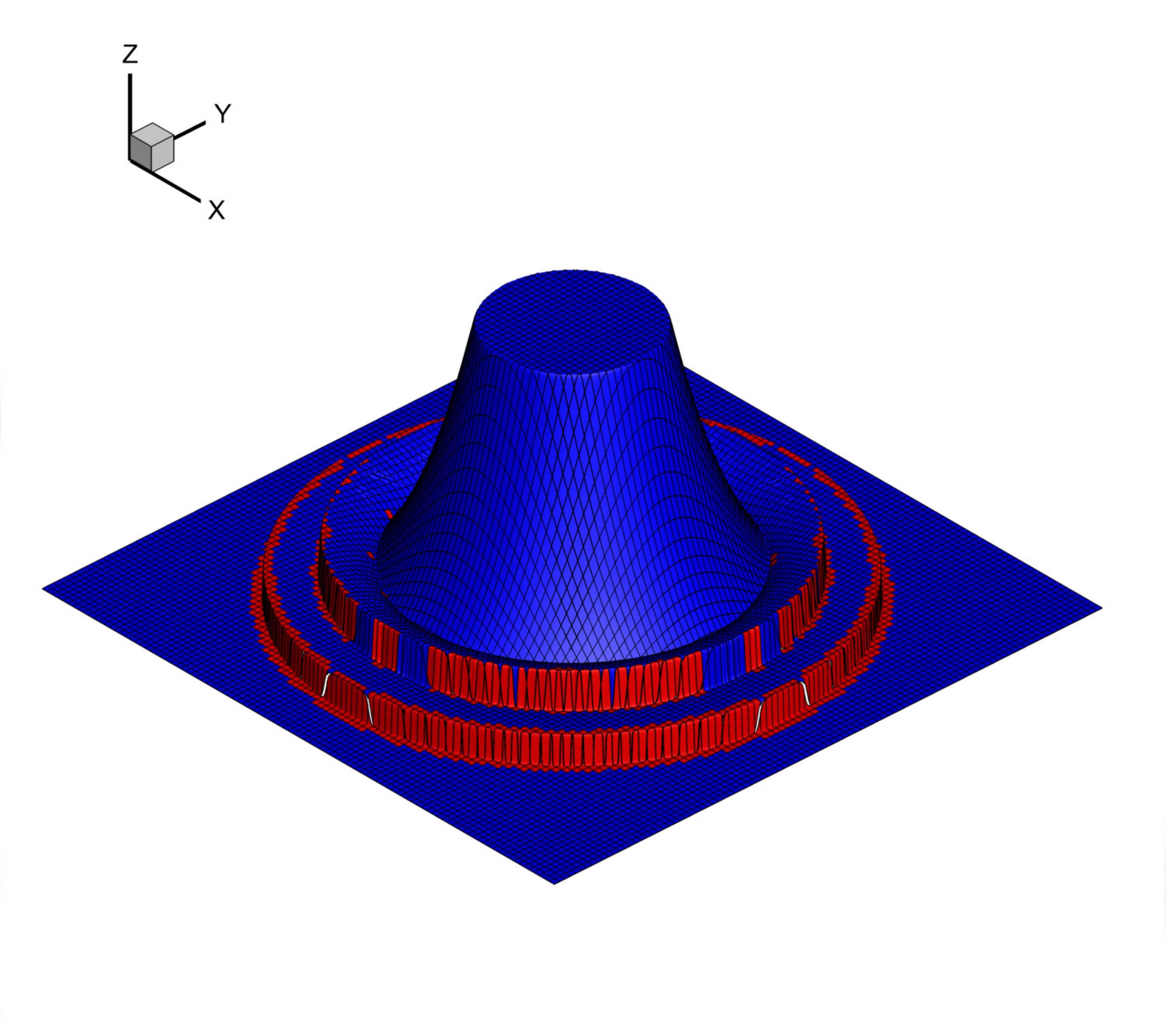}
					\vspace{-1cm} 
    \caption{ \label{fig:EP3D-limiter}
      Limited cells (red) updated with the subcell ADER-WENO3 finite volume scheme and unlimited DG cells (blue), together with the density 
			distribution on the plane $z=0$ for the three dimensional explosion problem at $t_{\text{final}}=0.2$ with ADER-DG-$\mathbb{P}_9$ 
			and $100^3$ elements, corresponding to 10 billion space-time degrees of freedom per time step.}
  \end{center}
\end{figure}
\begin{figure}[!htbp] 
  \begin{center} 
    \begin{tabular}{cc}
      \includegraphics[width=0.45\textwidth]{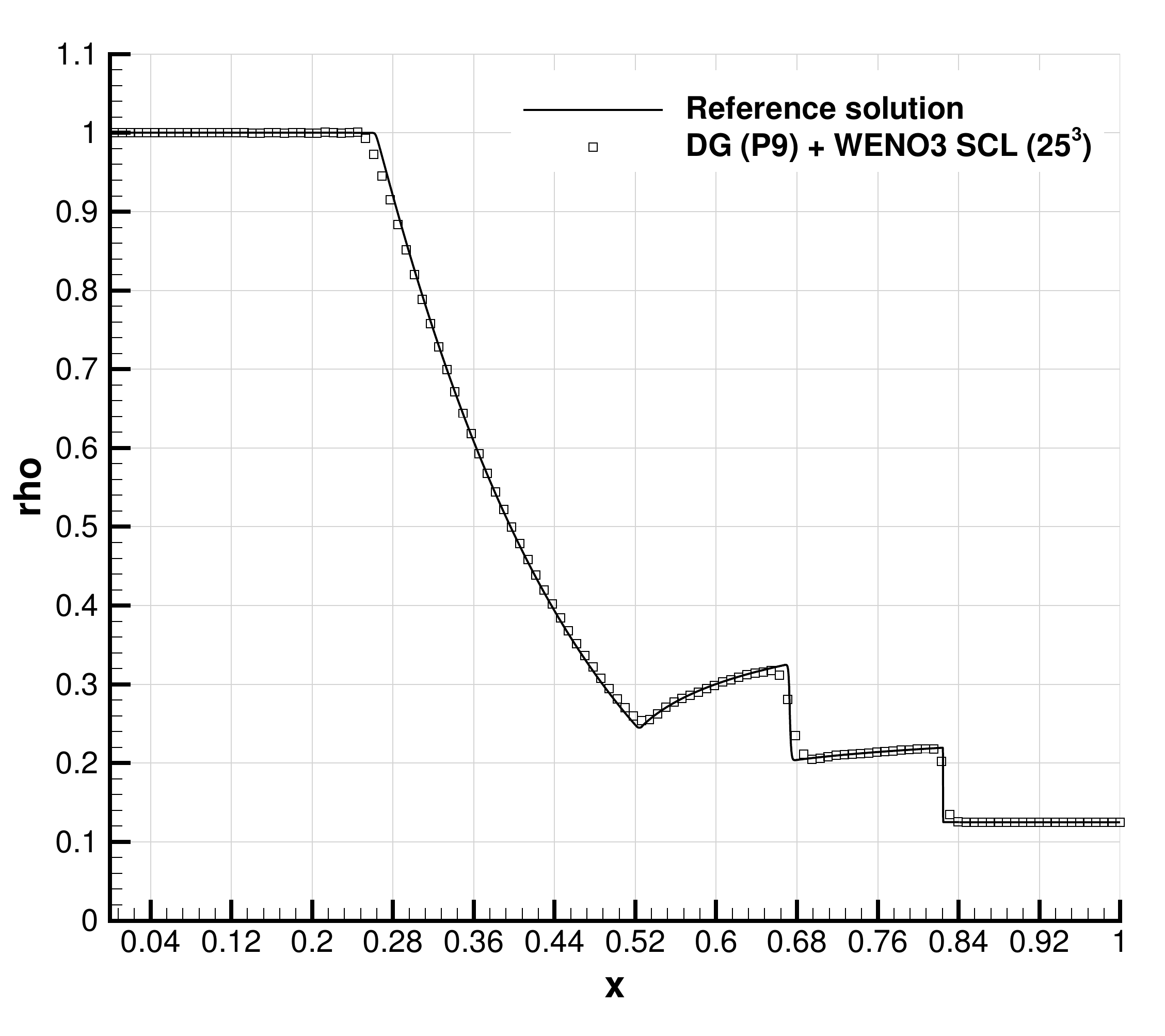}  &
      \includegraphics[width=0.45\textwidth]{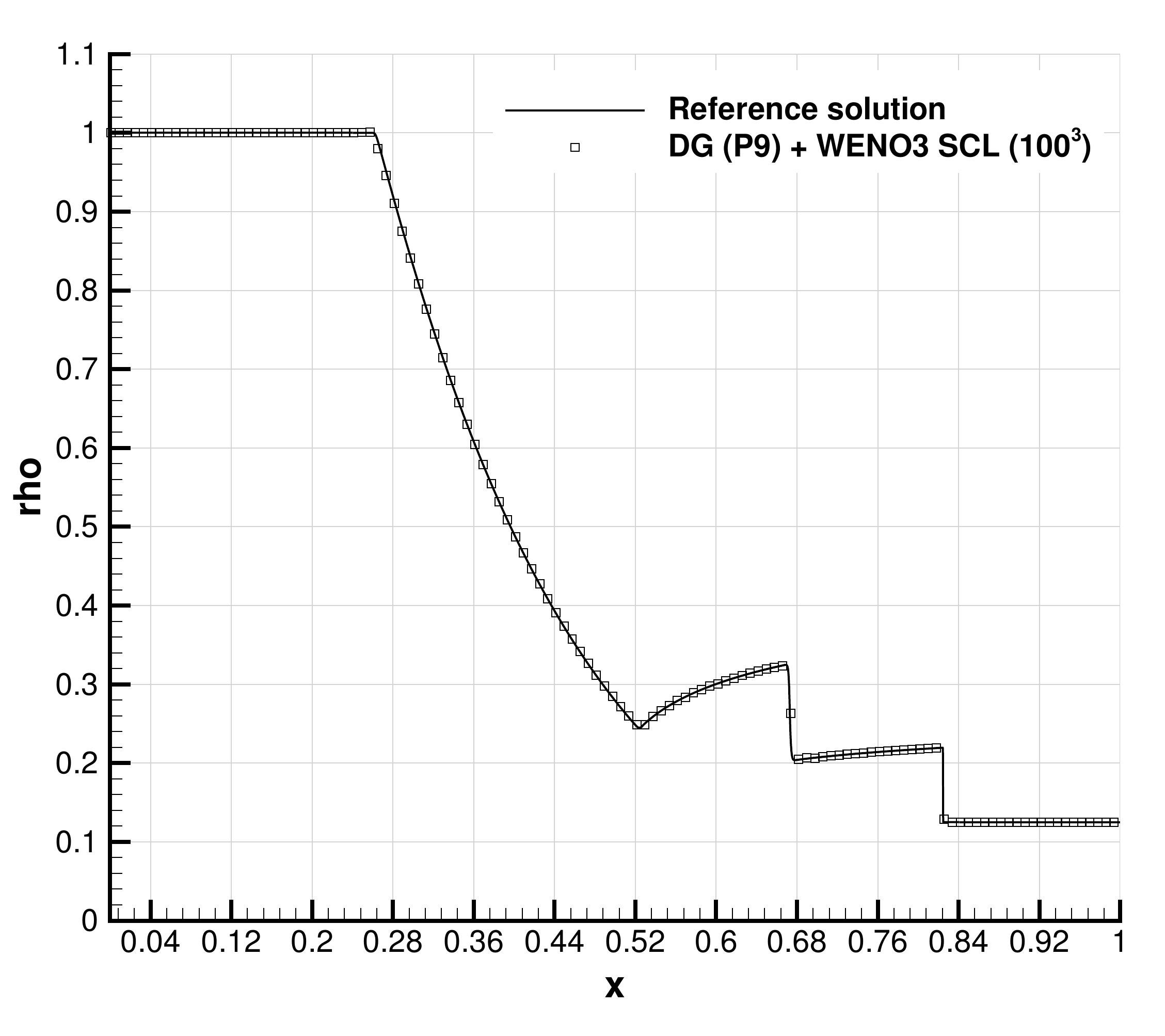}      
    \end{tabular} 
    \caption{ \label{fig:EP3D-1D-cuts}
      Three dimensional explosion problem at $t_{\text{final}}=0.2$ with ADER-DG-$\mathbb{P}_9$ supplemented 
      with \aposteriori ADER-WENO3 subcell limiter ---
      From left to right: comparison of the 1D reference solution with the numerical solution obtained 
      over two different grids with $25^3$ and $100^3$ elements, respectively.}
  \end{center}
\end{figure}
Even at the lowest resolution of $25^3$ elements, which means only $\sim 12$ elements for the interval $0\leq x \leq 1$, 
the high degree of the DG polynomial allows to capture the outgoing shock-wave and contact discontinuity with a very good 
level of accuracy. Note that the vertical gridlines displayed in the left panel of Fig. \ref{fig:EP3D-1D-cuts} correspond 
exactly to the cell size of the main grid, hence the shock wave as well as the contact discontinuity are each resolved in 
just \textit{one} single cell. 
The grid with $100^3$ elements provides an excellent result that is in perfect agreement with the 1D reference solution. 
In Figure \ref{fig:EP3D-limiter} we have reported a three-dimensional view of the mass density on the $z=0$ plane, 
while highlighting in red those elements which have been evolved with the ADER-WENO scheme on the subgrid. It is interesting 
to note that, even in this three dimensional problem, our \aposteriori limiter strategy has been required only for a rather
small number of elements that are crossed by the outward propagating waves, while no action of the limiter needs to been taken 
in the smooth regions of the solution. 
%
%



\section{Conclusion and Perspectives} \label{sec:conclusion}

In this paper we have presented a novel \aposteriori subcell limiter approach for the 
discontinuous Galerkin finite element method. The main building blocks of our scheme are 
rather simple. First, a so-called candidate solution is computed by using an \textit{unlimited} 
high order ADER-DG scheme. Second, we apply two simple \aposteriori detection criteria, 
namely the positivity of the candidate solution and a relaxed discrete maximum principle 
in the sense of piecewise polynomials. If the candidate solution satisfies both detection criteria
in a cell, then it is accepted, otherwise, the discrete solution is \textit{recomputed}. 
This is achieved by going back to the old time level $t^n$ in all troubled cells and by
projecting the discrete solution $\u_h(\x,t^n)$ onto a subgrid, obtaining and alternative
data representation $\v_h(\x,t^n)$ in terms of piecewise constant subcell averages. 
These subcell averages are then evolved in time using a robust and high order accurate 
ADER-WENO finite volume scheme on Cartesian grids, see \cite{titarevtoro,Balsara2013934,AMR3DCL}. 
For the troubled cells the new subcell averages $\v_h(\x,t^{n+1})$ are finally gathered back into 
the DG polynomials $\u_h(\x,t^{n+1})$ by using a subcell reconstruction operator and the entire 
process starts over again for the next time step. 

In this paper we have chosen the family of one-step ADER-DG schemes, but the \aposteriori
subcell limiter could be in principle also employed in combination with an explicit RK-DG scheme 
\cite{cbs4,CBS-convection-dominated}. Furthermore, also classical pointwise WENO schemes with 
TVD-Runge-Kutta time discretization as originally proposed by Jiang and Shu 
\cite{shu_efficient_weno} and their higher order extensions proposed by Balsara and Shu 
\cite{balsarashu} could be used to update the cell averages on the subgrid inside troubled
zones.   

\begin{figure}
  \begin{center} 
    \begin{tabular}{c} 
      \includegraphics[width=0.9\textwidth]{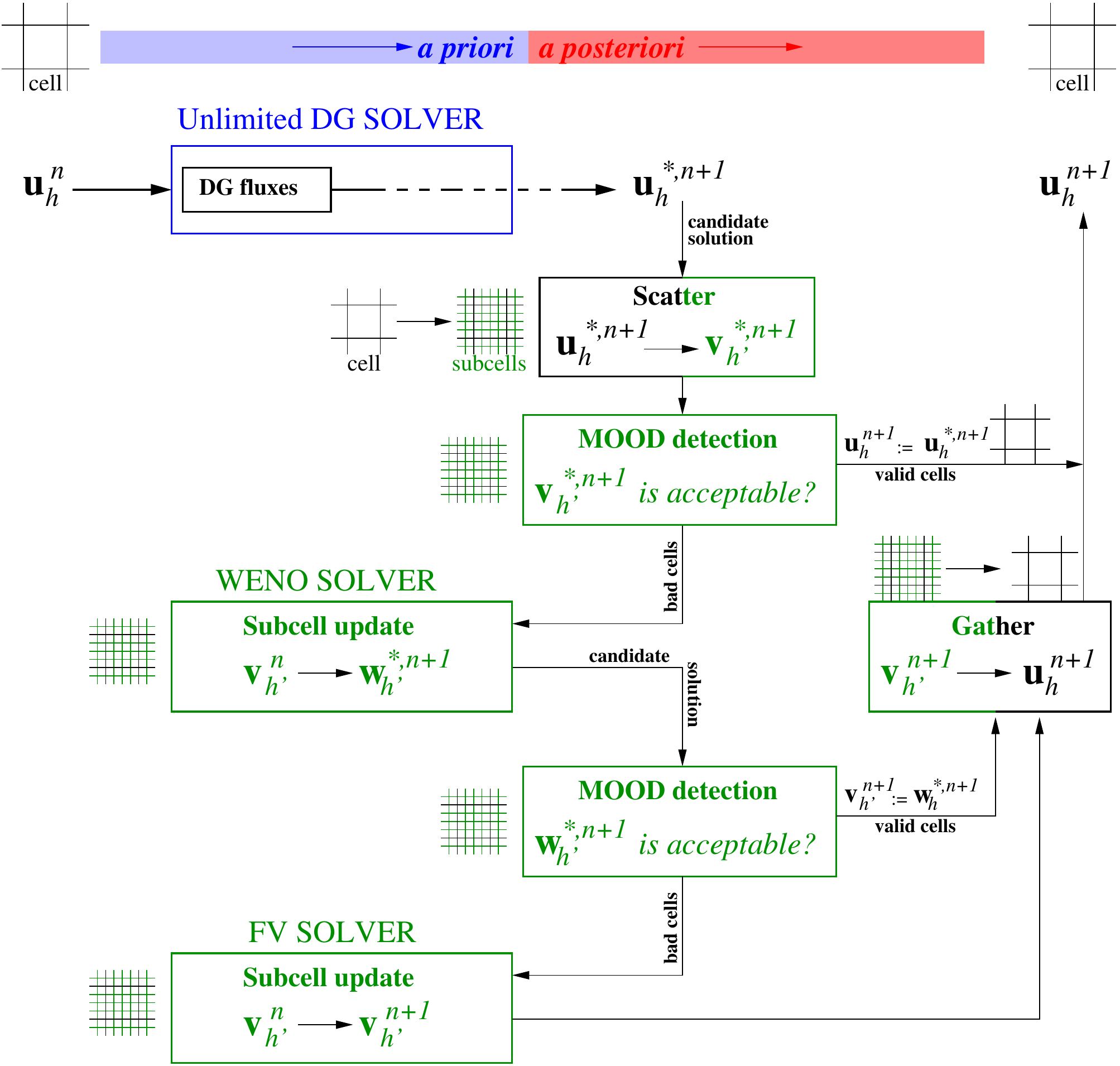} 
    \end{tabular}
    \caption{ \label{fig:algo2}
      Same figure as in Fig.~\ref{fig:algo} but with one more layer of an even lower order scheme (here 
      a standard first or second order TVD Finite Volume scheme) following the WENO procedure 
      if the latter fails to provide an acceptable solution according to the MOOD detection criteria. 
    }
  \end{center}
\end{figure}

Moreover, a possible future extension would be to add in the cascade after the DG-$\mathbb{P}_N$
scheme and the subcell WENO finite volume method, a second order TVD scheme or even a first order 
Godunov-type finite volume scheme on the subgrid, for the case where there are still remaining 
invalid subcells after the use of the subcell WENO procedure. 
The MOOD detection criteria would then be applied again after the subcell WENO solver on a 
WENO \textit{candidate solution}. In fact the WENO3 scheme may still have problems of positivity 
when the hyperbolic system of conservation laws becomes extremely complex, see \cite{ADER_MOOD_14} and
when the self-adjusting \textit{a priori} technique presented by Balsara in \cite{PositivityWENO2} 
is not used. 
A possibile solution to this problem either consists either in using the self-adjusting positivity 
preserving methodology presented in \cite{PositivityWENO2}, or to add one more MOOD layer in 
Fig.~\ref{fig:algo} with a more dissipative scheme to get for instance the sketch 
in Fig.~\ref{fig:algo2} 

%
%
In the future we plan to investigate the extension of the \aposteriori subcell limiter to ADER-DG schemes 
on unstructured meshes as well as to the more general family of $P_NP_M$ schemes proposed in \cite{Dumbser2008}. 
Finally, we would like to investigate the behavior of our new \aposteriori subcell limiter for DG schemes also 
on other systems of conservation laws and non-conservative hyperbolic PDE, such as for example 
compressible multi-material flows and classical or relativistic magneto-hydrodynamics (MHD).


\section*{Acknowledgments}
M.D. and O.Z. have been financed by the European Research Council (ERC) under the
European Union's Seventh Framework Programme (FP7/2007-2013) with the
research project \textit{STiMulUs}, ERC Grant agreement no. 278267.
R.L. has been partially funded by the ANR under the JCJC project
``ALE INC(ubator) 3D''.
This work has been authorized for publication under the reference LA-UR-14-24778. 

The authors would like to acknowledge PRACE for awarding access to the SuperMUC 
supercomputer based in Munich, Germany at the Leibniz Rechenzentrum (LRZ). 

Parts of the material contained in this work have been elaborated, gathered 
and tested during the SHARK-FV conference (\textit{Sharing Higher-order Advanced 
Know-how on Finite Volume} held in Ofir, Portugal in 2014, 
\texttt{www.math.univ-toulouse.fr/SHARK-FV/}.


\bibliographystyle{plain}
\bibliography{DG_MOOD}

\begin{thebibliography}{100}

\bibitem{Burbeau_2001}
A.Burbeau, P.Sagaut, and C.H. Bruneau.
\newblock A problem-independent limiter for high-order runge-kutta
  discontinuous galerkin methods.
\newblock {\em J. Comput. Phys.}, 169(1):111--150, May 2001.

\bibitem{Cook:2004}
A.W.Cook and W.H.Cabot.
\newblock A high-wavenumber viscosity for high-resolution numerical methods.
\newblock {\em J. Comput. Phys.}, 195(2):594--601, April 2004.

\bibitem{Mulet1}
A.~Baeza and P.~Mulet.
\newblock Adaptive mesh refinement techniques for high--order shock capturing
  schemes for multi--dimensional hydrodynamic simulations.
\newblock {\em International Journal for Numerical Methods in Fluids},
  52:455--471, 2006.

\bibitem{balsara2007}
D.~Balsara, C.~Altmann, C.D. Munz, and M.~Dumbser.
\newblock A sub-cell based indicator for troubled zones in {RKDG} schemes and a
  novel class of hybrid {RKDG}+{HWENO} schemes.
\newblock {\em Journal of Computational Physics}, 226:586--620, 2007.

\bibitem{balsarashu}
D.~Balsara and {C.W.} Shu.
\newblock Monotonicity preserving weighted essentially non-oscillatory schemes
  with increasingly high order of accuracy.
\newblock {\em Journal of Computational Physics}, 160:405--452, 2000.

\bibitem{balsarahlle2d}
D.S. Balsara.
\newblock {Multidimensional HLLE Riemann solver: Application to Euler and
  magnetohydrodynamic flows}.
\newblock {\em Journal of Computational Physics}, 229:1970--1993, 2010.

\bibitem{balsarahllc2d}
D.S. Balsara.
\newblock {A two-dimensional HLLC Riemann solver for conservation laws:
  Application to Euler and magnetohydrodynamic flows}.
\newblock {\em Journal of Computational Physics}, 231:7476--7503, 2012.

\bibitem{PositivityWENO2}
D.S. Balsara.
\newblock {Self--adjusting, positivity preserving high order schemes for
  hydrodynamics and magnetohydrodynamics}.
\newblock {\em Journal of Computational Physics}, 231:7504--7517, 2012.

\bibitem{BalsaraMultiDRS}
D.S. Balsara, M.~Dumbser, and R.~Abgrall.
\newblock {Multidimensional HLLC Riemann Solver for Unstructured Meshes - With
  Application to Euler and MHD Flows.}
\newblock {\em Journal of Computational Physics}, 261:172--208, 2014.

\bibitem{Balsara2013934}
D.S. Balsara, C.~Meyer, M.~Dumbser, H.~Du, and Z.~Xu.
\newblock Efficient implementation of {ADER} schemes for euler and
  magnetohydrodynamical flows on structured meshes – speed comparisons with
  runge–kutta methods.
\newblock {\em Journal of Computational Physics}, 235:934 -- 969, 2013.

\bibitem{BarthCharrier}
T.~Barth and P.~Charrier.
\newblock Energy stable flux formulas for the discontinuous {Galerkin}
  discretization of first-order nonlinear conservation laws.
\newblock Technical Report NAS-01-001, NASA, 2001.

\bibitem{BarthJespersen}
T.J. Barth and D.C. Jespersen.
\newblock The design and application of upwind schemes on unstructured meshes.
\newblock {\em AIAA Paper 89-0366}, pages 1--12, 1989.

\bibitem{caishu93}
W.~Cai and C.W. Shu.
\newblock Uniform high-order spectral methods for one and two dimen- sional
  euler equations.
\newblock {\em Journal of Computational Physics}, 104:427--443, 1993.

\bibitem{Cam01}
J.~C. Campbell and M.~J. Shaskhov.
\newblock A tensor artificial viscosity using a mimetic finite difference
  algorithm.
\newblock {\em J. Comput. Phys.}, 172(4):739--765, 2001.

\bibitem{Caramana-Burton-Shashkov-Whalen-98}
E.~J. Caramana, D.~E. Burton, M.~J. Shashkov, and P.~P. Whalen.
\newblock The construction of compatible hydrodynamics algorithms utilizing
  conservation of total energy.
\newblock {\em J. Comput. Phys.}, 146(1):227--262, 1998.

\bibitem{Csw98}
E.~J. Caramana, M.~J. Shashkov, and P.~P. Whalen.
\newblock Formulations of artificial viscosity for multidimensional shock wave
  computations.
\newblock {\em J. Comput. Phys.}, 144:70--97, 1998.

\bibitem{CasoniHuerta1}
E.~Casoni, J.~Peraire, and A.~Huerta.
\newblock One-dimensional shock-capturing for high-order discontinuous galerkin
  methods.
\newblock {\em International Journal for Numerical Methods in Fluids},
  71(6):737--755, 2013.

\bibitem{Casulli2009}
V.~Casulli.
\newblock A high-resolution wetting and drying algorithm for free-surface
  hydrodynamics.
\newblock {\em International Journal for Numerical Methods in Fluids},
  60:391--408, 2009.

\bibitem{CasulliStelling2011}
V.~Casulli and G.~S. Stelling.
\newblock Semi-implicit subgrid modelling of three-dimensional free-surface
  flows.
\newblock {\em International Journal for Numerical Methods in Fluids},
  67:441--449, 2011.

\bibitem{Feistauer4}
J.~Cesenek, M.~Feistauer, J.~Horacek, V.~Kucera, and J.~Prokopova.
\newblock Simulation of compressible viscous flow in time--dependent domains.
\newblock {\em Applied Mathematics and Computation}, 219:7139--7150, 2013.

\bibitem{CDL1}
S.~Clain, S.~Diot, and R.~Loub{\`e}re.
\newblock A high-order finite volume method for systems of conservation
  laws—multi-dimensional optimal order detection ({MOOD}).
\newblock {\em Journal of Computational Physics}, 230(10):4028 -- 4050, 2011.

\bibitem{cbs3}
B.~Cockburn, S.~Hou, and C.~W. Shu.
\newblock The {Runge}-{Kutta} local projection discontinuous {Galerkin} finite
  element method for conservation laws {IV}: the multidimensional case.
\newblock {\em Mathematics of Computation}, 54:545--581, 1990.

\bibitem{CBS-book}
B.~Cockburn, G.~E. Karniadakis, and {C.W.} Shu.
\newblock {\em Discontinuous {Galerkin} Methods}.
\newblock Lecture Notes in Computational Science and Engineering. Springer,
  2000.

\bibitem{cbs2}
B.~Cockburn, S.~Y. Lin, and C.W. Shu.
\newblock {TVB} {Runge}-{Kutta} local projection discontinuous {Galerkin}
  finite element method for conservation laws {III}: one dimensional systems.
\newblock {\em Journal of Computational Physics}, 84:90--113, 1989.

\bibitem{cbs1}
B.~Cockburn and C.~W. Shu.
\newblock {TVB} {Runge}-{Kutta} local projection discontinuous {Galerkin}
  finite element method for conservation laws {II}: general framework.
\newblock {\em Mathematics of Computation}, 52:411--435, 1989.

\bibitem{cbs0}
B.~Cockburn and C.~W. Shu.
\newblock The {Runge}-{Kutta} local projection {P1}-{Discontinuous} {Galerkin}
  finite element method for scalar conservation laws.
\newblock {\em Mathematical Modelling and Numerical Analysis}, 25:337--361,
  1991.

\bibitem{cbs4}
B.~Cockburn and C.~W. Shu.
\newblock The {Runge}-{Kutta} discontinuous {Galerkin} method for conservation
  laws {V}: multidimensional systems.
\newblock {\em Journal of Computational Physics}, 141:199--224, 1998.

\bibitem{CBS-convection-dominated}
B.~Cockburn and C.~W. Shu.
\newblock {Runge}-{Kutta} discontinuous {Galerkin} methods for
  convection-dominated problems.
\newblock {\em Journal of Scientific Computing}, 16:173--261, 2001.

\bibitem{CDL2}
S.~Diot, S.~Clain, and R.~Loub{\`e}re.
\newblock Improved detection criteria for the multi-dimensional optimal order
  detection ({MOOD}) on unstructured meshes with very high-order polynomials.
\newblock {\em Computers and Fluids}, 64:43 -- 63, 2012.

\bibitem{CDL3}
S.~Diot, R.~Loub{\`e}re, and S.~Clain.
\newblock The {MOOD} method in the three-dimensional case: Very-high-order
  finite volume method for hyperbolic systems.
\newblock {\em International Journal of Numerical Methods in Fluids},
  73:362--392, 2013.

\bibitem{Kuzmin2010}
D.Kuzmin.
\newblock A vertex-based hierarchical slope limiter for -adaptive discontinuous
  galerkin methods.
\newblock {\em Journal of Computational and Applied Mathematics}, 233(12):3077
  -- 3085, 2010.
\newblock Finite Element Methods in Engineering and Science (FEMTEC 2009).

\bibitem{Kuzmin2013}
D.Kuzmin.
\newblock Slope limiting for discontinuous galerkin approximations with a
  possibly non-orthogonal taylor basis.
\newblock {\em International Journal for Numerical Methods in Fluids},
  71(9):1178--1190, 2013.

\bibitem{Kuzmin2014}
D.Kuzmin.
\newblock Hierarchical slope limiting in explicit and implicit discontinuous
  galerkin methods.
\newblock {\em Journal of Computational Physics}, 257, Part B(0):1140 -- 1162,
  2014.
\newblock Physics-compatible numerical methods.

\bibitem{Feistauer6}
V.~Dolejsi, M.~Feistauer, and C.~Schwab.
\newblock On some aspects of the discontinuous galerkin finite element method
  for conservation laws.
\newblock {\em Mathematics and Computers in Simulation}, 61(3-6):333--346,
  2003.

\bibitem{Dumbser2008}
M.~Dumbser, D.~Balsara, E.F. Toro, and C.D. Munz.
\newblock A unified framework for the construction of one-step finite-volume
  and discontinuous {Galerkin} schemes.
\newblock {\em Journal of Computational Physics}, 227:8209--8253, 2008.

\bibitem{DumbserEnauxToro}
M.~Dumbser, C.~Enaux, and E.F. Toro.
\newblock Finite volume schemes of very high order of accuracy for stiff
  hyperbolic balance laws.
\newblock {\em Journal of Computational Physics}, 227:3971--4001, 2008.

\bibitem{DumbserKaeser06b}
M.~Dumbser and M.~K\"aser.
\newblock Arbitrary high order non-oscillatory finite volume schemes on
  unstructured meshes for linear hyperbolic systems.
\newblock {\em Journal of Computational Physics}, 221:693--723, 2007.

\bibitem{DumbserKaeser07}
M.~Dumbser, M.~K\"aser, V.A Titarev, and E.F. Toro.
\newblock Quadrature-free non-oscillatory finite volume schemes on unstructured
  meshes for nonlinear hyperbolic systems.
\newblock {\em Journal of Computational Physics}, 226:204--243, 2007.

\bibitem{dumbser_jsc}
M.~Dumbser and {C.D.} Munz.
\newblock Building blocks for arbitrary high order discontinuous {Galerkin}
  schemes.
\newblock {\em Journal of Scientific Computing}, 27:215--230, 2006.

\bibitem{OsherUniversal}
M.~Dumbser and E.~F. Toro.
\newblock On universal {Osher}--type schemes for general nonlinear hyperbolic
  conservation laws.
\newblock {\em Communications in Computational Physics}, 10:635--671, 2011.

\bibitem{DumbserZanotti}
M.~Dumbser and O.~Zanotti.
\newblock Very high order {PNPM} schemes on unstructured meshes for the
  resistive relativistic {MHD} equations.
\newblock {\em Journal of Computational Physics}, 228:6991--7006, 2009.

\bibitem{AMR3DCL}
M.~Dumbser, O.~Zanotti, A.~Hidalgo, and D.S. Balsara.
\newblock {ADER-WENO Finite Volume Schemes with Space-Time Adaptive Mesh
  Refinement}.
\newblock {\em Journal of Computational Physics}, 248:257--286, 2013.

\bibitem{Shock_vortex_95}
J.L. Ellzey, M.R. Henneke, J.M. Picone, and E.S.Oran.
\newblock The interaction of a shock with a vortex: Shock distortion and the
  production of acoustic waves.
\newblock {\em Phys. Fluids}, 7:172--184, 1995.

\bibitem{Feistauer7}
M.~Feistauer, V.~Dolejsi, and V.~Kucera.
\newblock On the discontinuous galerkin method for the simulation of
  compressible flow with wide range of mach numbers.
\newblock {\em Computing and Visualization in Science}, 10(1):17--27, 2007.

\bibitem{Feistauer5}
M.~Feistauer, V.~Kucera, and J.~Prokopov\'a.
\newblock Discontinuous galerkin solution of compressible flow in
  time--dependent domains.
\newblock {\em Mathematics and Computers in Simulation}, 80(8):1612--1623,
  2010.

\bibitem{GassnerDumbserMunz}
G.~Gassner, M.~Dumbser, F.~Hindenlang, and C.D. Munz.
\newblock Explicit one--step time discretizations for discontinuous {Galerkin}
  and finite volume schemes based on local predictors.
\newblock {\em Journal of Computational Physics}, 230:4232--4247, 2011.

\bibitem{GassnerKopriva}
G.~Gassner and D.A. Kopriva.
\newblock {A comparison of the dispersion and dissipation errors of Gauss and
  Gauss-Lobatto discontinuous Galerkin spectral element methods}.
\newblock {\em SIAM Journal on Scientific Computing}, 33(5):2560--2579, 2011.

\bibitem{Barter_2010}
G.E.Barter and D.L.Darmofal.
\newblock Shock capturing with pde-based artificial viscosity for dgfem: Part
  i. formulation.
\newblock {\em J. Comput. Phys.}, 229(5):1810--1827, March 2010.

\bibitem{godunov}
{S.K.} Godunov.
\newblock Finite difference methods for the computation of discontinuous
  solutions of the equations of fluid dynamics.
\newblock {\em Mathematics of the USSR - Sbornik}, 47:271--306, 1959.

\bibitem{eno}
A.~Harten, B.~Engquist, S.~Osher, and S.~Chakravarthy.
\newblock Uniformly high order essentially non-oscillatory schemes, {III}.
\newblock {\em Journal of Computational Physics}, 71:231--303, 1987.

\bibitem{HidalgoDumbser}
A.~Hidalgo and M.~Dumbser.
\newblock {ADER} schemes for nonlinear systems of stiff
  advection-diffusion-reaction equations.
\newblock {\em Journal of Scientific Computing}, 48:173--189, 2011.

\bibitem{Luo_2007}
H.Luo, J.D.Baum, and R.L\"{o}hner.
\newblock A hermite weno-based limiter for discontinuous galerkin method on
  unstructured grids.
\newblock {\em J. Comput. Phys.}, 225(1):686--713, July 2007.

\bibitem{HouLiu}
S.~Hou and X.~D. Liu.
\newblock Solutions of multi-dimensional hyperbolic systems of conservation
  laws by square entropy condition satisfying discontinuous {Galerkin} method.
\newblock {\em Journal of Scientific Computing}, 31:127--151, 2007.

\bibitem{CasoniHuerta2}
A.~Huerta, E.~Casoni, and J.~Peraire.
\newblock A simple shock-capturing technique for high-order discontinuous
  galerkin methods.
\newblock {\em International Journal for Numerical Methods in Fluids},
  69(10):1614--1632, 2012.

\bibitem{Zhu_hadap_2013}
H.Zhu and J.Qiu.
\newblock An h-adaptive rkdg method with troubled-cell indicator for
  two-dimensional hyperbolic conservation laws.
\newblock {\em Advances in Computational Mathematics}, 39(3-4):445--463, 2013.

\bibitem{GrothAMR}
L.~Ivan and C.P.T. Groth.
\newblock {High-order solution-adaptive central essentially non-oscillatory
  (CENO) method for viscous flows}.
\newblock {\em Journal of Computational Physics}, 257:830--862, 2014.

\bibitem{jiangshu}
G.~Jiang and C.W. Shu.
\newblock On a cell entropy inequality for discontinuous {Galerkin} methods.
\newblock {\em Mathematics of Computation}, 62:531--538, 1994.

\bibitem{shu_efficient_weno}
{G.-S.} Jiang and {C.W.} Shu.
\newblock Efficient implementation of weighted {ENO} schemes.
\newblock {\em Journal of Computational Physics}, 126:202--228, 1996.

\bibitem{Qiu_2004}
J.Qiu and C-W.Shu.
\newblock Hermite weno schemes and their application as limiters for
  runge-kutta discontinuous galerkin method: One-dimensional case.
\newblock {\em J. Comput. Phys.}, 193(1):115--135, January 2004.

\bibitem{Qiu_2005}
J.Qiu and C-W.Shu.
\newblock A comparison of troubled-cell indicators for runge--kutta
  discontinuous galerkin methods using weighted essentially nonoscillatory
  limiters.
\newblock {\em SIAM J. Sci. Comput.}, 27(3):995--1013, October 2005.

\bibitem{Zhu_2009}
J.Zhu and J.Qiu.
\newblock Hermite weno schemes and their application as limiters for
  runge-kutta discontinuous galerkin method, iii: Unstructured meshes.
\newblock {\em J. Sci. Comput.}, 39(2):293--321, May 2009.

\bibitem{Zhu_2008}
J.Zhu, J.~Qiu, C.-W.Shu, and M.Dumbser.
\newblock Runge-kutta discontinuous galerkin method using weno limiters ii:
  Unstructured meshes.
\newblock {\em J. Comput. Phys.}, 227(9):4330--4353, April 2008.

\bibitem{Zhu_13}
C.W.~Shu J.Zhu, X.Zhong and J.~Qiu.
\newblock Runge–kutta discontinuous galerkin method using a new type of weno
  limiters on unstructured meshes.
\newblock {\em J. Comp. Phys.}, 248:200--220, 2013.

\bibitem{kaeserjcp}
M.~K\"aser and A.~Iske.
\newblock {ADER} schemes on adaptive triangular meshes for scalar conservation
  laws.
\newblock {\em Journal of Computational Physics}, 205:486--508, 2005.

\bibitem{KlaijVanDerVegt}
C.~Klaij, J.J.W.~Van der Vegt, and H.~Van der Ven.
\newblock {Space-time discontinuous Galerkin method for the compressible
  Navier-Stokes equations}.
\newblock {\em Journal of Computational Physics}, 217:589--611, 2006.

\bibitem{Kopriva}
D.A. Kopriva.
\newblock Metric identities and the discontinuous spectral element method on
  curvilinear meshes.
\newblock {\em Journal of Scientific Computing}, 26(3):301--327, 2006.

\bibitem{KoprivaGassner}
D.A. Kopriva and G.~Gassner.
\newblock On the quadrature and weak form choices in collocation type
  discontinuous galerkin spectral element methods.
\newblock {\em Journal of Scientific Computing}, 44(2):136--155, 2010.

\bibitem{krivodonova2013analysis}
L.~Krivodonova and R.Qin.
\newblock An analysis of the spectrum of the discontinuous galerkin method.
\newblock {\em Applied Numerical Mathematics}, 64:1--18, 2013.

\bibitem{Krivodonova2004}
L.~Krivodonova, J.~Xin, J.-F. Remacle, N.~Chevaugeon, and J.E. Flaherty.
\newblock Shock detection and limiting with discontinuous galerkin methods for
  hyperbolic conservation laws.
\newblock {\em Applied Numerical Mathematics}, 48(3–4):323 -- 338, 2004.
\newblock Workshop on Innovative Time Integrators for \{PDEs\}.

\bibitem{kurganovtadmor}
A.~Kurganov and E.~Tadmor.
\newblock Solution of two-dimensional {Riemann} problems for gas dynamics
  without {Riemann} problem solvers.
\newblock {\em Numer. Methods Partial Differential Equations}, 18:584--608,
  2002.

\bibitem{Landau-Lifshitz6}
L.~D. Landau and E.~M. Lifshitz.
\newblock {\em Fluid Mechanics, Course of Theoretical Physics, Volume 6}.
\newblock Elsevier Butterworth-Heinemann, Oxford, 2004.

\bibitem{spectralfv3d}
Y.~Liu, M.~Vinokur, and Z.J. Wang.
\newblock {Spectral (Finite) Volume Method for Conservation Laws on
  Unstructured Grids V: Extension to Three-Dimensional Systems}.
\newblock {\em Journal of Computational Physics}, 212:454--472, 2006.

\bibitem{Loubere_HDR_2013}
R.~Loub{\`e}re.
\newblock {\em Contribution to Lagrangian and Arbitrary-Lagrangian-Eulerian
  numerical schemes}.
\newblock PhD thesis, University of Toulouse, France, 2013.
\newblock Habilitation {\`a} diriger des recherches.

\bibitem{ADER_MOOD_14}
R.~Loub{\`e}re, M.~Dumbser, and S.~Diot.
\newblock A new family of high order unstructured mood and ader finite volume
  schemes for multidimensional systems of hyperbolic conservation laws.
\newblock {\em Communication in Computational Physics}, 16:718--763, 2014.

\bibitem{Yang_AIAA_2009}
Z.J.~Wang M.~Yang.
\newblock A parameter-free generalized moment limiter for high-order methods on
  unstructured grids,.
\newblock {\em AIAA}, 605, 2009.

\bibitem{Michoski2011}
C.~Michoski, C.~Mirabito, C.~Dawson, D.~Wirasaet, E.J. Kubatko, and J.J.
  Westerink.
\newblock Adaptive hierarchic transformations for dynamically p-enriched
  slope-limiting over discontinuous galerkin systems of generalized equations.
\newblock {\em Journal of Computational Physics}, 230(22):8028 -- 8056, 2011.

\bibitem{vonneumann50}
J.~Von Neumann and R.~D. Richtmyer.
\newblock A method for the numerical calculation of hydrodynamic shocks.
\newblock {\em Journal of Applied Physics}, 21:232--237, 1950.

\bibitem{Nguyen11}
N.~C. Nguyen and J.~Peraire.
\newblock An adaptive shock-capturing hdg method for compressible flows
  presented at aiaa conference, 2011.

\bibitem{Persson_06}
P.-O. Persson and J.~Peraire.
\newblock Sub-cell shock capturing for discontinuous galerkin methods.
\newblock {\em AIAA Paper 2006-112}, 2006.

\bibitem{QiuDumbserShu}
J.~Qiu, M.~Dumbser, and {C.W.} Shu.
\newblock The discontinuous {Galerkin} method with {Lax}-{Wendroff} type time
  discretizations.
\newblock {\em Computer Methods in Applied Mechanics and Engineering},
  194:4528--4543, 2005.

\bibitem{QiuShu2}
J.~Qiu and {C.W.} Shu.
\newblock Hermite {WENO} schemes and their application as limiters for
  {Runge}-{Kutta} discontinuous {Galerkin} method: one-dimensional case.
\newblock {\em Journal of Computational Physics}, 193:115--135, 2003.

\bibitem{QiuShu3}
J.~Qiu and {C.W.} Shu.
\newblock Hermite {WENO} schemes and their application as limiters for
  {Runge}-{Kutta} discontinuous {Galerkin} method {II}: two dimensional case.
\newblock {\em Computers and Fluids}, 34:642--663, 2005.

\bibitem{QiuShu1}
J.~Qiu and {C.W.} Shu.
\newblock {Runge}-{Kutta} discontinuous {Galerkin} method using {WENO}
  limiters.
\newblock {\em SIAM Journal on Scientific Computing}, 26:907--929, 2005.

\bibitem{donat}
A.~Rault, G.~Chiavassa, and R.~Donat.
\newblock Shock-vortex interactions at high mach numbers.
\newblock {\em Journal of Scientific Computing}, 19:347--371, 2003.

\bibitem{Biswas_94}
R.Biswas, K.~D. Devine, and J.~E. Flaherty.
\newblock Parallel, adaptive finite element methods for conservation laws.
\newblock {\em APPL. NUMER. MATH}, 14:255--283, 1994.

\bibitem{reed}
{W.H.} Reed and {T.R.} Hill.
\newblock Triangular mesh methods for neutron transport equation.
\newblock Technical Report LA-UR-73-479, Los Alamos Scientific Laboratory,
  1973.

\bibitem{Hartman_02}
R.Hartmann and P.Houston.
\newblock Adaptive discontinuous {G}alerkin finite element methods for the
  compressible {E}uler equations.
\newblock {\em J. Comp. Phys.}, 183(2):508--532, 2002.

\bibitem{Rider_2001}
W.J. Rider and L.~G. Margolin.
\newblock Simple modifications of monotonicity-preserving limiter.
\newblock {\em Journal of Computational Physics}, 174(1):473 -- 488, 2001.

\bibitem{riemann1}
B.~Riemann.
\newblock {\"Uber die Fortpflanzung ebener Luftwellen von endlicher
  Schwingungsweite}.
\newblock {\em G\"ottinger Nachrichten}, 19, 1859.

\bibitem{riemann2}
B.~Riemann.
\newblock {\"Uber die Fortpflanzung ebener Luftwellen von endlicher
  Schwingungsweite}.
\newblock {\em Abhandlungen der K\"oniglichen Gesellschaft der Wissenschaften
  zu G\"ottingen}, 8:43--65, 1860.

\bibitem{Rusanov:1961a}
V.~V. Rusanov.
\newblock {Calculation of Interaction of Non--Steady Shock Waves with
  Obstacles}.
\newblock {\em J. Comput. Math. Phys. USSR}, 1:267--279, 1961.

\bibitem{schulzrinne}
C.~W. Schulz-Rinne.
\newblock Classification of the {Riemann} problem for two-dimensional gas
  dynamics.
\newblock {\em SIAM J. Math. Anal.}, 24:76--88, 1993.

\bibitem{schwartzkopff}
T.~Schwartzkopff, {C.D.} Munz, and {E.F.} Toro.
\newblock {ADER}: A high order approach for linear hyperbolic systems in 2d.
\newblock {\em Journal of Scientific Computing}, 17(1-4):231--240, 2002.

\bibitem{Shu1}
{C.W.} Shu.
\newblock Essentially non-oscillatory and weighted essentially non-oscillatory
  schemes for hyperbolic {Conservation} {Laws}.
\newblock {\em NASA/CR-97-206253 ICASE Report No.97-65}, November 1997.

\bibitem{shuosher2}
C.W. Shu and S.~Osher.
\newblock Efficient implementation of essentially non-oscillatory shock
  capturing schemes {II}.
\newblock {\em Journal of Computational Physics}, 83:32--78, 1989.

\bibitem{sod}
G.~A. Sod.
\newblock A survey of several finite difference methods for systems of
  non-linear hyperbolic conservation laws.
\newblock {\em Journal of Computational Physics}, 27:1--31, 1978.

\bibitem{Sonntag}
M.~Sonntag and C.D. Munz.
\newblock Shock capturing for discontinuous galerkin methods using finite
  volume subcells.
\newblock In J.~Fuhrmann, M.~Ohlberger, and C.~Rohde, editors, {\em Finite
  Volumes for Complex Applications VII}, pages 945--953. Springer, 2014.

\bibitem{stroud}
{A.H.} Stroud.
\newblock {\em Approximate Calculation of Multiple Integrals}.
\newblock Prentice-Hall Inc., Englewood Cliffs, New Jersey, 1971.

\bibitem{spectralfv.dg}
Y.~Sun and Z.J. Wang.
\newblock Evaluation of discontinuous galerkin and spectral volume methods for
  scalar and system conservation laws on unstructured grids.
\newblock {\em International Journal for Numerical Methods in Fluids},
  45(8):819--838, 2004.

\bibitem{SureshHuynh}
A.~Suresh and H.T. Huynh.
\newblock Accurate monotonicity-preserving schemes with runge-kutta time
  stepping.
\newblock {\em Journal of Computational Physics}, 136:83--99, 1997.

\bibitem{taube_jsc}
A.~Taube, M.~Dumbser, D.~Balsara, and {C.D.} Munz.
\newblock Arbitrary high order discontinuous {Galerkin} schemes for the
  magnetohydrodynamic equations.
\newblock {\em Journal of Scientific Computing}, 30:441--464, 2007.

\bibitem{toro3}
{V.A.} Titarev and {E.F.} Toro.
\newblock {ADER}: Arbitrary high order {Godunov} approach.
\newblock {\em Journal of Scientific Computing}, 17(1-4):609--618, December
  2002.

\bibitem{titarevtoro}
V.A. Titarev and E.F. Toro.
\newblock {ADER} schemes for three-dimensional nonlinear hyperbolic systems.
\newblock {\em Journal of Computational Physics}, 204:715--736, 2005.

\bibitem{toro-book}
E.F. Toro.
\newblock {\em {Riemann} Solvers and Numerical Methods for Fluid Dynamics}.
\newblock Springer, second edition, 1999.

\bibitem{toro4}
{E.F.} Toro and {V. A.} Titarev.
\newblock Solution of the generalized {Riemann} problem for advection-reaction
  equations.
\newblock {\em Proc. Roy. Soc. London}, pages 271--281, 2002.

\bibitem{spacetimedg1}
J.~J.~W. van~der Vegt and H.~van~der Ven.
\newblock Space--time discontinuous {Galerkin} finite element method with
  dynamic grid motion for inviscid compressible flows {I}. general formulation.
\newblock {\em Journal of Computational Physics}, 182:546--585, 2002.

\bibitem{spacetimedg2}
H.~van~der Ven and J.~J.~W. van~der Vegt.
\newblock Space--time discontinuous {Galerkin} finite element method with
  dynamic grid motion for inviscid compressible flows {II}. efficient flux
  quadrature.
\newblock {\em Comput. Methods Appl. Mech. Engrg.}, 191:4747--4780, 2002.

\bibitem{neu}
J.~von Neumann and R.~D. Richtmyer.
\newblock A method for the numerical calculations of hydrodynamical shocks.
\newblock {\em J. Appl. Phys.}, 21:232--238, 1950.

\bibitem{wang2011adaptive}
Z.J. Wang.
\newblock {\em Adaptive High-order Methods in Computational Fluid Dynamics}.
\newblock Advances in computational fluid dynamics. World Scientific, 2011.

\bibitem{spectralfv.bnd}
Z.J. Wang and Y.~Liu.
\newblock Extension of the spectral volume method to high-order boundary
  representation.
\newblock {\em Journal of Computational Physics}, 211:154--178, 2006.

\bibitem{spectralfv2d}
Z.J. Wang, L.~Zhang, and Y.~Liu.
\newblock Spectral (finite) volume method for conservation laws on unstructured
  grids iv: Extension to two-dimensional euler equations.
\newblock {\em Journal of Computational Physics}, 194:716--741, 2004.

\bibitem{woodwardcol84}
P.~Woodward and P.~Colella.
\newblock The numerical simulation of two-dimensional fluid flow with strong
  shocks.
\newblock {\em Journal of Computational Physics}, 54:115--173, 1984.

\bibitem{ShuPositivity2}
Y.~Xing, X.~Zhang, and C.W. Shu.
\newblock {Positivity preserving high order well balanced discontinuous
  Galerkin methods for the shallow water equations}.
\newblock {\em Advances in Water Resources}, 33:1476--1493, 2010.

\bibitem{Zhong_2013}
X.Zhong and C-W.Shu.
\newblock A simple weighted essentially nonoscillatory limiter for runge-kutta
  discontinuous galerkin methods.
\newblock {\em J. Comput. Phys.}, 232(1):397--415, January 2013.

\bibitem{Yang_parameterfree_09}
M.~Yang and Z.~Wang.
\newblock {A Parameter-Free Generalized Moment Limiter for High-Order Methods
  on Unstructured Grids}.
\newblock {\em Advances in Applied Mathematics and Mechanics}, 2009.

\bibitem{ShuPositivity1}
X.~Zhang and C.W. Shu.
\newblock {On positivity preserving high order discontinuous Galerkin schemes
  for compressible Euler equations on rectangular meshes}.
\newblock {\em Journal of Computational Physics}, 229:8918--8934, 2010.

\bibitem{ShuPositivity3}
X.~Zhang and C.W. Shu.
\newblock {Maximum-principle-satisfying and positivity-preserving high order
  schemes for conservation laws: Survey and new developments}.
\newblock {\em Proceedings of the Royal Society A}, 467:2752--2776, 2011.

\bibitem{Zhu_3D_12}
J.~Zhu and J.~Qiu.
\newblock Runge-kutta discontinuous galerkin method using weno type limiters:
  Three dimensional unstructured meshes.
\newblock {\em Commun.Comput. Phys.}, pages 985--1005, 2012.

\end{thebibliography}


\end{document}